\documentclass[a4paper,11pt]{article}

\usepackage{amsfonts,amssymb,amsbsy,amsmath,amsthm,mathrsfs,enumerate,verbatim}

\usepackage{graphics} 
\usepackage{epsfig} 
\usepackage{graphicx}
\usepackage{algorithmic,algorithm}
\usepackage{epstopdf}
\usepackage[colorlinks=true]{hyperref}
\hypersetup{urlcolor=blue, citecolor=red}

\usepackage{amsfonts}
\usepackage{amsmath}
\usepackage{amssymb}
\usepackage{graphicx}
\usepackage{epsf}
\usepackage{epsfig}

\usepackage{color}
\usepackage{afterpage}
\usepackage{verbatim}

\usepackage{times}

\newtheorem{theo}{Theorem}

\newcommand{\bel}{\begin{equation} \label}
\newcommand{\ee}{\end{equation}}

\def\beq{\begin{equation}}
\def\eeq{\end{equation}}
\newcommand{\bea}{\begin{eqnarray}}
\newcommand{\eea}{\end{eqnarray}}
\newcommand{\beas}{\begin{eqnarray*}}
\newcommand{\eeas}{\end{eqnarray*}}

{

 \usepackage{graphicx,color}
 \definecolor{mygreen}{cmyk}{1,0,1,0.1}

\usepackage{epsf}
\usepackage{epstopdf}

\usepackage{pdfsync}

\graphicspath{
   {FIGURES/}
}

\begin{document}


\title{A time-adaptive optimization approach for reconstructing  immune response  in a mathematical model \\ of  acute HIV infection using clinical data}

\author{L. Beilina  \thanks{
Department of Mathematical Sciences, Chalmers University of Technology and
 University of Gothenburg, SE-42196 Gothenburg, Sweden, e-mail: \texttt{\
larisa@chalmers.se}}
  \and
I. Gainova \thanks{Sobolev Institute of Mathematics SB RAS, 630090 Novosibirsk, Russia, e-mail: \texttt{\
gajnova@math.nsc.ru}}
  \and
G. Bocharov \thanks{Marchuk Institute of Numerical Mathematics RAS, 119333
Moscow, Russia, e-mail: \texttt{\
gbocharov@gmail.com}}
}

\date{}

\maketitle

\begin{abstract}

  The paper proposes a time-adaptive optimization approach for determining the
time-dependent immune response function in a mathematical model of
acute HIV infection, using clinical data from four untreated patients.
We formulate the problem as a parameter identification problem
for  an immune response system of ODE which includes
 novel component integrated into the third equation of the classical three-equation HIV model.
 Tikhonov's regularization method, Lagrangian approach, from which we derive the
optimality conditions, and a numerical scheme to solve the forward
 and adjoint problems, as well as parameter identification
problem, are presented.

 Three different a posteriori error
estimates are derived and
based on these estimates a time adaptive
 optimization algorithm
is formulated.
Numerical experiments demonstrate the effectiveness of the proposed
adaptive method in reconstructing the immune response function during the acute
phase of HIV infection, using patient-specific clinical data.
 Computational results show improvement of
  reconstruction of
 immune response function
using the local
  time-adaptive mesh refinement method compared to the standard conjugate gradient method
  applied on a uniform time mesh.

\end{abstract}

\maketitle


\section{Introduction}

\label{sec:0}

The mathematical modeling of immune dynamics is a foundational
component of immunological research \cite{Boch2012, Boch2018,
Eftimie}.  Although advanced experimental methods have significantly
enhanced the analysis of immune function, mathematical modeling
continues to play a crucial role in clinical applications—especially
in tailoring individualized treatment strategies for pathological
processes  such that bacterial/viral infections or tumor growth.
It is well established that physiological parameters differ across
individuals. Consequently, developing robust and efficient parameter
estimation methods is essential for integrating patient-specific data
into mathematical models, thereby enabling truly personalized
treatment strategies \cite{Banks}.

Parameter identification problem (PIP) for mathematical models
frequently involves solution of nonlinear and ill-posed problems, making them
particularly challenging to address through numerical methods
\cite{Tikhonov, T1, T2}.
The computational algorithms developed in this work are based on
adaptive time-mesh refinement techniques  for PIP for ODE or PDE -- see recent advances in this
area in \cite{BErG,BookBK,Sprg-1,Sprg-2}.

The primary objective of our work is to develop a time-adaptive
algorithm for determining the immune response function within a
mathematical model of untreated HIV infection. The model comprises a
standard system of three ordinary differential equations (ODEs),
enhanced by a novel component integrated into the third equation.
In  the recent work \cite{BErG}
the adaptive finite element method has been shown to significantly
enhance parameter reconstruction when applied to assessing drug
efficacy in the mathematical model of HIV infection  proposed in \cite{Sriv09}.

In comparison to other optimal control algorithms used for solving
parameter identification problems (PIPs), notable examples can be
found in the works of \cite{alb, Arruda3, Hatz} and references
therein, our time-mesh refinement algorithm
 is based on rigorous finite
element  analysis
for a posteriori error in the reconstructed parameter.

The approach proposed in this work is applicable to a wide range of
parameter identification problems (PIPs), including more complex
models of HIV infection that incorporate additional unknown functions
and parameters  --see, for example,  HIV models  proposed  in \cite{Arruda1, Arruda2, Shu,
Nowak, Kepler, Smith, Oui, Yan}.  However, these models are
substantially more complex than the three-equation system studied here
and may be more appropriately addressed in future research.

The acute phase of HIV infection is associated with an exponential
increase in viral load and, in most HIV infected individuals, begins
with fever, headache, increased heart rate, and symptoms of viral
spread to the lymphoid tissues (lymphadenopathy). The symptoms occur
both before and at the peak of viremia and are relatively short in
duration. The peak of viremia is more than $10^{6}$ copies of RNA/ml
of blood. After the peak, the viral load begins to decrease, both due
to the body's immune response and due to the limited population of
target cells (macrophages and CD4+ T-lymphocytes), and reaches a point
of stabilization of the infectious process (so called as "viral
equilibrium point"). The infection enters a chronic, latent phase,
which can last several years.

Events occurring during the acute phase as well as quantitative characteristics, such as the level
viral load at the point of stabilization of the infectious process and the level
Target T cells at the end of the acute phase of HIV infection, are
one of the factors determining the long-term development forecast
diseases. In addition, interest in this research is due to the fact that
the risk of HIV transmission is significantly higher during the acute phase.
In the context of HIV infection, mathematical models play a crucial role in understanding disease dynamics, predicting outcomes, and informing treatment strategies. Evaluation of the parameters which characterized the model is essential for fitting the one to clinical data and improving our understanding of HIV infection pathogenesis.

Our mathematical model is based on consideration of a standard three-component model of the acute HIV infection dynamics, which describes, in terms of differential equations, the interaction of uninfected target cells, infected cells and viral particles. This model allows us to present a qualitative picture of the acute phase: a sharp increase in viral load, its peak and further decrease until the viral load is established at an
equilibrium, and then calculate the basic rate of virus reproduction ($R_0$). The magnitude of $R_0$ ($R_0 =1;\, < 1;\, > 1,\, << 1;\, >> 1$) is an important flow indicator for HIV infections in the future.

Our main contributions of this paper can be  summarized as follows:
\begin{itemize}

\item Formulation of a novel immune response system of ODE which includes
 novel component integrated into the third equation of the classical three-equation HIV model, enhancing the biological realism of the system.

\item Derivation of optimality conditions using a Lagrangian approach, leading to a complete numerical scheme for solving forward, adjoint, and parameter identification problems.

\item  Proof of three a posteriori error estimates to guide  local time-mesh refinement and assess reconstruction accuracy.

\item Introduction of an adaptive time-mesh refinement strategy that locally adjusts the time discretization based on residual analysis.

\item Demonstration of method efficiency via clinical data, showing improved reconstruction of immune dynamics during the acute phase of HIV in four untreated patients.

\end{itemize}

The outline of the paper is as follows.  The 
  new  mathematical model is proposed in section \ref{sec:1-0}.
  The PIP is formulated in section 
\ref{sec:1}.
The Lagrangian approach and optimality conditions are derived in
 section \ref{sec:2}.
Numerical methods for solution of optimization problem are formulated
 in section  \ref{sec:fem}.
The a posteriori error estimates
are derived in section \ref{sec:apostframework}:
an a posteriori error estimate for the Tikhonov functional is formulated in section 
  \ref{sec:errorfunc} and
 a posteriori error estimate for the error in the reconstructed immune response function is derived in section
 \ref{sec:adaptrelax}.
 Conjugate gradient algorithm (CGA) and adaptive conjugate gradient algorithm  (ACGA)  for solution of PIP are formulated in section
\ref{sec:fem_IP}.
Finally, in section \ref{sec:numex} numerical examples
illustrate the effectiveness of the proposed ACGA using clinical data for four untreated HIV-infected   patients taken from \cite{Lit_new}.

\begin{figure}
  \begin{center}
    {\includegraphics[scale=0.3, clip=]{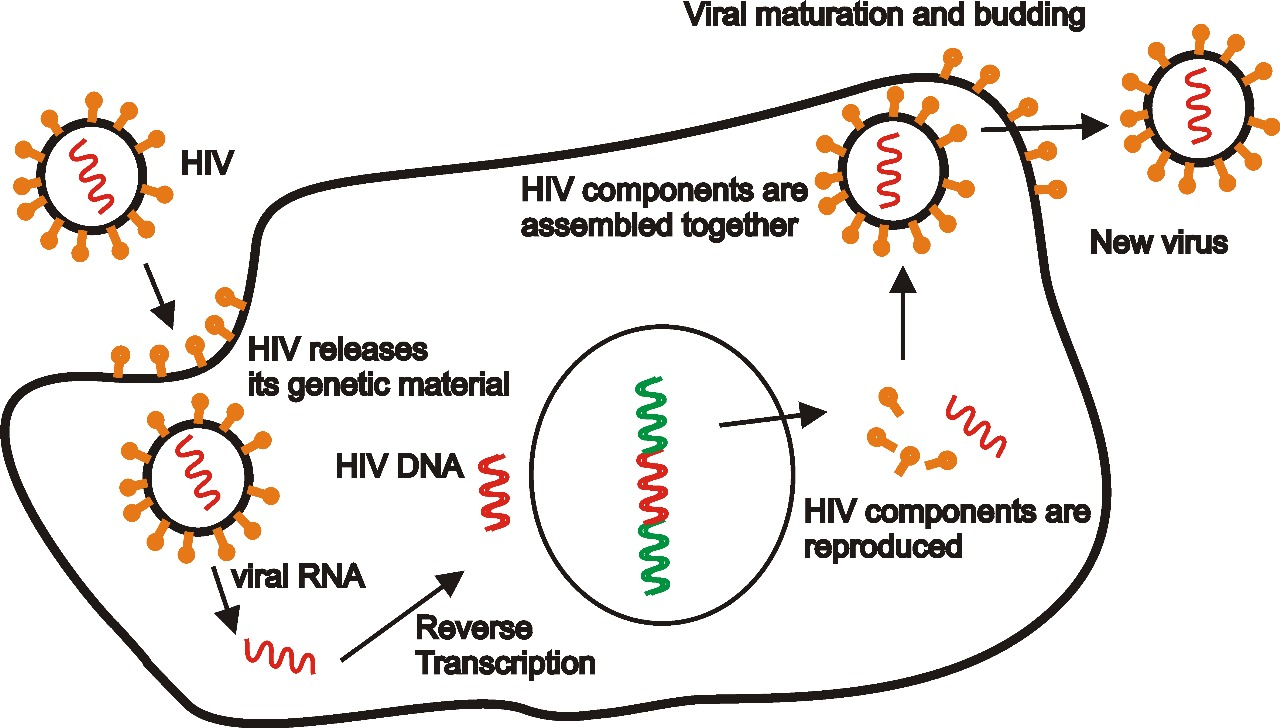}}
    \end{center}
\caption{
  \emph{The HIV life cycle.}}
\label{fig:hivc}
\end{figure}

\section{Description of the mathematical model}

\label{sec:1-0}

Recently in \cite{BErG}  was  developed time-adaptive algorithm for determination of parameters in the mathematical model of HIV infection proposed in \cite{Sriv09}. This model describes the effect of drug (Reverse Transcriptase
Inhibitor, RTI) on the dynamics of HIV infection.  Figure \ref{fig:hivc}  presents the HIV life circle described by this model.

In the current work we consider the case of HIV-1 infection without treatment
which will result in transformation of the system from \cite{Sriv09} to three equations
combined with the new component in the third equation.

Let us denote:
\begin{itemize}
\item $T(t)$ -- population of uninfected target cells,
\item $I(t)$ -- population of infected target cells,
\item $V(t)$ -- population of viral particles,
\item $E(t)$ -- function that describes effectivity of immune response by CTLs (cytotoxic T-cells).
\end{itemize}
We also adjust the dimensions involved:
$mm^3$ to $ml$ and $mm^{-3}$ to $ml^{-1}$, respectively, such that:
\[
1\, mm^3 =1 \, \mu l = 10^{-3}ml, \quad 1\, mm^{-3} = 1\, {\mu l}^{-1} = 10^3 {ml}^{-1}.
\]

Let us  denote by $\Omega_{t}= [0,T_{obs}]$ the time domain for $t>0$, where $T_{obs}$ is the final observation time.
The  model problem considered in this works is:
\begin{equation}
  \label{state_forw_new}
        \begin{cases}
{\displaystyle\frac{dT}{dt}}   = {f_1(T(t), I(t), V(t), E(t))}= s - \beta_1 T(t) V(t) - \mu T(t),\\[8pt]
{\displaystyle\frac{dI}{dt}}   = {f_2(T(t), I(t), V(t), E(t))}= \beta_1 T(t) V(t) - d(t) E(t)  I(t),\\[8pt]
{\displaystyle\frac{dV}{dt}}   = {f_3(T(t), I(t), V(t), E(t))}=  \rho I(t) - \beta_2 T(t) V(t) - c V(t),
        \end{cases}
\end{equation}
with initial conditions
\begin{equation}\label{inidata_new}
\begin{array}{lll}
  T(0)=T^0 \ cell/ml,& I(0)=I^0 \ cell/ml, & V(0)=V^0 \ virion/ml.
\end{array}
\end{equation}
The parameters used in system \eqref{state_forw_new} are described in Table~1.

\begin{table}
{\small
\begin{tabular}{llll}
  \hline
  Parameter & Value & Units & Description \\ \hline \\
  $s$ & $10^4$ & $cell/(ml\cdot day)$ & inflow rate of T cells \\[4pt]
  $\mu$ & $0.01$ & $1/day$ &natural death rate of T cells  \\[4pt]
  $\beta_1$ & $2.4$e-8  & $ml/(virion \cdot day)$ & virus infectivity rate \\[4pt]
  $\beta_2$ & $2.4$e-8  & $ml/(cell \cdot day)$ & rate of viral decline on infection of T-cells  \\[4pt]
  $d$ & $0.26$ & $1/day$ & death rate of infected cells \\[4pt]
  $c$ & $2.4$ & $1/day$ & clearance rate of virus\\[4pt]
  $\rho$ &  $ 10^3 $ & $virion/(cell \cdot day) $ & total number of viral particles produced by an infected cell \\
  \hline
\end{tabular}
}
\caption{Parameters dataset for system \eqref{state_forw_new}.}
\label{param2}
\end{table}

To simplify notations
for further analytical investigations
we again use change of variables in system  \eqref{state_forw_new}
as follows: $u_1(t) =  T(t)$, $u_2(t) = I(t)$, $u_3(t) =  V(t)$.
Then  the system \eqref{state_forw_new} with initial conditions \eqref{inidata_new}
can be  written   for $t \in \Omega_t$ as
\begin{equation}
  \label{state_forw_new2}
        \begin{cases}
{\displaystyle\frac{du_1}{dt}}  = s - \beta_1 u_1(t) u_3(t) - \mu u_1(t),\\[8pt]
{\displaystyle\frac{d u_2}{dt}}   =\beta_1 u_1(t) u_3(t) - d(t) E(t)  u_2(t),\\[8pt]
{\displaystyle\frac{d  u_3}{dt}}  = \rho u_2(t) - \beta_2 u_1(t) u_3(t) - c u_3(t), \\
u_1(0)=u_1^0 \ cell/ml,\\
u_2(0)=u_2^0 \ cell/ml, \\
u_3(0)=u_3^0 \ virion/ml.
        \end{cases}
\end{equation}

\begin{figure}
  \begin{center}
    {\includegraphics[scale=1.5, clip=]{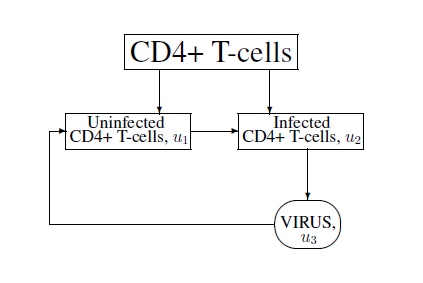}}
    \end{center}
\caption{
  \emph{Interaction between HIV and T-cells.}}
\label{fig:hivc2}
\end{figure}

Figure \ref{fig:hivc2}  shows schematic behaviour of
interactions between HIV and T-cells  described by system
\eqref{state_forw_new2}.

The problem \eqref{state_forw_new2} can be presented  in the following compact form:
\begin{equation} \label{forward_bio1}
    \begin{cases}
\dfrac{du}{dt} &= f(u(t),E(t)) ~~t \in \Omega_t, \\[6pt]
u(0) & =  u^0,
\end{cases}
    \end{equation}
 where all  involved functions  are denoted as
\begin{equation}\label{ode}
  \begin{split}
    u\,\, & =u(t)=(u_1(t),u_2(t),u_3(t))^T, \\
    u^0 &= (u_1(0),u_2(0),u_3(0))^T, \\
    \frac{du}{dt} &= \left(\frac{\partial u_1}{ \partial t}, \frac{\partial u_2}{ \partial t},\frac{\partial u_3}{ \partial t}\right)^T, \\
    f(u(t),E(t)) &=(f_1,f_2,f_3)^T(u(t),E(t))=\\
    & = ( f_1(u_1,u_2,u_3, E(t)), ...,f_3(u_1, u_2, u_3, E(t)))^T.
\end{split}
  \end{equation}
Here, $(\cdot)^T$ denotes transposition operator.

%
%
%
%
%
%
%
%
%
%
%
%
%
%
%
%
%
%


In the current work we are interesting in the development of time-adaptive algorithm
for determination of the function $E(t)$ in system
\eqref{state_forw_new2} using clinical data of four patients taken
from \cite{Lit_new}. The function $E(t)$ describes effectivity of
immune response during treatment of  HIV infection.
Determination of this parameter is very important for doctors for
analysis of each treated patient since according to \cite{Boer-2010},
a single HIV-1 infected cell produces between $10^3$ and more than
$10^4$ viral particles over its life span.  In the work \cite{Lit_new}
(see Fig.~1-b)  in  \cite{Lit_new})  clinical data  for the function $V(T)$ -- viral load in the plasma -- and  the function
$\Sigma =T(t)+I(t)$~-- total CD4+ T cell count -- are presented for
four patients with HIV without treatment. Data were measured during
one year period of the hyperacute and acute HIV infection in the 8
points (pre-infection, 0 weeks, 1 weeks, 2 weeks, 3 weeks, 4 weeks, 6
months, 1 year).
Table~\ref{tab:table3} and Table~\ref{tab:table4} summarize the clinical data.


\begin{table}[H]
  \begin{tabular}{cc}
{\small
\begin{tabular}{|c|c|c|}
  \hline
 N & $\log_{10}$ V(t)  & $\Sigma = T(t)+I(t)$ \\
 point  &  $copies \cdot ml^{-1}$ & $\times 10^3\, cells \cdot ml^{-1}$ \\ \hline
 1. & 0.0 & 1125. \\
 2. & 5.5 & 825. \\
 3. & 7.75 & 675.\\
 4. & 6.5 & 525. \\
 5. & 5.5 & 540.  \\
 6. & 5.3 & 525. \\
 7. & 4.5 & 550.\\
 8. & 5.1 & 600.  \\
  \hline
\end{tabular}
}
&
{\small
\begin{tabular}{|c|c|c|}
  \hline
 N & $\log_{10}$ V(t)  & $\Sigma = T(t)+I(t)$ \\
 point  &  $copies \cdot ml^{-1}$ & $\times 10^3\, cells \cdot ml^{-1}$  \\ \hline
 1.  & 0.0 & 750. \\
 2. & 5.5 & 630. \\
 3. & 7.0 & 450. \\
 4.  & 5.75 & 280. \\
 5. & 4.8 & 750. \\
 6.  & 3.8 & 570. \\
 7. & 4.3 & 450. \\
 8. & 4.0 & 530. \\
  \hline
\end{tabular}
}
\\
 Patient 1 &   Patient 2 \\
\end{tabular}
\caption{\textit{Clinical data corresponding to Fig.~1-b of \cite{Boer-2010}.}}
\label{tab:table3}
\end{table}


\begin{table}[H]
  \begin{tabular}{cc}
{\small
\begin{tabular}{|c|c|c|}
  \hline
 N & $\log_{10}$ V(t)  & $\Sigma = T(t)+I(t)$ \\
 point  &  $copies \cdot ml^{-1}$ & $\times 10^3\, cells \cdot ml^{-1}$ \\ \hline
 1. & 0.0 & 700.\\
 2. & 5.7  & 430.\\
 3. & 7.4 & 300.\\
 4. & 6.8 & 480.\\
 5. & 4.2 & 570.\\
 6. & 3.8 & 450. \\
 7. & 3.2 & 630.\\
 8. & 3.7 & 510. \\
  \hline
\end{tabular}
}
&
{\small
\begin{tabular}{|c|c|c|}
  \hline
 N & $\log_{10}$ V(t)  & $\Sigma = T(t)+I(t)$ \\
 point  &  $copies \cdot ml^{-1}$ & $\times 10^3\, cells \cdot ml^{-1}$  \\ \hline
  1. & 0.0 & 615. \\
  2. &3.7 & 700. \\
  3.  & 5.7 & 450. \\
  4. & 6. & 615. \\
  5. & 3.9 & 575. \\
  6. & 3.5 & 520. \\
  7. & 3. & 450. \\
  8. & 3. & 615. \\
  \hline
\end{tabular}
}
\\
Patient 3 & Patient 4 \\
\end{tabular}
\caption{\textit{Clinical data  corresponding to Fig.~1-b of \cite{Boer-2010}.}}
\label{tab:table4}
\end{table}


To proceed further we introduce several mathematical assumptions which are necessary for development of
 a stable
 reconstruction algorithm.
 First, we assume that in the  system
\eqref{forward_bio1}
 the function $f \in C^1(\Omega_{t})$
 is Lipschitz continuous. Next, accordingly to work \cite{Stafford}
 we assume that
 the function $E(t) \in C(\Omega_{t})$  belongs to the following set of admissible parameters $M_{E}$:
\begin{equation}\label{6.1}
  M_{E} = \{ E(t) \in C(\Omega_{t}) : E(t)~ \in~ \left[1, 10 \right]~
\textrm{in} ~\Omega_{t}
\}.
\end{equation}



\section{The parameter identification problem}

\label{sec:1}

To formulate the parameter identification problem we assume that all
parameters in system \eqref{forward_bio1} are known except the control
parameter $E(t)$ which presents time-dependent distribution of viral
particles produced by one infected cell.
The values of other parameters $\{s, \mu, \beta, \rho, c \}$, and   $d(t) = d$ in
 the model \eqref{forward_bio1} correspond to values in the  Table~\ref{param2}.
 We note that  though  in the Table~1
   the  death rate function of infected T-cells $d(t)$
 is
defined as a constant, i.e.,  $d(t) = d$, in numerical experiments  we take the function $d(t)$  as a time-dependent function of a special form - see more details in section \ref{sec:numex}.

\medskip

\textbf{Parameter Identification Problem (}PIP\textbf{)}. Assume
that the condition (\ref{6.1}) for the function $E(t)$
 holds and parameters
 $\{s , \mu, \beta_1, \beta_2, d, \rho, c \}$ in system
\eqref{forward_bio1} are known. Assume further that the function $E(t) \in M_E$
is unknown inside the domain $\Omega_{t}$. The PIP is: determine $E(t)$ for $t
\in \Omega_{t}$,  under the conditions  that the population function $g_1(t)$ of total number of the uninfected and infected T cells is known:
\begin{equation}
u_1(t)+u_2(t) \approx g_1(t),~~ t \in  \Omega_t,  \label{6.4a}
\end{equation}
as well as the virus population function $g_2(t)$ is known:
\begin{equation}
u_3(t) \approx g_2(t),~~ t \in \Omega_t.  \label{6.4b}
\end{equation}
Here, the function $g_1\left(t\right) $ presents observations of the sum $u_1\left(t\right) + u_2\left(t\right)$,  or unifected and infected T cells, and $g_2\left(t\right) $ presents observations of the  virus function $u_3\left(t\right)$
over the domain  $\Omega_t$.



\section{Optimization method}

\label{sec:2}

Let $H$ be a Hilbert space of functions defined in $\Omega_t$.
To determine $E(t)$, $t\in \Omega_t$,  we  construct the
Tikhonov functional  which is more appropriate for the data presented in Tables
\ref{tab:table3}  and \ref{tab:table4}:
\begin{equation}\label{Tikh_f0}
\begin{array}{ll}
  J(E) & =
\displaystyle{\frac{1}{2}} \int\limits_{\Omega_t}
  [
  \log_{10}(u_1(t)  + u_2(t)) - \log_{10}(g_1(t))]^{2} {z_1}_{\zeta}\left(  t\right)
~\mathrm{d}t
\\[12pt]
 & +
\displaystyle{\frac{1}{2}} \int\limits_{\Omega_t} [\log_{10} u_3(t) - \log_{10} g_2(t)]^{2}{z_2}_{\zeta}\left(  t\right)
~\mathrm{d}t
+\displaystyle{\frac{1}{2}}\, {\gamma}
\int\limits_{\Omega_t} (E(t) - E^{0})^{2}dt.
\end{array}
\end{equation}

Here,  $u_1(t)$, $u_2(t)$, $u_3(t)$  are the solutions of the system
\eqref{forward_bio1}  which are depended on the
  function $E(t)$, $g_1(t)$ is the observed sum of uninfected and infected T-cells, $g_2(t)$ is the observed virus population function,
$E^0$ is the initial guess for the parameter $E(t)$ and
$\gamma\in (0,1)$ is the regularization parameter, ${z_1}_{\zeta}(t),
\zeta \in \left( 0,1\right) $ and ${z_2}_{\zeta}(t),
\zeta \in \left( 0,1\right) $ are smoothing functions for data which can be
defined similarly  as in \cite{Sprg-2}.

To  find the function $E(t) \in H$ which minimizes the Tikhonov functional (\ref{Tikh_f0}) we
 seek for a stationary point of (\ref{Tikh_f0})  with respect to $E$  such that
\begin{equation}\label{minimum}
J^{\prime }(E)(\bar{E})=0, ~~\forall \bar{E} \in H.
\end{equation}

To find minimum of (\ref{Tikh_f0}) we use constrained optimization
with the standard Lagrangian approach \cite{alb, Polak} and introduce the  Lagrangian  with constrains
corresponding to the model problem:
\begin{equation}\label{Lagran_eta}
  L(v)=J(E)+  \sum_{i=1}^3 \int\limits_{\Omega_t} \lambda_i
  \left(\frac{du_i}{dt} - f_i \right) ~dt,
\end{equation}
where $u(t)=(u_1(t),u_2(t),u_3(t))$ is the solution of the system
\eqref{forward_bio1}, $\lambda(t)$ is the vector of Lagrange multipliers $\lambda(t)=(\lambda_1(t),\lambda_2(t),\lambda_3(t))$, and $v= (\lambda,u,E)$.

We now introduce  following spaces   which we are using in analysis  of PIP:
\begin{equation}
\begin{array}{rl}
H_{u}^{1}(\Omega_t)  &  =\{ u \in H^{1}(\Omega_t): u(0)= u_0\},\\
H_{\lambda}^{1}(\Omega_t)  &  =\{ \lambda \in H^{1}(\Omega_t):
\lambda(T)=0\},\\
U  &=H_{u}^{1}(\Omega_t)\times H_{\lambda}^{1}(\Omega_t)\times C(\Omega_t),
\end{array}
\label{spaces}
\end{equation}
for all   real valued functions.

To derive the Fr\'{e}chet derivative of the Lagrangian (\ref{Lagran_eta})
 we use such called all-at-once approach when
we assume that functions $v=(\lambda,u,E)$ can be varied independently
of each other  such that
\begin{equation}\label{Frechet_der}
  L'(v)(\bar{v})=0,\quad \forall  \bar{v} =(\bar{\lambda},\bar{u},\bar{E}) \in U.
\end{equation}
The   optimality condition \eqref{Frechet_der} means also that for all
  $\bar{v} \in U$ we  have
\begin{equation}
  L'(v; \bar{v}) = \frac{\partial L}{\partial \lambda}(v)(\bar{\lambda}) +  \frac{\partial L}{\partial u}(v)(\bar{u})
  + \frac{\partial L}{\partial E}(v)(\bar{E}) = 0  \label{scalar_lagrang1}
\end{equation}
i.e., every component of \eqref{scalar_lagrang1} should be zero
out.
Thus, the optimality condition $\frac{\partial L}{\partial \lambda}(v)(\bar{\lambda}) = 0$   yields
\begin{equation}\label{forward1}
\begin{split}
  0 = \frac{\partial L}{\partial \lambda_1}(v)(\bar{\lambda_1}) &=
   \int\limits_{\Omega_t}( \frac{\partial{u_1}}{\partial t}  - s +  \beta_1 u_1 u_3 + \mu u_1   ) \bar{\lambda_1} dt  ~~~\forall \bar{\lambda_1} \in H_\lambda^1(\Omega_t), \\
  0 = \frac{\partial L}{\partial \lambda_2}(v)(\bar{\lambda_2})
  &= \int\limits_{\Omega_t} ( \frac{\partial{u_2}}{\partial t}  - \beta_1 u_1 u_3 +  d E u_2  ) \bar{\lambda_2} dt,~~~\forall \bar{\lambda_2} \in H_\lambda^1(\Omega_t)\\
 0 = \frac{\partial L}{\partial \lambda_3}(v)(\bar{\lambda_3})
  &= \int\limits_{\Omega_t} ( \frac{\partial{u_3}}{\partial t}  -  \rho u_2 + c u_3 + \beta_2 u_1 u_3 ) \bar{\lambda_3} dt,
~~~\forall \bar{\lambda_3} \in H_\lambda^1(\Omega_t),
\end{split}
\end{equation}

Next, the optimality condition $\frac{\partial L}{\partial u}(v)(\bar{u}) = 0$
means that
\begin{equation} \label{control1}
\begin{split}
  0 = \frac{\partial L}{\partial u_1}(v)(\bar{u_1}) &=
  \int\limits_{\Omega_t}
  \frac{1}{(u_1 + u_2) \ln{10} }   [  \log_{10}(u_1(t)  + u_2(t)) - \log_{10}(g_1(t))] {z_1}_{\zeta}\left(  t\right)
  ~ \bar{u_1}\mathrm{d}t \\
  &+  \int\limits_{\Omega_t} (- \frac{\partial{\lambda_1}}{\partial t}  + \lambda_1  \beta_1  u_3 + \mu \lambda_1 - \lambda_2 \beta_1 u_3  + \beta_2 \lambda_3 u_3 ) \bar{u_1} dt  ~~~\forall \bar{u_1} \in H_u^1(\Omega_t), \\
   0 = \frac{\partial L}{\partial u_2}(v)(\bar{u_2}) &=
  \int\limits_{\Omega_t}
  \frac{1}{(u_1 + u_2) \ln{10}  }   [  \log_{10} (u_1(t)  + u_2(t)) - \log_{10}(g_1(t))] {z_1}_{\zeta}\left(  t\right)
  ~ \bar{u_2}\mathrm{d}t \\
  &+  \int\limits_{\Omega_t} (- \frac{\partial{\lambda_2}}{\partial t} + d E \lambda_2 - \lambda_3 \rho ) \bar{u_2} dt  ~~~\forall \bar{u_2} \in H_u^1(\Omega_t), \\
  0 = \frac{\partial L}{\partial u_3}(v)(\bar{u_3}) &=
  \int\limits_{\Omega_t}
  \frac{1}{u_3 \ln{10}  }   [  \log_{10} (u_3) - \log_{10} (g_2)] {z_1}_{\zeta}\left(  t\right)
  ~ \bar{u_3}\mathrm{d}t \\
  &+  \int\limits_{\Omega_t} (- \frac{\partial{\lambda_3}}{\partial t}  + \lambda_1 \beta_1 u_1 - \lambda_2 \beta_1 u_1 + c \lambda_3  + \lambda_3 \beta_2 u_1) \bar{u_3} dt  ~~~\forall \bar{u_3} \in H_u^1(\Omega_t).
\end{split}
\end{equation}

Finally, the optimality condition $\frac{\partial L}{\partial E}(v)(\bar{E}) = 0$  yields:
\begin{equation} \label{grad1new}
\begin{split}
0 &= \frac{\partial L}{\partial E}(v)(\bar{E})
= \gamma \int\limits_{\Omega_t}  (   E -  E^{0})\bar{E} dt
 + \int_{\Omega_t} d \lambda_2 u_2 \bar{E}    dt
~~~\forall \bar{E} \in C\left(\Omega_t \right).
\end{split}
\end{equation}

\medskip


Using optimality conditions \eqref{forward1}-\eqref{control1}
we observe that the  equations  \eqref{forward1} correspond to the system of model equations
\eqref{state_forw_new2}  which is  also  our forward problem,  and equations
  \eqref{control1} --- to  the  following  adjoint problem
\begin{equation}\label{adjoint}
    \begin{split}
      \frac{\partial \lambda_1}{\partial t}  = \tilde{f}_1(\lambda(t),\eta(t)) &=  \lambda_1(t) \beta_1 u_{3}(t)
      + \lambda_1(t) \mu  - \lambda_2(t) \beta_1 u_{3}(t)  + \beta_2 \lambda_3(t) u_3(t)  \\
      &+  \frac{1}{(u_1(t) + u_2(t)) \ln{10}  }   [  \log_{10}(u_1(t)  + u_2(t)) - \log_{10} (g_1)]  {z_1}_{\zeta}, \\
      \frac{\partial \lambda_2}{\partial t}  = \tilde{f}_2(\lambda(t),\eta(t)) &=  -\lambda_3(t) \rho + d \lambda_2(t) E(t)  \\
      &+   \frac{1}{(u_1(t) + u_2(t))\ln{10}   }   [  \log_{10}(u_1(t)  + u_2(t)) - \log_{10}(g_1)]  {z_1}_{\zeta},\\
      \frac{\partial \lambda_3}{\partial t}  = \tilde{f}_3(\lambda(t),\eta(t)) &= \lambda_1(t) \beta_1 u_1(t) - \lambda_2(t) \beta_1 u_1(t) + c \lambda_3(t) + \lambda_3(t)\beta_2 u_1(t) \\
      &+  \frac{1}{u_3(t) \ln{10}  }   [  \log_{10} (u_3(t)) - \log_{10} (g_2)] {z_1}_{\zeta}, \\
\lambda_i(T) &= 0,\quad i=1,2,3,
\end{split}
    \end{equation}
which we rewrite  as
\begin{equation}\label{adjointcompact}
    \begin{cases}
    \frac{\partial \lambda}{\partial t} =  \tilde{f}(\lambda(t),E(t)),\\
    \lambda_i(T) = 0,\quad i=1,2,3,
    \end{cases}
\end{equation}
with
\begin{equation}
  \begin{split}
    \lambda= \lambda(t) &=(\lambda_1(t), \lambda_2(t), \lambda_3(t))^T, \\
    0&= (\lambda_1(T), \lambda_2(T), \lambda_3(T))^T, \\
    \frac{d\lambda}{dt} &= { \left(\frac{\partial \lambda_1}{ \partial t}, \frac{\partial \lambda_2}{ \partial t},\frac{\partial \lambda_3}{ \partial t} \right)^T}, \\
    \tilde{f}(\lambda(t)) &=( \tilde{f}_1,  \tilde{f}_2,  \tilde{f}_3)( \lambda(t),E(t))^T.
\end{split}
  \end{equation}

The adjoint   problem should be  solved backwards in time with already known forward problem solution $u(t)$ and a given  data $g_1(t), g_2(t)$.


From definition of the Lagrangian  \eqref{Lagran_eta} 
 it follows that
 for the  known  exact  function $u$,  known parameters in the forward problem
and known solution $\lambda$ of the  adjoint  problem (\ref{adjoint})
  we get that
\begin{equation}
L(v(E)) = J(E).
\end{equation}
In this case the Fr\'{e}chet derivative of the Tikhonov functional can be
written as
\begin{equation}\label{derfunc}
\begin{split}
 J'(E) := &J_E(u(E), E) =  \frac{\partial J}{\partial E}(u(E), E)
  =  \frac{\partial L}{\partial E}(v(E)).
\end{split}
\end{equation}
Using (\ref{grad1new})  in  (\ref{derfunc}),
we get the following expression for the Fr\'{e}chet derivative of the Tikhonov functional
which will be used later in the conjugate gradient update of the parameter $E$ in time :
\begin{equation}\label{Frder2}
  J'(E)(t)  = \gamma ( E -  E^{0})(t) +  d  \lambda_2(t) u_2(t).
\end{equation}


\section{Numerical methods  for solution of  optimization problem}
\label{sec:fem}

To  solve the minimization problem \eqref{Frechet_der}  we will employ the finite element  method.
To do this
we discretize the computational domain $\Omega_t = [0,T]$  and    consider a partition $ J_{\tau} = \{J_k\}$ of the time
domain $\Omega_t$
 into time subintervals $J_k=(t_{k-1},t_k]$
  with the  time step $\tau_k =  t_k - t_{k-1}$. Let us define  the piecewise-constant
  function $\tau(t)$ such that
\begin{equation}\label{neshfunction}
\tau(t) = \tau_k, ~~\forall J_k \in  J_{\tau}.
\end{equation}

We   define also the
finite element spaces $V_{\tau}^{u}\subset H_{u}^{1}\left(\Omega_t\right)
$ and $V_{\tau}^{\lambda}$ $\subset H_{\lambda} ^{1}\left(\Omega_t\right)
$ for approximations $u_\tau, \lambda_\tau$ of $u$ and $\lambda$, respectively, as
\begin{equation}
\begin{array}{rl}
V_{\tau}^{u}  &=\{ u \in  L_2(\Omega_t):  u|_J \in P^0(J), u(0) = u^0, ~~~ \forall J \in J_{\tau}\},\\
V_{\tau}^{\lambda} &=\{\lambda \in   L_2(\Omega_t): \lambda|_J \in P^0(J), \lambda(T)=0, ~~~ \forall J \in J_{\tau} \}
\end{array}
\label{femspaces}
\end{equation}
for discretization of the state and adjoint problems.
We also introduce the finite element space
$V_{\tau}^{E}\subset L_2\left(\Omega_t \right) $ for the function $E(t)$ consisting of piecewise
constant functions
\begin{equation}
V_{\tau}^{E} =\{f\in  L_{2}\left( \Omega_t \right) :  f|_J \in P^0(J)~~ \forall J \in J_{\tau}\}.
\end{equation}
 Let $U_{\tau}= V_{\tau}^{u}\times V_{\tau}^{\lambda}\times V_{\tau}^{E}$
such that $U_{\tau}\subset U$.
We can formulate now the finite element method for (\ref{Frechet_der}):
find $v_{\tau}\in U_{\tau}$ such that
\begin{equation}
L^{\prime}(v_{\tau})(\bar{v})  =0,~~\forall\overline{v}\in U_{\tau}.
\label{discr_lagr_ac}
\end{equation}

We note that the forward \eqref{state_forw_new2} and adjoint \eqref{adjoint}  problems are nonlinear.
For solution of these problems we employ  Newton's method as follows.
The variational formulation of the forward problem \eqref{forward_bio1}   is:
\begin{equation} \label{varform1}
  \begin{split}
\left (\dfrac{du}{dt},\bar{u} \right)  &= (f(u(t)), \bar{u} ) ~~\forall \bar{u} \in H_u^1(\Omega_t).
\end{split}
\end{equation}
Let us use  the following discretization  for the derivative in time
\begin{equation*}
\frac{\partial u}{\partial t} = \frac{u^{k+1} - u^{k}}{\tau_k}
\end{equation*}
in the variational formulation \eqref{varform1} to obtain
\begin{equation} \label{varform2}
  \begin{split}
 (u^{k+1},\bar{u})  - \tau_k (f(u^{k+1}), \bar{u} ) -  (u^{k},\bar{u}) = 0 ~~\forall \bar{u} \in H_u^1(\Omega_t).
\end{split}
\end{equation}
Approximating $u^{k+1}$  by its finite element approximation  $u_\tau^{k+1}$  the above equation can be written as
\begin{equation} \label{varform3}
  \begin{split}
 (u_\tau^{k+1} - \tau_k f(u_\tau^{k+1})  -  u_\tau^{k},\bar{u}) = 0 ~~\forall \bar{u} \in
  V_\tau^u.
\end{split}
\end{equation}
We denote now
\begin{equation}\label{4.26a}
  \begin{split}
  \tilde{u} &= u_\tau^{k+1}, \\
  V(\tilde{u}) & = \tilde{u} - \tau_k f(\tilde{u}) - u_\tau^k.
\end{split}
  \end{equation}
Next, we rewrite the variational formulation \label{varform3} as
\begin{equation} \label{varform4}
( V(\tilde{u}),\bar{u}) = 0 ~~\forall \bar{u} \in   V_\tau^u.
  \end{equation}
To solve nonlinear equation $V(\tilde{u})=0$ we are using  the Newton's method  \cite{Burden} for iterations in the Newtons' method $n = 1,....M_n$:
\begin{equation}\label{4.31}
\tilde{u}^{n+1} = \tilde{u}^n - [ V'(\tilde{u}^n)]^{-1} \cdot V(\tilde{u}^n).
  \end{equation}
Here,  the Jacobian $V'(\tilde{u}^n)$ is computed via  definition of  $V(\tilde{u})$  in \eqref{4.26a} as
\begin{equation*}
V'(\tilde{u}^n) = I - \tau_k f'(\tilde{u}^n),
 \end{equation*}
where $I$ is the identity matrix, $f'(\tilde{u}^n)$ is the Jacobian of
$f$  which represents the right hand side of the forward problem   \eqref{state_forw_new2} at $\tilde{u}^n$,  $n$ is the iteration number in Newton's method and $M_n$ is the final iteration in Newton's method.
The entries in the Jacobian  $f'(\tilde{u}^n)$ for system \eqref{state_forw_new2} can be   explicitly computed and they are  defined as:
\begin{equation*}\label{eq20}
 f'(\tilde{u}^n) =\left[
\begin{array}{cccc}
  \dfrac{\partial f_1}{\partial u_1} & \dfrac{\partial f_1}{\partial u_2} & \dfrac{\partial f_1}{\partial u_3}  \\[6pt]
 \dfrac{\partial f_2}{\partial u_1} & \dfrac{\partial f_2}{\partial u_2} & \dfrac{\partial f_2}{\partial u_3}  \\[6pt]
 \dfrac{\partial f_3}{\partial u_1} & \dfrac{\partial f_3}{\partial u_2} & \dfrac{\partial f_3}{\partial u_3} \\[6pt]
\end{array}
 \right](\tilde{u}^n) = \left[
\begin{array}{cccc}
    - \beta_1 {u_3}_\tau^{n} - \mu  & 0 &  -\beta_1  {u_1}_\tau^{n} \\[6pt]
 \beta_1 {u_3}_\tau^{n}   & - d ~E & \beta_1  {u_1}_\tau^{n} \\[6pt]
   -\beta_2 {u_3}_\tau^n &  \rho  & -\beta_2 {u_1}_\tau^n -c
\end{array}
 \right].
\end{equation*}

We derive the Newtons's method  for the solution the adjoint problem \eqref{adjoint} in a similar way.
Since we solve the adjoint problem backwards in time starting from  the  known  $\lambda(T)=0$,
we discretize the
time derivative as
\begin{equation}\label{4.32}
\frac{\partial \lambda}{\partial t} = \frac{\lambda^{k+1} - \lambda^{k}}{\tau_k}
  \end{equation}
for the already known  $\lambda^{k+1}$, and write the variational formulation of the adjoint problem for all
$\bar{\lambda} \in H_\lambda^1(\Omega_t)$  as
\begin{equation}\label{4.33old}
 ( \lambda^{k+1} - \lambda^{k}    - \tau_k \tilde{f}(\lambda^k) ,\bar{\lambda}) = 0.
\end{equation}
We  rewrite \eqref{4.33old} in the form
\begin{equation}\label{4.33new}
 ( \lambda^{k} - \lambda^{k+1}    + \tau_k \tilde{f}(\lambda^k) ,\bar{\lambda}) = 0
\end{equation}
and
formulate the finite element method  for the  variational formulation
 \eqref{4.33new}:
 find
 $\lambda_\tau^k \in V_\tau^\lambda$   such that for all $\bar{\lambda} \in  V_\tau^\lambda$
\begin{equation}\label{4.33}
( \lambda_\tau^{k}  - \lambda_\tau^{k+1} + \tau_k \tilde{f}(\lambda_\tau^k) ,\bar{\lambda}) = 0.
\end{equation}

Denoting
\begin{equation}\label{4.34}
  \begin{split}
  \tilde{\lambda} &= \lambda_\tau^{k}, \\
  R(\tilde{\lambda}) &= \tilde{\lambda} + \tau_k \tilde{f}(\tilde{\lambda})  -  \lambda_\tau^{k+1},
\end{split}
  \end{equation}
  we can rewrite \eqref{4.33}  for all
$\bar{\lambda} \in  V_\tau^\lambda$ as
\begin{equation}\label{4.35}
(R(\tilde{\lambda}),  \bar{\lambda}) = 0.
  \end{equation}
For solution $ R(\tilde{\lambda})=0$  we use  again Newton's method:
\begin{equation}\label{4.36}
\tilde{\lambda}^{n+1} = \tilde{\lambda}^n - [ R'(\tilde{\lambda}^n)]^{-1} \cdot R(\tilde{\lambda}^n), n=1,...., K_n.
  \end{equation}
We compute  $R'(\tilde{\lambda}^n)$ using the  definition of  $R(\tilde{\lambda})$  in \eqref{4.34} as
\begin{equation*}
 R'(\tilde{\lambda}^n) = I + \tau_k \tilde{f}'(\tilde{\lambda}^n),
 \end{equation*}
where $I$ is the identity matrix, $\tilde{f}'(\tilde{\lambda}^n)$ is the Jacobian of
$\tilde{f}$ at $\tilde{\lambda}^n$,  $n$ is the iteration number in the
Newton's method and $K_n$ is the final number of iteration.
The  entries in the Jacobian  $\tilde{f}'(\tilde{\lambda}^n)$ for the adjoint
system \eqref{adjoint} are  explicitly given by
\begin{equation*}
  \tilde{f}'(\tilde{\lambda}^n) =
  \left[
  \begin{array}{cccc}
  \dfrac{\partial \tilde{f}_1}{\partial \lambda_1} & \dfrac{\partial \tilde{f}_1}{\partial \lambda_2} & \dfrac{\partial \tilde{f}_1}{\partial \lambda_3}  \\[6pt]
 \dfrac{\partial \tilde{f}_2}{\partial \lambda_1} & \dfrac{\partial \tilde{f}_2}{\partial \lambda_2} & \dfrac{\partial \tilde{f}_2}{\partial \lambda_3} \\[6pt]
 \dfrac{\partial \tilde{f}_3}{\partial \lambda_1} & \dfrac{\partial \tilde{f}_3}{\partial \lambda_2} & \dfrac{\partial \tilde{f}_3}{\partial \lambda_3}  \\[6pt]
\end{array}
  \right](\tilde{\lambda}^n)=
  \left[
\begin{array}{cccc}
    \beta_1 {u_3}_\tau^n +\mu &  -\beta_1 {u_3}_\tau^n &  \beta_2 {u_3}_\tau^n \\[6pt]
   0 & d E  & -  \rho \\[6pt]
   \beta_1 {u_1}_\tau^n & -\beta_1 {u_1}_\tau^n  & c + \beta_2 {u_1}_\tau^n
\end{array}
 \right].
\end{equation*}

We note that $\det f'(\tilde{u}^n) \neq 0$ as well as $ \det
\tilde{f}'(\tilde{\lambda}^n) \neq 0$ since all parameters used in the
model problem are non-zero.
Thus, schemes given by
formulas \eqref{4.31}, \eqref{4.36} will converge to the solutions of
forward and adjoint problems for the appropriate  initial guesses for values
$\tilde{u}^1$ and $\tilde{\lambda}^1$, correspondingly. 


\section{A Posteriori Error Estimates }

\label{sec:apostframework}

In this section we derive a posteriori error estimates for the error
in The Tikhonov's functional $||J(E) - J(E_h)||_{L_2(\Omega_t)} $ and for the error in
the reconsructed function $||E - E_h||_{L_2(\Omega_t)}$, correspondingly.  Both error
estimators can be used as the mesh refinement   criteriors in the adaptive conjugate gradient algorithm.

Before  deriving  error estimates we make some assumptions
which are usual for parameter identification problems since they are
ill-posed problems \cite{Tikhonov,T1,T2}.
We assume that  the function $E(t) \in C(\Omega_t)$ as a minimizer of the
Lagrangian \eqref{Lagran_eta}, and
$E_{\tau} \in V_{\tau}^E$  is
its  finite element
approximation.
We also assume that $E^*  \in  C(\Omega_t)$ is the exact  function describing the  effectivity of the immune response 
and that
there exists a good approximation $E_{\tau}$  to the exact function
 $E^* \in  C(\Omega_t)$.
Let $g_{1,2}^{\ast }(t)$ be the exact data
 for the measured functions  $g_{1,2}(t)$   defined in (\ref{6.4a}), (\ref{6.4b}), correspondingly,
and the functions
${g_{1,2}}_{\sigma }(t)$ represent the error levels in the 
 measured data. This means that we assume
that measurements $g_{1,2}(t)$ in (\ref{6.4a}), (\ref{6.4b})  are given with some small noise
levels $\sigma_1, \sigma_2$ such that 
\begin{equation}\label{6.18}
\begin{split}
g_1(t) &= g_1^{\ast }(t) + {g_1}_{\sigma}(t);\,~ g_1^{\ast }, {g_1}_{\sigma }\in L_{2}\left( \Omega_{t}\right) ,~ \left\Vert {g_1}_{\sigma }\right\Vert _{L_{2}\left( \Omega_{t}\right) }\leq \sigma_1,  \\
g_2(t) &= g_2^{\ast }(t) + {g_2}_{\sigma}(t);\,~ g_2^{\ast }, {g_2}_{\sigma }\in L_{2}\left( \Omega_{t}\right) ,~ \left\Vert {g_2}_{\sigma }\right\Vert _{L_{2}\left( \Omega_{t}\right) }\leq \sigma_2.
\end{split}
\end{equation}


We also assume that the regularization parameter in the Tikhonov's functional is such that
\begin{equation}\label{6.19}
\gamma = \gamma(\sigma_1,\sigma_2) = (\sigma_1 + \sigma_2)^{2 \mu},~~\mu \in(0,1/4), ~~\sigma_{1,2} \in (0,1)
\end{equation}
and for the exact solution  $E^*$ corresponding to the exact data  $g_{1,2}^*(t)$ it holds that
\begin{equation}\label{6.20}
\| E_0 - E^* \|_{L_2(\Omega_t)} \leq \frac{(\sigma_1 + \sigma_2)^{3\mu}}{3}.
\end{equation}
The condition \eqref{6.20} means that $E_0$  is located in the close neighborhood of the exact solution $E^*$.

 Let us define the set  for any $\varepsilon > 0$
\begin{equation}
V_\varepsilon(E) = \{ \tilde{E} \in  C(\Omega_t): \| E - \tilde{E}  \| < \varepsilon ~~~\forall E \in C(\Omega_t)  \}.
  \end{equation}
Assume  that
for all $E \in V_1(E^*)$ the operator
\begin{equation*}
\begin{split}
F(E) &=  \displaystyle{\frac{1}{2}} \int\limits_{\Omega_t}
  [ \log_{10}(u_1(E,t)  + u_2(E,t)) - \log_{10}(g_1(t))]^{2} {z_1}_{\zeta}\left(  t\right)
~\mathrm{d}t
\\[12pt]
 & +
\displaystyle{\frac{1}{2}} \int \limits_{\Omega_t} [\log_{10} u_3(E,t) - \log_{10} g_2(t)]^{2}{z_2}_{\zeta}\left(  t\right)
~\mathrm{d}t
\end{split}
 \end{equation*}
has the
Fr\'{e}chet derivative  $F'(E)$ which is bounded and Lipshitz continuous
in  $V_1(E^*)$ for $ N_1, N_2 = const. > 0$
\begin{equation}\label{6.21}
  \begin{split}
  \| F'(E) \|_{L_2(\Omega_t)} &\leq N_1 ~ ~~ \forall E \in V_1(E^*), \\
  \| F'(E_1) -  F'(E_2)  \|_{L_2(\Omega_t)} &\leq N_2 \| E_1 - E_2 \|_{L_2(\Omega_t)}  ~~ \forall E_1, E_2 \in  V_1(E^*).
\end{split}
  \end{equation}

\subsection{An a posteriori error estimate for the Tikhonov functional}

\label{sec:errorfunc}

In this section we derive an a posteriori error estimate for the error  $\| J(E) - J(E_h)\|_{L_2(\Omega_t)} $
in
the Tikhonov functional (\ref{Tikh_f0}) on the finite element time partition $ J_{\tau}$.

\textbf{Theorem 1}.
\emph{\ Let   $E \in
  C(\Omega_t)$ is the minimizer of the functional $J(E)$ defined by (\ref{Tikh_f0}) and
  there exists finite element approximation of a minimizer
  $E_{\tau} \in V_{\tau}^{E}$ of
  $J(E)$.  Then the following
   a posteriori error estimate for the error $ e=|| J(E) - J(E_{\tau}) ||_{L^2(\Omega_t)}$ in the Tikhonov functional (\ref{Tikh_f0}) holds  true}
\begin{equation}\label{theorem2}
  e= || J(E) - J(E_{\tau}) ||_{L^2(\Omega_t)} \leq C_I \left\| J^{\prime }(E_{\tau})\right\| _{L^2(\Omega_t)} \left( || \tau E_{\tau} ||_{L_2(\Omega_t)}
+
  \sum_{J_k} \left \Vert [E_\tau] \right \Vert_{L_2(J_k)} \right) ~ \forall E_\tau \in  V_\tau^E,
\end{equation}
where $C_I  = const. > 0$
and 
\begin{equation}
J^{\prime }(E_{\tau}) = \gamma (E_{\tau} -E^{0}) +
 d  {\lambda_2}_{\tau} {u_2}_\tau.
 \label{theorem2_1}
\end{equation}

\textbf{Proof}  Using first two terms in the Taylor's expansion we get
\begin{equation}\label{theorem2_a}
J(E)  = J(E_{\tau}) +  J'(E_{\tau})(E - E_{\tau}) + R(E, E_{\tau}),
\end{equation}
where $R(E, E_{\tau}) =O(\| E - E_{\tau} \|^2),~~ (E -
E_{\tau}) \to 0 ~~\forall E, E_{\tau} \in
 V_{\tau}^{E}$.
 We neglect the remainder term $R(E, E_{\tau})$ in
(\ref{theorem2_a})
since  it is small because of
assumption \eqref{6.20}.
Next, we use the following splitting  for $E - E_{\tau}$ in \eqref{theorem2_a} 
\begin{equation}\label{splitting}
E - E_{\tau} =  E - E_\tau^I +  E_\tau^I -  E_{\tau}
\end{equation}
  together with  the Galerkin orthogonality  property
\begin{equation}
J'(E_{\tau})( E_\tau^I -  E_{\tau}) = 0,\quad \forall E_\tau^I, E_{\tau} \in
V_{\tau}^{E}
\end{equation}
to get
\begin{equation}\label{5.37}
J(E) -J(E_{\tau}) \leq J'(E_{\tau})(E - E_\tau^I).
\end{equation}
Here, $E_\tau^I$ is a standard interpolant of $E$ on the mesh  $ J_{\tau}$ \cite{ERJ}.
Taking norms in \eqref{5.37}, we obtain
\begin{equation}\label{5.38}
||J(E) -J(E_{\tau}) ||_{L^2(\Omega_t)} \leq ||J'(E_{\tau})||_{L^2(\Omega_t)} ||E - E_\tau^I||_{L^2(\Omega_t)}.
\end{equation}
For  estimation of the term $||E - E_\tau^I||_{L^2(\Omega_t)}$
 we are using standard interpolation estimate  with the interpolation 
 constant $C_I$
\begin{equation}\label{5.39}
||E - E_\tau^I||_{L^2(\Omega_t)} \leq
C_I \left \| \tau E \right \|_{H^1(\Omega_t)}.
\end{equation}

 We can estimate  $ || \tau ~E||_{H^1(\Omega_t)}$   in \eqref{5.39} as
\begin{equation}\label{5.40}
\begin{split}
  || \tau ~E||_{H^1(\Omega_t)}  &\leq \sum_{J_k} || \tau_k E||_{H^1(J_k)} =
  \sum_{J_k}   \left \| \left (E   +  \frac{ \partial E}{ \partial t} \right) \tau_k \right\|_{L_2(J_k)} \\
  &\leq  \sum_{J_k} \left( || E_\tau  \tau_k ||_{L_2(J_k)} +
  \left \Vert \frac{[E_\tau]}{\tau_k} \tau_k \right \Vert_{L_2(J_k)} \right) \\
  &\leq  || \tau E_\tau ||_{L_2(\Omega_t)} +
  \sum_{J_k} \left \Vert [E_\tau] \right \Vert_{L_2(J_k)}.
\end{split}
\end{equation}
Here,   $[E_\tau]$  denote the jump of the function $E_\tau$  over the time
intervals $[t_{k-1}, t_k]$ and $[t_k, t_{k+1}]$ which we define as
\begin{equation*}
[E_{\tau}] = E_{\tau}^+ - E_{\tau}^-,
\end{equation*}
where functions $E_{\tau}^-, E_{\tau}^+$ are computed on two neighboring time intervals $[t_{k-1}, t_k]$
and $[t_k, t_{k+1}]$, respectively.

Now we substitute the estimate \eqref{5.40} into \eqref{5.38} to get
 final a posteriori error estimate
\begin{equation}\label{theorem2_3}
 ||J(E) -J(E_{\tau}) ||_{L^2(\Omega_t)}  \leq C_I  \left\| J^{\prime }(E_{\tau})\right\| _{L^2(\Omega_t)} \left (|| \tau E_\tau ||_{L_2(\Omega_t)} +
  \sum_{J_k} \left \Vert [E_\tau] \right \Vert_{L_2(J_k)} \right) ~ \forall E_\tau \in  V_\tau^E.
\end{equation}

$\square $

\subsection{A posteriori error estimate of the minimizer on  refined meshes}

\label{sec:adaptrelax}

In the
Theorem 2  and Theorem 3
we derive  two  different a posteriori  error estimates for
 the error $ \| E - E_\tau \|_{L^2(\Omega_t)} $.
 Proof of the Theorem 2 follows from the proof of Theorem 5.1 of \cite{KB}.

\textbf{Theorem 2}

 \emph{
Let   $E \in L^2(\Omega_t)$ be a minimizer of the Tikonov's functional 
  \eqref{Tikh_f0} and
   $E_\tau \in V_\tau^E$
   be a finite element approximation of the minimizer on the finite element mesh $J_\tau$.
Then there exists a
  Lipschitz constant   $D = const.>0$ defined by
\begin{equation}
  \left\| J^{\prime }(E_1) - J^{\prime }(E_2) \right\|_{L_2(\Omega_t)} 
  \leq D \left\| E_1 - E_2\right\|_{L_2(\Omega_t)}  ,\forall E_1, E_2 \in L_2(\Omega_t),
  \label{2.10}
\end{equation}
  and the interpolation constant $C_I$  independent on the mesh size $\tau$
   such that the following
  a posteriori error estimate for the minimizer   $E $ holds  true  for $\gamma > 0$
\begin{equation}
|| E  - E_\tau ||_{L_2(\Omega_t)} \leq  \frac{D}{\gamma } C_I
 \left (|| \tau E_\tau ||_{L_2(\Omega_t)} +
  \sum_{J_k} \left \Vert [E_\tau] \right \Vert_{L_2(J_k)} \right)
 ~ \forall E_\tau \in V_\tau^E.
 \label{theorem1}
\end{equation}
 }

\textbf{Proof.}

 Let $E_\tau$ be the minimizer of the Tikhonov functional
(\ref{Tikh_f0}). The existence and uniqueness of this minimizer is
guaranteed by   conditions \eqref{6.20} and thus follows from Theorem 1.9.1.2 of \cite{BookBK}. Using this theorem we can conclude that the functional
(\ref{Tikh_f0}) is strongly convex on the space $L_2(\Omega_t)$ with the strong
convexity constant $\gamma$. This  yields that
\begin{equation}
\gamma \left\Vert E - E_\tau \right\Vert_{L_2(\Omega_t)} ^{2}\leq \| \left(
J^{\prime }\left( E \right)  - J^{\prime }\left( E_\tau \right) , E -
 E_\tau \right)\|_{L_2(\Omega_t)} .  \label{4.222}
\end{equation}
Here, $J^{\prime }(E_\tau), J^{\prime }\left(E\right)$ are the Fr\'{e}chet derivatives of the functional
(\ref{Tikh_f0})  given by the formula \eqref{Frder2}  for
  $E_\tau, E$, respectively.

Since $E$ is the minimizer of the Tikhonov functional
(\ref{Tikh_f0}) then
\begin{equation*}
\left( J^{\prime }\left(E \right), v \right) =0,\text{ }
\forall v \in L_2(\Omega_t).
\end{equation*}

Using again  the  splitting \eqref{splitting} in \eqref{4.222}
together with  the Galerkin orthogonality principle  for all
$E_\tau,  E_\tau^I  \in V_\tau^E$
\begin{equation}
  \left( J^{\prime }\left( E \right) - J^{\prime}\left(E_\tau\right),
  E_\tau^I -  E_\tau \right) =0
\label{4.223}
\end{equation}
  we get
\begin{equation}
  \gamma \left\Vert E - E_\tau \right\Vert_{L_2(\Omega_t)}^{2}\leq
  \| \left( J^{\prime }\left( E \right)  - J^{\prime }\left(E_\tau \right),
  E - E_\tau^I \right) \|_{L_2(\Omega_t)}  .  \label{4.224}
\end{equation}

We can estimate the
right hand side of (\ref{4.224})
 using (\ref{2.10}) as
\begin{equation*}
\|  \left( J^{\prime }\left( E \right) - J^{\prime }\left( E_\tau \right),
 E-   E_\tau^I  \right) \|_{L_2(\Omega_t)} 
    \leq D || E - E_\tau  ||_{L_2(\Omega_t)} || E - E_\tau^I   ||_{L_2(\Omega_t)}.
\end{equation*}
Substituting above equation into (\ref{4.224})  we obtain
\begin{equation}
  || E -  E_\tau  ||_{L_2(\Omega_t)} \leq \frac{D}{\gamma}
  ||  E -  E_\tau^I  ||_{L_2(\Omega_t)} .  \label{theorem1_1}
\end{equation}
Using again the interpolation property \eqref{5.39}
 we obtain a posteriori error estimate for the regularized solution :
\begin{equation}\label{5.48}
  || E - E_\tau  ||_{L_2(\Omega_t)} \leq \frac{D}{\gamma }
  || E -  E_\tau^I   ||_{L_2(\Omega_t)} \leq   \frac{D}{\gamma}
  C_I || \tau ~E||_{H^1(\Omega_t)}.
\end{equation}.

We can estimate the term $ || \tau ~E||_{H^1(\Omega_t)}$   in the right hand side of \eqref{5.48} similarly as in   \eqref{5.40}.
Substituting the estimate \eqref{5.40}  into the right hand side of (\ref{5.48})
we get
\begin{equation*}
  || E - E_\tau ||_{L_2(\Omega_t)} \leq
  \frac{D}{\gamma } C_I \left (|| \tau E_\tau ||_{L_2(\Omega_t)} +
  \sum_{J_k} \left \Vert [E_\tau] \right \Vert_{L_2(J_k)} \right)
 ~ \forall E_\tau \in V_\tau^E.
\end{equation*}

$\square $

\textbf{Theorem 3}

 \emph{Let $E_\tau \in V_\tau^E$
   be a finite element approximation
 of the minimizer
   $E \in L^2(\Omega_t)$ 
on the finite element mesh $J_\tau$.
Then there exists
  an interpolation constant $C_I$  independent on the mesh function $\tau$
   such that the following
  a posteriori error estimate   holds  for $\gamma > 0$:
\begin{equation}
|| E - E_\tau   ||_{L_2(\Omega_t)} \leq 
\frac{\|R(E_\tau)\|_{L_2(\Omega_t)}  }{\gamma}  ~~~~~ \forall E_\tau \in V_\tau^E,  \label{theorem3}
\end{equation}
where $R(E_\tau)$ is the residual defined as
\begin{equation}\label{resid}
  R(E_\tau)  = \gamma  (E_\tau - E^{0}) +
 d  {\lambda_2}_{\tau} {u_2}_\tau.
\end{equation}
}

\textbf{Proof.}

Let again  $E_\tau$ be the minimizer of the Tikhonov functional
(\ref{Tikh_f0}). Strong convexity  of the functional (\ref{Tikh_f0})  on the space $L_2(\Omega_t)$ implies that
\begin{equation}
  \gamma \left\Vert E -  E_\tau \right\Vert_{L_2(\Omega_t)} ^{2}\leq
  \| \left(
J^{\prime }\left(  E \right)  - J^{\prime }\left( E_\tau \right) , E -
 E_\tau \right) \|_{L_2(\Omega_t)}.  \label{T31}
\end{equation}

Applying   splitting \eqref{splitting} to  \eqref{T31}
we obtain \eqref{4.224} where the term   $J'(E_\tau)$
can be written via
\eqref{Frder2}. More precisely, when $u(t), \lambda(t)$ are exact  solutions of the forward and adjoint problems, respectively,  we  have  for $E_\tau$:
\begin{equation*}
L(v(E_\tau)) = J(E_\tau),
  \end{equation*}
and thus,   one can write
\begin{equation}\label{T31n}
J'(E_\tau) = L'(E_\tau) = \gamma  (E_\tau - E^{0}) +
 d  {\lambda_2} {u_2}.
  \end{equation}

We can  write estimate  \eqref{T31}  as
\begin{equation*}
  \gamma \left\Vert E -  E_\tau \right\Vert_{L_2(\Omega_t)} \leq
  \|  J^{\prime }\left(  E \right)  - J^{\prime }\left( E_\tau \right)\|_{L_2(\Omega_t)}
\end{equation*}

Then  noting that $J'(E)=0$  and
  using  \eqref{T31n}
 for approximations ${u_2}_\tau, {\lambda_2}_\tau$  of  ${u_2}, {\lambda_2}$,
 respectively,
we get
\begin{equation}\label{T32}
  ||  E - E_\tau  ||_{L_2(\Omega_t)} \leq \frac{\|R(E_\tau)\|_{L_2(\Omega_t)}}{\gamma}
\end{equation}
  where $R(E_\tau)$  is the residual defined as in \eqref{resid}.

 $\square $


\begin{algorithm}[hbt!]
  \centering
  \caption{Conjugate Gradient  Algorithm  (CGA) \label{alg:cga}}
  \begin{algorithmic}[1]
  
    \STATE Choose time partition $J_{\tau}$ of the time
  interval $(0, T)$. Choose initial value
of the regularization parameter $\gamma_0$, step size $\beta_0$  and
the initial approximation
  ${E_\tau}^{0}$. Compute approximations  of ${E_\tau}^{m}$,   $m = 1, ..., M$ in the following steps.

  \STATE  Use the Newton's method on the time partition
  $ J_{\tau}$ to
compute approximate solutions $u_\tau^m = u_\tau\left({E_\tau}^{m}, t\right)
  , \lambda_\tau^m = \lambda_\tau\left({E_\tau}^{m}, t \right) $ of the state (\ref{state_forw_new2})
  and adjoint (\ref{adjoint}) problems via iterative schemes \eqref{4.31}, \eqref{4.36}, respectively.

\STATE Compute gradient $G^{m}(t_i)$  by (\ref{gradient2}).
\STATE Update the unknown parameter $E:= E_\tau^{m+1}$   as
\begin{equation*}
\begin{split}
E_\tau^{m+1}(t_i) &=  E_\tau^{m}(t_i)  + r^m d^m(t_i)
\end{split}
\end{equation*}
with
\begin{equation*}
\begin{split}
 d^m(t_i)&=  -G^m(t_i)  + \beta^m  d^{m-1}(t_i),
\end{split}
\end{equation*}
and
\begin{equation}
r^m = -\frac{(G^m, d^m)}{\gamma_m \| d^m\|^2},
\end{equation}
\begin{equation*}
\begin{split}
 \beta^m &= \frac{|| G^m(t_i)||^2}{|| G^{m-1}(t_i)||^2},
\end{split}
\end{equation*}
where $d^0(t_i)= -G^0(t_i)$  and $G^{m}(t_i)$ is the gradient  vector  which is computed by (\ref{gradient2}) in time moments $t_i$.

\STATE  Update the regularization   parameter $\gamma_m$ for any $p \in (0,1)$   via iterative rule of \cite{BKok}

\begin{equation}
\gamma_m = \frac{\gamma_0}{(m+1)^p}
\end{equation}

\STATE Choose the tolerance  $0 <\theta < 1$ and stop computing the functions
  $E_\tau^{m}$ if either
  $||G^m||_{L_2(\Omega_t)} \leq \theta$, or  norm of the computed gradient
  $||G^m||_{L_2(\Omega_t)} $ abruptly grow, or  relative norms  of the computed parameter
$\frac{  \| E_\tau^{m} - E_\tau^{m-1} \|}{\| E_\tau^m\| }$ 
  are stabilized. Otherwise, set $m:=m+1$ and go to Step 2.
  \end{algorithmic}
\end{algorithm}

\section{Optimization algorithms}

\label{sec:fem_IP}

This sections presents two different optimization
 algorithms for solution of PIP:  conjugate gradient algorithm  (CGA)
 and time-adaptive conjugate gradient algorithm  (ACGA).
The CGA algorithm is a standard one and can be used on a mesh with the
equidistant time step $\tau_k$ on every time interval $J_k$.  The ACGA is used for minimization of
the functional \eqref{Tikh_f0} on a locally adaptivelly refined meshes in time
and thus, the computational mesh may have different mesh sizes $\tau_k$
 on every time interval $J_k$ - see   definition of  partition $J_\tau$ in section \ref{sec:fem}.

 Let us denote the nodal value of the gradient at every iteration $m$ of the CGA at the observation points
$\{t_i\}$ by $G^{m}(t_i)$ and compute it accordingly to
\eqref{Frder2} as
\begin{equation}
   G^{m}(t_i) = \gamma   ( E_\tau^m(t_i) -  E_\tau^{0}(t_i)) + d(t_i) ~{ \lambda_3}_\tau^m(t_i) {u_2}_\tau^m(t_i).
\label{gradient2}
\end{equation}
The  approximations of solutions ${u_2}_\tau^{m} $and ${\lambda_{3}}_\tau^m$ in the expression for
gradient \eqref{gradient2}
are
obtained  via solution of the forward and adjoint problems, respectively,  via Newton's method
taking $E:={E_\tau^{m}}$.
Algorithm 1 describes explicitly all steps of CGA algorithm.

In the formulation of an adaptive algorithm ACGA we have used  Theorem 3  
for the estimation of the error $e^k = \| E_\tau^{M,k} - E \|_{L_2(\Omega_t)}$ computed on the locally $k$ times  refined meshes
$J_\tau^k$.
We  choose
 the tolerance  $0 <\theta < 1$  and run ACGA  algorithm until
\begin{equation*}
e^k = \| E_\tau^{M,k} - E \|_{L_2(\Omega_t)} \leq \theta,
\end{equation*}
where   $e^k $ is estimated  via a posteriori  error estimate \eqref{theorem3} of the Theorem 3.
Here,  $E_\tau^{M,k}$ is the computed parameter  obtained at the final  iteration   $M$ of the CGA algorithm on k times refined mesh $J_\tau^k$.
For the time-mesh refinements we use the following   criterion.

\begin{algorithm}[hbt!]
  \centering
  \caption{Adaptive Conjugate Gradient Algorithm (ACGA) \label{alg:acga}}
  \begin{algorithmic}[1]
  
\STATE Choose initial time partition $J_{\tau}^0$ of the time
  interval $(0, T)$. Choose initial value
 of the regularization parameter $\gamma_0^0$, step size in CGA $\beta_0^0$,
and the initial approximation
  ${E_\tau}^{0,0}$. Initialize data $g_1^0 := g_1, g_2^0 = g_2$.
 For mesh refinements $k=0,...,N$  perform  following steps:


\STATE Compute $E_\tau^{M,k}$ on the time mesh $ J_{\tau}^k$ using steps
2-6 in CGA algorithm. Here,  $E_\tau^{M,k}$ is the computed  parameter obtained
in the CGA algorithm at the final
 iteration $M$ on the  $k$ times refined mesh  $ J_{\tau}^k$.

\STATE  Refine the time mesh $ J_{\tau}^k$ at all points where
\begin{equation}
\frac{\left \vert R(E_\tau^{M,k})(t) \right \vert}{\gamma}
  \geq \tilde{\beta_k} \max_{\Omega_t}
  \frac{ \left\vert R(E_\tau^{M,k}) (t) \right\vert }{\gamma}.  \label{6.70}
\end{equation}
Here, $ \tilde{\beta_k} \in (0,1)$ is the  mesh refinement parameter chosen by the user.

\STATE  Construct a new  time
partition $J_{\tau}^{k+1}$ of the time interval $\left( 0,T\right) $.

\STATE  Interpolate ( alternative is chosen by the user):

\begin{itemize}

\item or the initial approximation ${E_\tau}^{0,k}$ from the
 previous  time partition $J_{\tau}^{k}$ to the new  time partition $J_{\tau}^{k+1}$.

\item or   computed ${E_\tau}^{M,k}$ from the
 previous  time partition  $J_{\tau}^{k}$ to the new  time partition $J_{\tau}^{k+1}$.

\end{itemize}

\STATE  Interpolate data $g_1^k, g_2^k$  from the
 previous  time partition $J_{\tau}^{k}$  to the new  time partition $J_{\tau}^{k}$.
 
\STATE  Stop    computations if norms  of the gradients
$||G({E_\tau}^{M,k})||_{L_2(\Omega_t)} $  either increase or stabilize, compared with the
  previous time partition. Otherwise set $k=k+1$ and go to Step 2.

  \end{algorithmic}
\end{algorithm}

 \textbf{The  Refinements Criterion}

Refine the time-mesh $ J_{\tau}$ in neighborhoods of those time-mesh
  points  $t\in {\Omega_t}$  where the  weighted  residual $\left\vert
 R(E_{\tau})(t) \right\vert$  defined in \eqref{resid}
 attains its maximal values:
\begin{equation*}
\frac{ \left\vert R(E_\tau)(t) \right\vert}{\gamma}
\geq \tilde{\beta}  \max_{\Omega_t} \frac{\left\vert R(E_\tau)(t) \right\vert }{\gamma}.
\end{equation*}
Here, $\tilde{\beta} \in \left(0,1\right) $ is   the tolerance number chosen by user.
See discussion in \cite{BErG}  how this number can be chosen in optimal way.
Algorithm 2 uses the above mesh refinement recommendation and implements ACGA algorithm.

\section{Numerical studies}

\label{sec:numex}

This section presents several numerical results  which demonstrate
performance  of the time-adaptive reconstruction of
 the immune response function $E(t)$ in PIP using ACGA
algorithm.
  Numerical tests are performed in Matlab R2023b.
  The code is available for download at \cite{source}.

 In all computations we have used clinical data for the virus population function $u_3(t) \approx  g_2(t)$ and
 for  the total number of uninfected and infected T-cells  $u_1(t) + u_2(t) \approx g_1(t)$
  for four patients presented in the Table  3.
The initial time-mesh $J_\tau^0$ of the observation interval $[0,T] = [0,363]$ is generated  such that the time step is $\tau=1$.
We note that only  data at 8 points are  given in the Table 3.  Since we are running CGA  algorithm on the mesh $J_\tau^0$ consisting of
  364 points we first use  linear interpolation to interpolate data given in the Table 3 to the data on the time mesh  $J_\tau^0$.
  Next, we use interpolated data   in all numerical tests for solution of PIP.


The value of the mesh refinement parameter $\tilde{\beta_k}$ in ACGA is chosen such that it allows local
refinements  of the time-mesh $J_\tau^k$ where $k$ is number of mesh refinements.
All computations presented below are  performed for the mesh refinement
 parameter $\tilde{\beta_k}   =0.875$  in \eqref{6.70}  for all mesh refinements $k$.

We compute also relative residuals  $R_1, R_2$  between computed  functions ${u_3}_\tau^{M,k}$,  $ {u_1}_\tau^{M,k} + {u_2}_\tau^{M,k} $  and interpolated
measured  $V_\tau^k = {g_2}_\tau^k, \Sigma_\tau^k = {g_1}_\tau^k $ which are defined as
 \begin{equation}\label{residR}
 \begin{split}
R_1 &= \frac{  |  \log_{10} {u_3}_\tau^{M,k}  -  \log_{10} {g_2}_\tau^{M,k}    | }{  \rm nno},\\
R_2 &= \frac{ | \log_{10} ({u_1}_\tau^{M,k}  + {u_2}_\tau^{M,k})  -  \log_{10} {g_1}_\tau^k  |}{  \rm nno}
\end{split}
 \end{equation}

In all computations we choose parameters in the model system \eqref{forward_bio1} as in the Table 1.
 The initial data  is   defined as follows:
\begin{equation}\label{initdata}
  \begin{split}
u_1(0) &= 1125000; \\
u_2(0) &=0;\\
u_3(0) &=1. 
\end{split}
\end{equation}

Since CGA and ACGA algorithms  are locally convergent
    it is of vital importance that the initial guess $E^0(t)$  for approximation of the function $E(t)$
    satisfy  conditions
 \eqref{6.1} what means that $E^0$ is located in the close
neighborhood of the exact   function $E(t)$.
We model the initial guess  $E^0(t)$  for the function $E(t)$ with  inclusion of such called cytotoxic T lymphocyte
 (CTL) response in  $E^0(t)$,
   as it was  proposed in \cite{Stafford}.

 In  all our tests
 we have modelled  effect of the CTL response for both  functions $d(t), E^0(t)$ as follows:
 \begin{equation}\label{CTLTest4}
  \begin{split}
  d(t) &= d_0 +  d_1(t, V),\\
  d_1(t, u_3) &=
  \begin{cases}
     0, t < t_1,\\
     f(t) \log_{10} V, t \geq t_1,
  \end{cases}
  \\
 E^0(t) &= E_0^0 +  d_1(t, V), \\
  f(t) &= \frac{\beta_{CTL}}{1 + \kappa e^{(-(t-t_1)/ \delta T_1)}} -
  \frac{\beta_{CTL}}{1 + \kappa e^{(-(t-t_2)/ \delta T_2)}}, \\
  \kappa &= 1 + 10^5\beta_{CTL}.
  \end{split}
  \end{equation}

 Parameter values  for all parameters  in \eqref{CTLTest4}  are chosen differently for all patients such that they give as good as possible  fitting to data.
 Detailed values of parameters are given in the next subsections.

\subsection{Results for Patient 1}

\begin{figure}
\begin{center}
\begin{tabular}{ccc}
  {\includegraphics[scale=0.5, clip = true, trim = 0.0cm 0.0cm 0.0cm 0.0cm ]{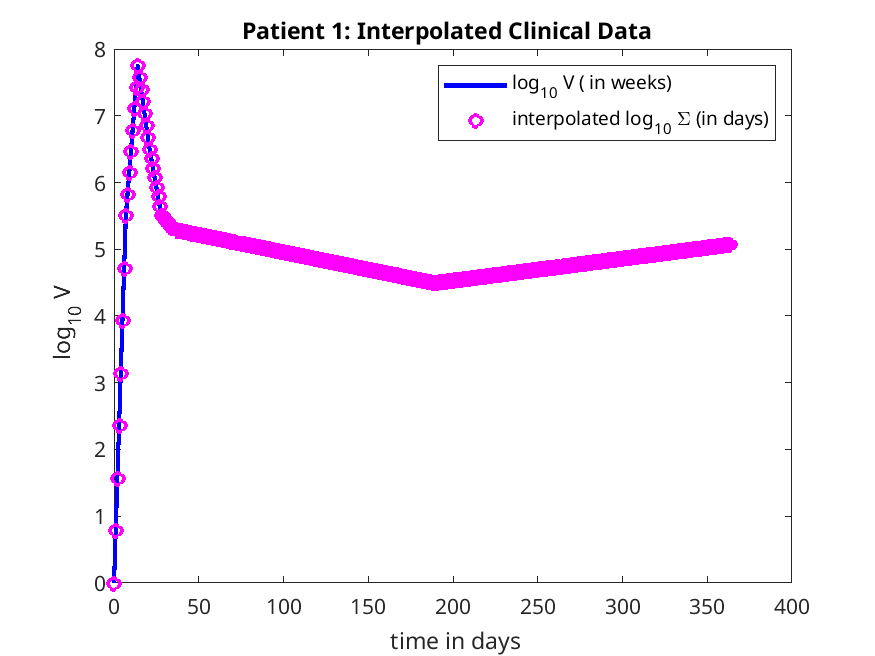}} &
  {\includegraphics[scale=0.5, clip = true, trim = 0.0cm 0.0cm 0.0cm 0.0cm]{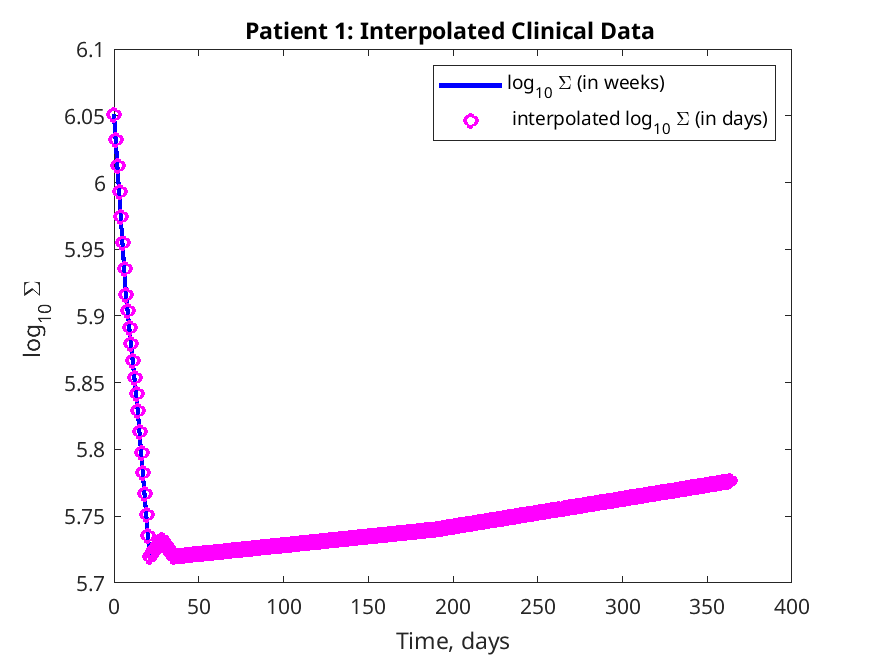}} \\
 $ V = g_2^0$ & $ \Sigma = g_1^0$  \\
\end{tabular}
\end{center}
\caption{Patient 1:  Interpolated clinical data by linear spline.}
 \label{fig:Pat1data}
 \end{figure}

In computations of boths, CGA and ACGA, algorithms for patient 1 we used interpolated clinical data  of the Table 2 presented in the Figure  \ref{fig:Pat1data}. 
  Next, we define the initial guess for the immune response function $E_0$ in optimization algorithms using \eqref{CTLTest4}.
For patient 1 the parameters  in \eqref{CTLTest4} are chosen as follows:
\begin{equation}\label{paramCTLPat1}
  \begin{split}
d_0 & = 0.26; \\
\beta_{CTL} &= 0.015;\\
\kappa  &=  1 + 10^5 \beta_{CTL}; \\
\delta T_1 &= 2.5; \\
\delta T_2 &= 5.0; \\
E_0^0 &= 1,\\
t_1 &= 1,\\
t_2 &=1.
\end{split}
\end{equation}

\begin{figure}
\begin{center}
\begin{tabular}{ccc}
  {\includegraphics[scale=0.35, clip = true, trim = 0.0cm 0.0cm 0.0cm 0.0cm ]{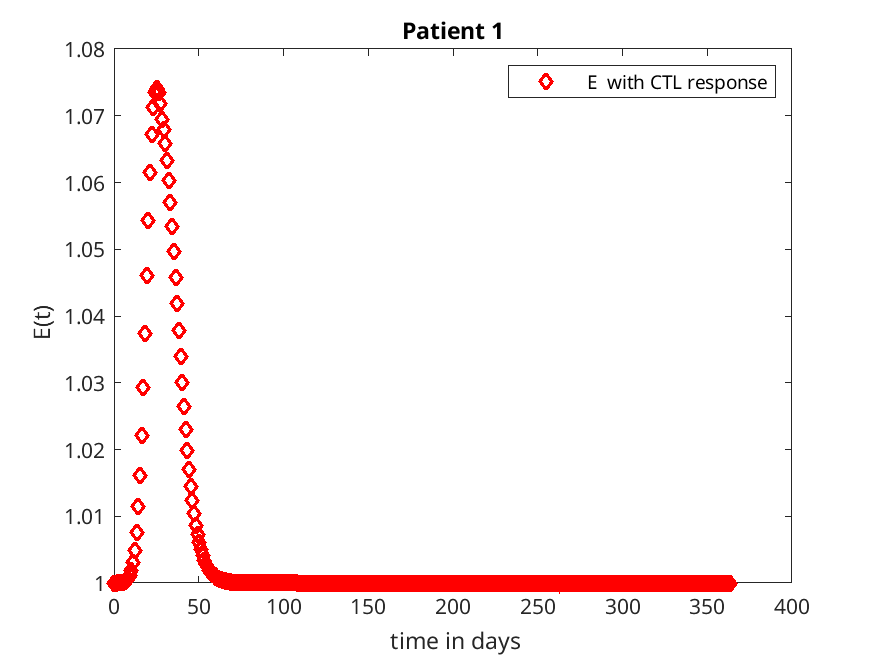}} &
  {\includegraphics[scale=0.35, clip = true, trim = 0.0cm 0.0cm 0.0cm 0.0cm]{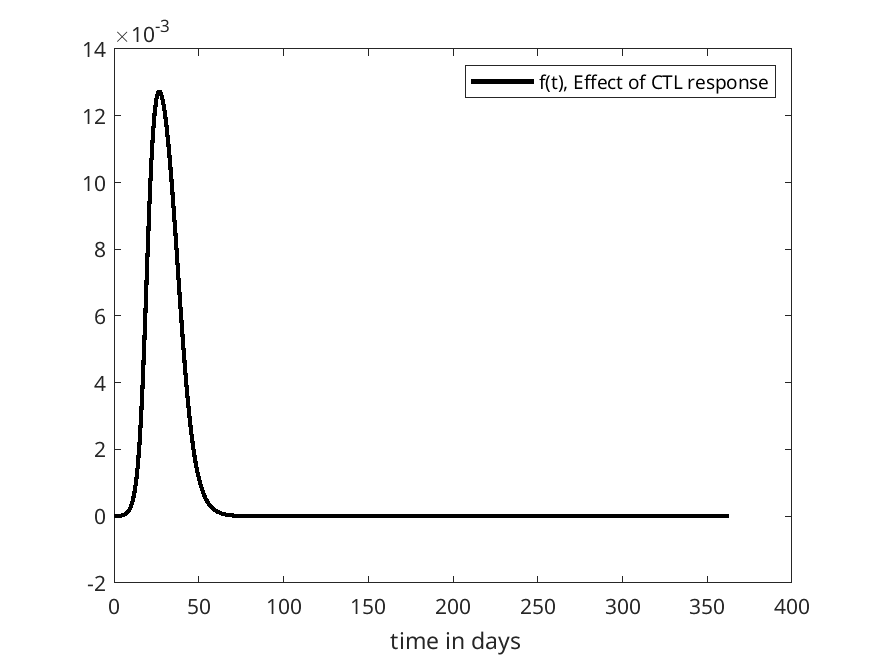}} &
    {\includegraphics[scale=0.35, clip = true, trim = 0.0cm 0.0cm 0.0cm 0.0cm ]{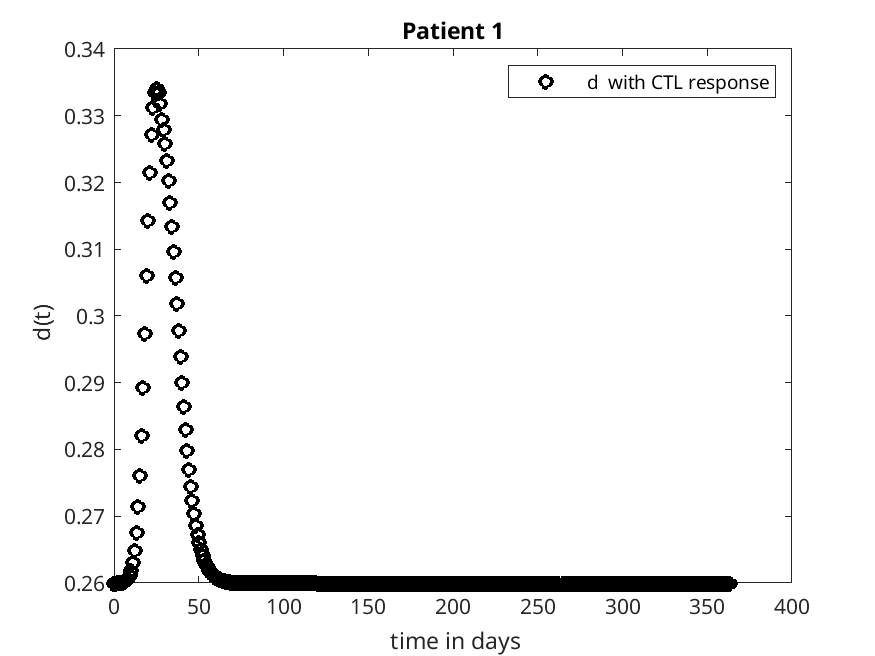}} \\
 $ E^0(t)$ & $f(t)$  & $d(t)$ \\
\end{tabular}
\end{center}
\caption{Patient 1: functions $E^0(t), f(t), d(t)$ with CTL response modelled as in \eqref{CTLTest4}  using parameters defined in \eqref{paramCTLPat1}.}
 \label{fig:Pat1CTLfunc}
 \end{figure}

\begin{figure}
\begin{center}
\begin{tabular}{cc}
  {\includegraphics[scale=0.5, clip = true, trim = 0.0cm 0.0cm 0.0cm 0.0cm ]{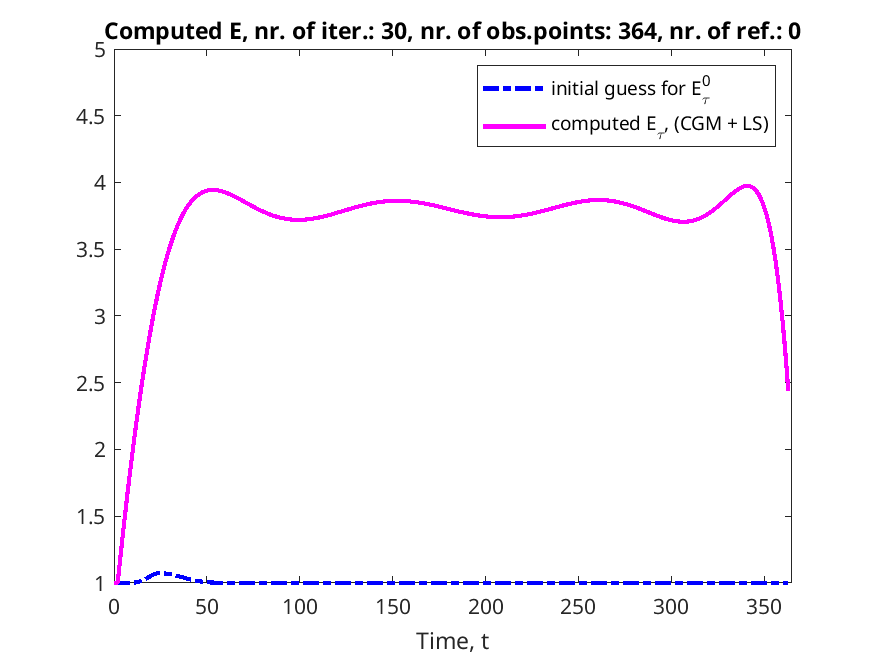}} &
  {\includegraphics[scale=0.5, clip = true, trim = 0.0cm 0.0cm 0.0cm 0.0cm]{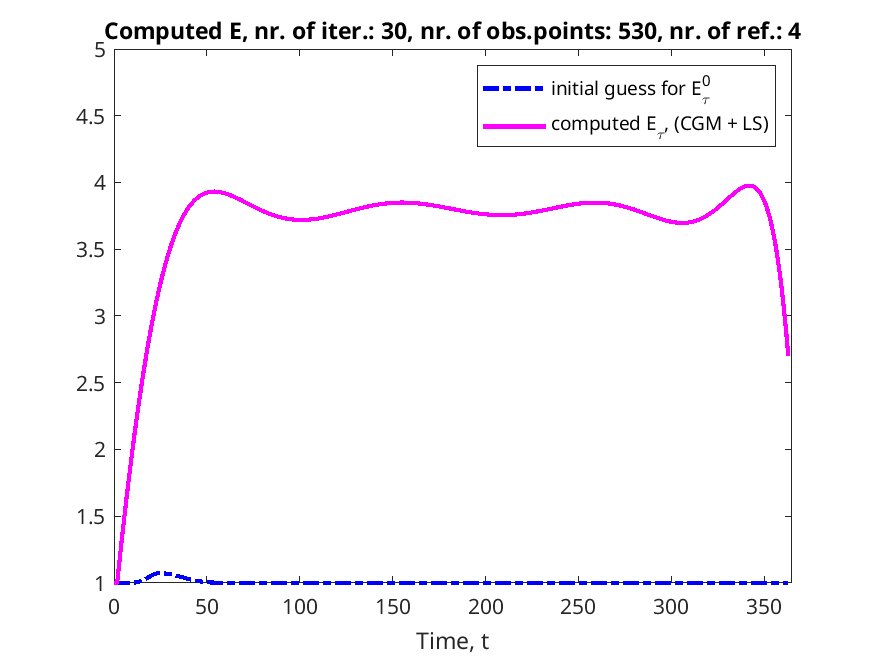}} \\
  $k=0$  & $k=4$ 
\end{tabular}
\end{center}
\caption{Patient 1:  results of reconstruction of the function $E(t)$ on   $k$ times adaptivelly refined meshes $J_\tau^k$
 in ACGA algorithm.
 Computations are  performed for the mesh refinement parameter $\tilde{\beta_k}
  =0.875$  for all mesh refinements $k$.}
 \label{fig:Pat1RecE}
 \end{figure}

Figures \ref{fig:Pat1CTLfunc}  show  time-dependent  functions for initial guess $E^0(t)$ as well as $ f(t), d(t)$
 with CTL response  for parameters defined in \eqref{paramCTLPat1}.
These functions are used  now in optimization algorithms  for reconstruction of the immune response function $E(t)$.
 Figures \ref{fig:Pat1RecE} show
 computed functions $E_\tau^k$
 on   $k$ times adaptivelly refined meshes $J_\tau^k$
 in ACGA algorithm.
 
\begin{figure}
\begin{center}
\begin{tabular}{cc}
  {\includegraphics[scale=0.5, clip = true, trim = 0.0cm 0.0cm 0.0cm 0.0cm ]{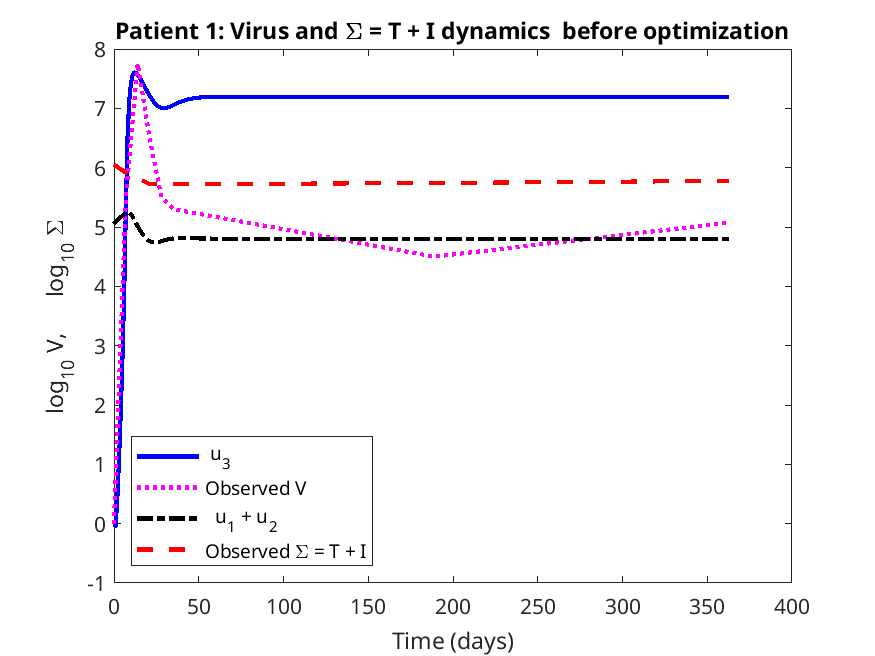}} &
  {\includegraphics[scale=0.5, clip = true, trim = 0.0cm 0.0cm 0.0cm 0.0cm]{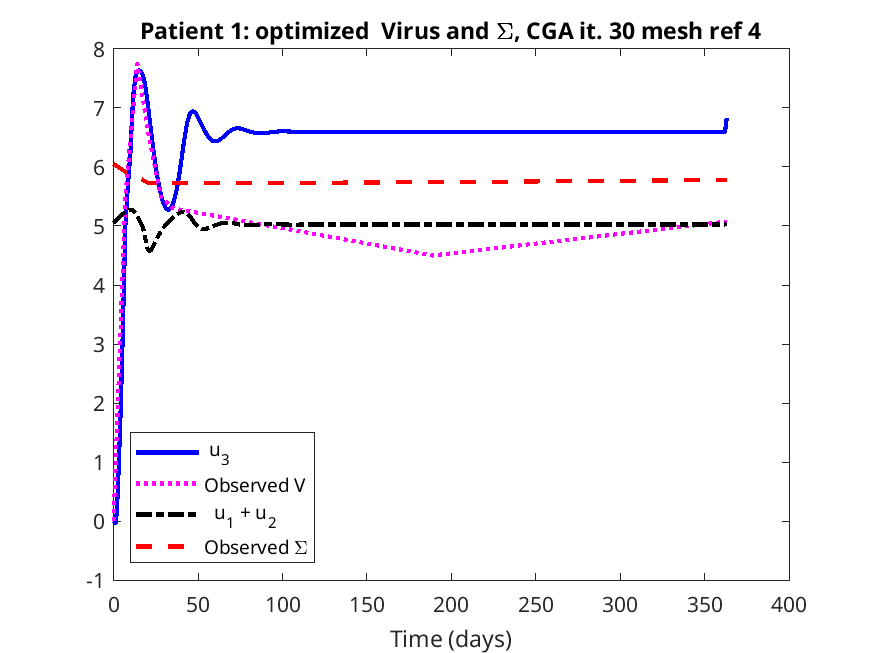}} \\
  $k=0$ & $k=4$ 
\end{tabular}
\end{center}
\caption{ Patient 1: Dynamics of the computed virus function $V_\tau^k$ before and after optimization corresponding to the
   computed $E_\tau^k$ on   $k$ times adaptivelly refined meshes versus interpolated clinical data ${g_1^0}_\tau,{g_2^0}_\tau$.}
 \label{fig:Pat1virus}
 \end{figure}

Figures \ref{fig:Pat1virus}  show  time-dependent behaviour  of the computed virus function $V_\tau^k, k=0,4$ before and after applying ACGA.
Here,  the computed virus function $V_\tau^k$ corresponds to the
   computed $E_\tau^k$ presented in Figure \ref{fig:Pat1RecE}.
Results are compared with  clinical data  for the observed virus function  $V = g_2^0$  and observed total number of  the uninfected
and infected cells $\Sigma = g_1^0$.
Figure \ref{fig:Pat1virus}-b) shows computed virus function on the four times locally refined time mesh  $J_\tau^4$ . Computed optimized
virus function on  $J_\tau^k, k=0,1,2,3$ has similar behaviour and is not presented here.

\begin{figure}
\begin{center}
\begin{tabular}{cc}
  {\includegraphics[scale=0.5, clip = true, trim = 0.0cm 0.0cm 0.0cm 0.0cm ]{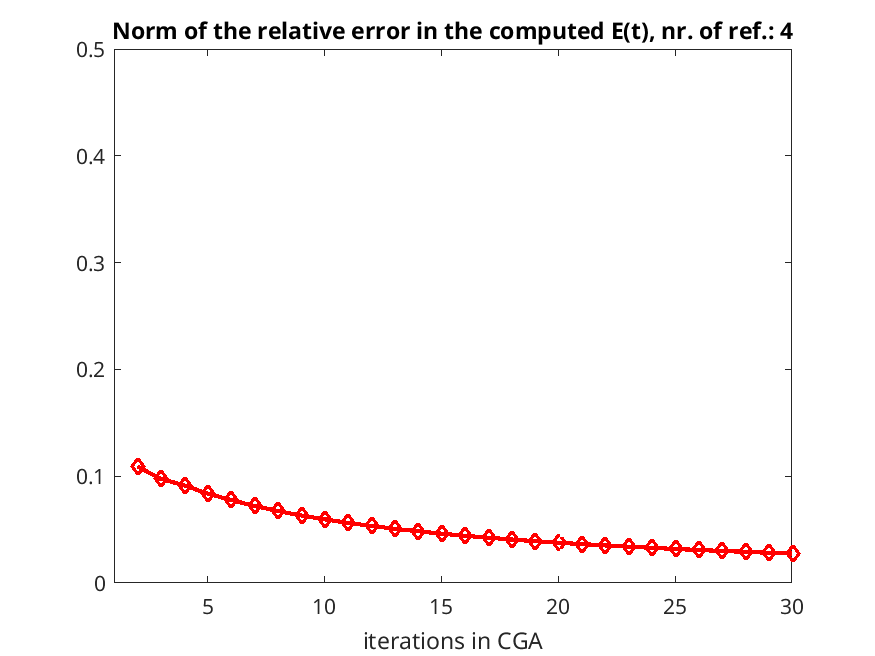}} &
  {\includegraphics[scale=0.5, clip = true, trim = 0.0cm 0.0cm 0.0cm 0.0cm]{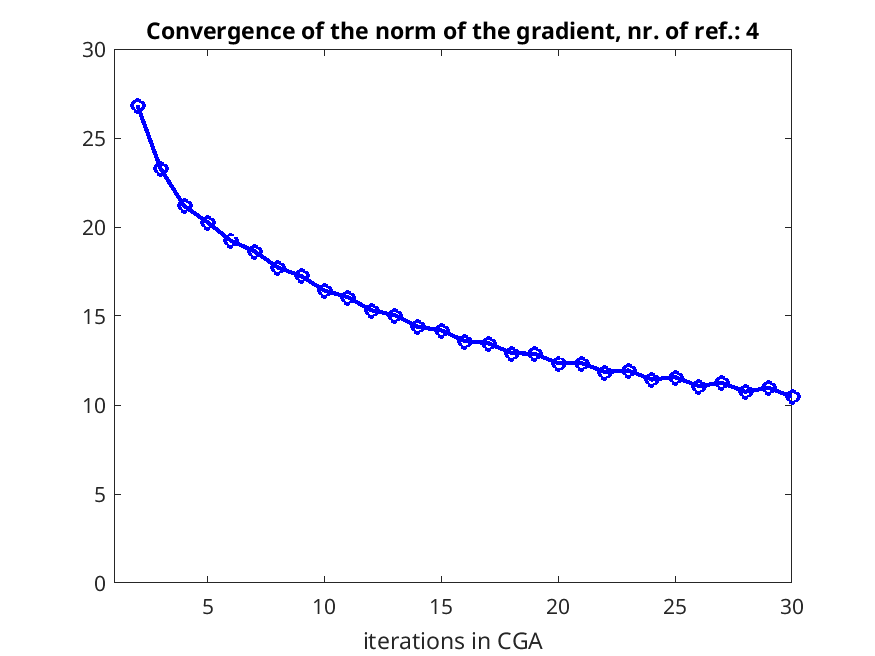}} \\
  $\frac{  \| E_\tau^{m} - E_\tau^{m-1} \|_{L_2(\Omega_t)}}{\| E_\tau^m\|_{L_2(\Omega_t)}}$   &  $\|G^m(t) \|_{L_2(\Omega_t)} $ 
\end{tabular}
\end{center}
\caption{Patient 1:   computed relative norms  $\frac{  \| E_\tau^{m} - E_\tau^{m-1} \|_{L_2(\Omega_t)}}{\| E_\tau^m\|_{L_2(\Omega_t)}}$
and $\|G^m(t) \|_{L_2(\Omega_t)}$  on the mesh $J_\tau^4$.}
 \label{fig:Pat1Residuals}
 \end{figure}

\begin{figure}
\begin{center}
\begin{tabular}{cc}
  {\includegraphics[scale=0.5, clip = true, trim = 0.0cm 0.0cm 0.0cm 0.0cm ]{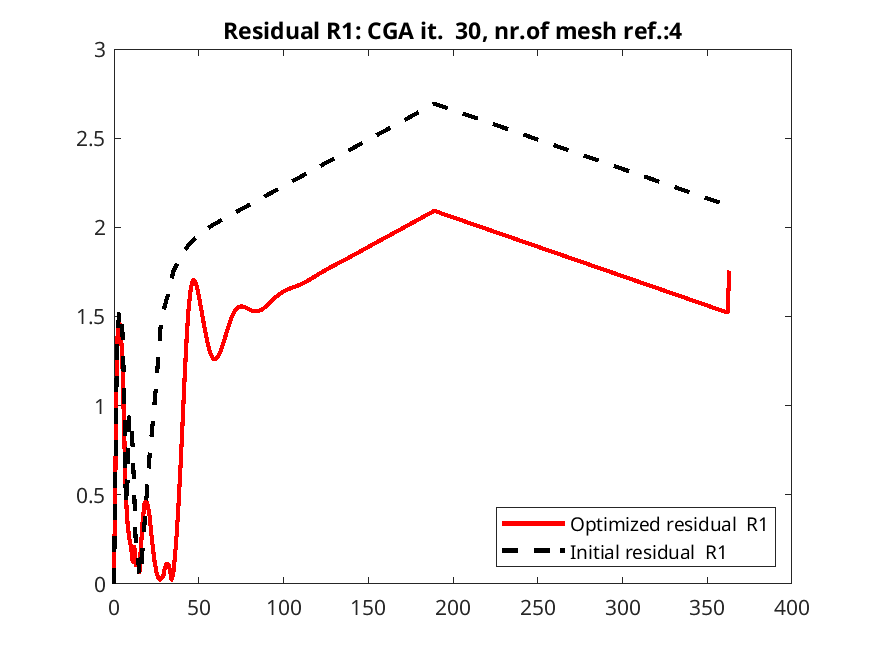}} &
  {\includegraphics[scale=0.5, clip = true, trim = 0.0cm 0.0cm 0.0cm 0.0cm]{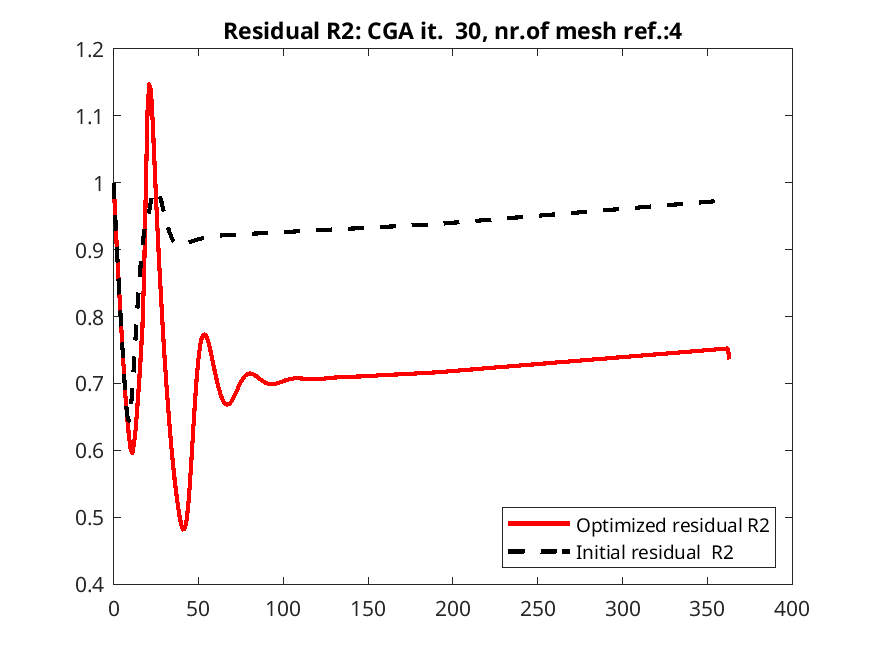}} \\
   a) $  R_1 $  & b) $ R_2  $  \\
 {\includegraphics[scale=0.5, clip = true, trim = 0.0cm 0.0cm 0.0cm 0.0cm]{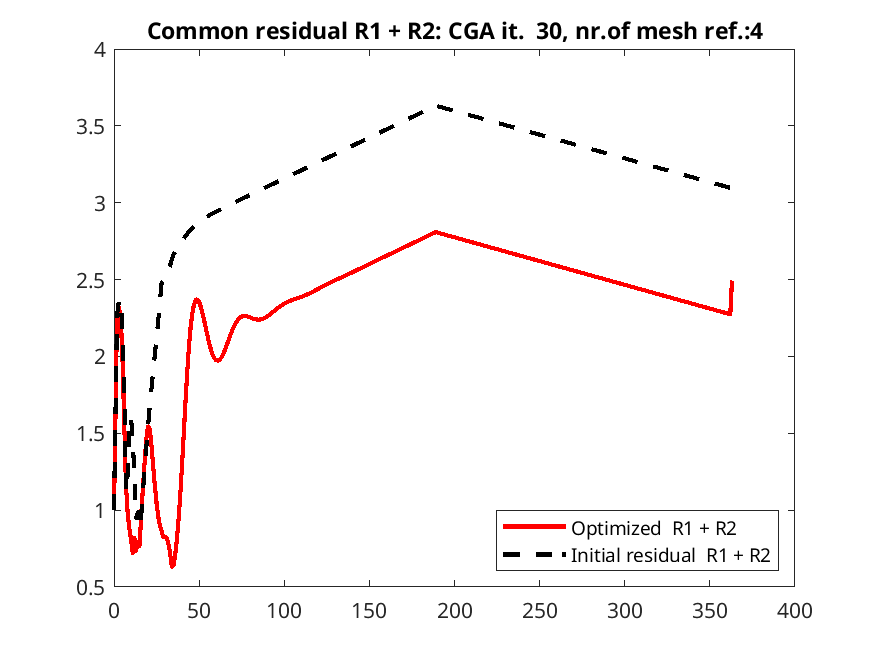}} &
  {\includegraphics[scale=0.5, clip = true, trim = 0.0cm 0.0cm 0.0cm 0.0cm]{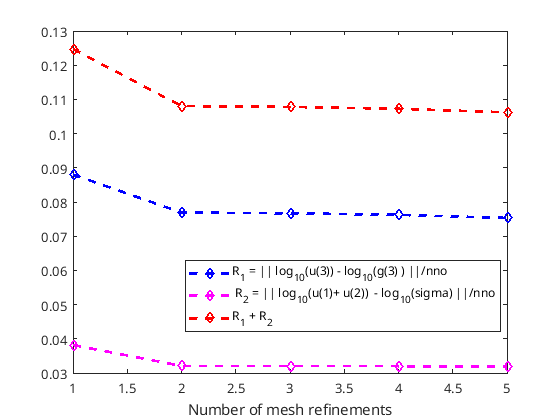}}  \\
 c)  $ R_1 + R_2  $ &  d) 
\end{tabular}
\end{center}
\caption{Patient 1: a), b), c)  computed residuals $R_1$ and $R_2$ on the mesh $J_\tau^4$; d) Comparison of related residuals
 $\|R_1\|, \|R_2\|$ on different refined meshes
$J_\tau^k, k = 0,...,4$. }
 \label{fig:Pat1R1R2}
 \end{figure}

Figures   \ref{fig:Pat1Residuals} show
stabilization of the computed relative norms  
$\frac{  \| E_\tau^{m} - E_\tau^{m-1} \|_{L_2(\Omega_t)}}{\| E_\tau^m\|_{L_2(\Omega_t)}}$ 
  and behaviour of norms of the computed gradient $\|G^m(t) \|_{L_2(\Omega_t)} $  for patient 1.
 Relative residuals $R_1$ and $R_2$ and $\|R_1\|_{L_2(\Omega_t)} , \|R_2\|_{L_2(\Omega_t)} $  which are  computed as in
 \eqref{residR}  on the locally refined meshes
 $J_\tau^k$ are presented in the Table 4 as well as in the Figures
 \ref{fig:Pat1R1R2}.
 All these norms are used in ACGA  in the stopping criterion.

\begin{table}[H]
\begin{center}
\begin{tabular}{cc}
{\small
\begin{tabular}{|c|c|c|c|}
  \hline
 k &\rm  nno & $\| R_1 \|$  & $\| R_2 \|$ \\
 \hline
 0 & 364 & 0.0880  & 0.0382 \\
 1 & 516 & 0.0769  & 0.0321 \\
 2 & 518 & 0.0767 &  0.0320\\
 3 & 522 & 0.0762  & 0.0320 \\
 4 & 530 & 0.0753  & 0.0320  \\
  \hline
\end{tabular}
}
&
{\small
\begin{tabular}{|c|c|c|c|}
  \hline
 $k$ & \rm  nno & $\| R_1 \|$  & $\| R_2 \|$ \\
 \hline
 0 &  364 & 0.1321    & 0.0412 \\
 1 &  410  & 0.1215    & 0.0389 \\
   &   &     &       \\
   &   &     &  \\
   &   &     &   \\
  \hline
\end{tabular}
}
\\
a)  Patient 1 & b)   Patient 2 \\ 
{\small
\begin{tabular}{|c|c|c|c|}
  \hline
 $k$ &\rm  nno  & $\| R_1 \|$  & $\| R_2\|$ \\
 \hline
 0 &  364 & 0.1643    & 0.0434 \\
 1 &  385  & 0.1560    & 0.0421 \\
 2  &  402 &  0.1495   &  0.0411     \\
   &   &     &  \\
   &   &     &   \\
  \hline
\end{tabular}
}
&
{\small
\begin{tabular}{|c|c|c|c|}
  \hline
 $k$ &\rm  nno  & $\| R_1 \|$  & $ \|R_2 \|$ \\
 \hline
 0 &  364 & 0.1802    & 0.0407 \\
 1 &  413  & 0.165    & 0.0384 \\
 2  &  478 &  0.1464   & 0.0358      \\
   &   &     &  \\
   &   &     &   \\
  \hline
\end{tabular}
}
\\
c) Patient 3 & d) Patient 4  \\
\end{tabular}
\caption{\textit{ Computed  norms of the residuals $R_1$ and $R_2$ defined by
 \eqref{residR}  on the meshes
 $J_\tau^k$. Here, $nno$ is the number of nodes in the time-mesh.}}
\label{tab:table5}
\end{center}
\end{table}

\subsection{Results for Patient 2}

\begin{figure}
\begin{center}
\begin{tabular}{ccc}
  {\includegraphics[scale=0.5, clip = true, trim = 0.0cm 0.0cm 0.0cm 0.0cm ]{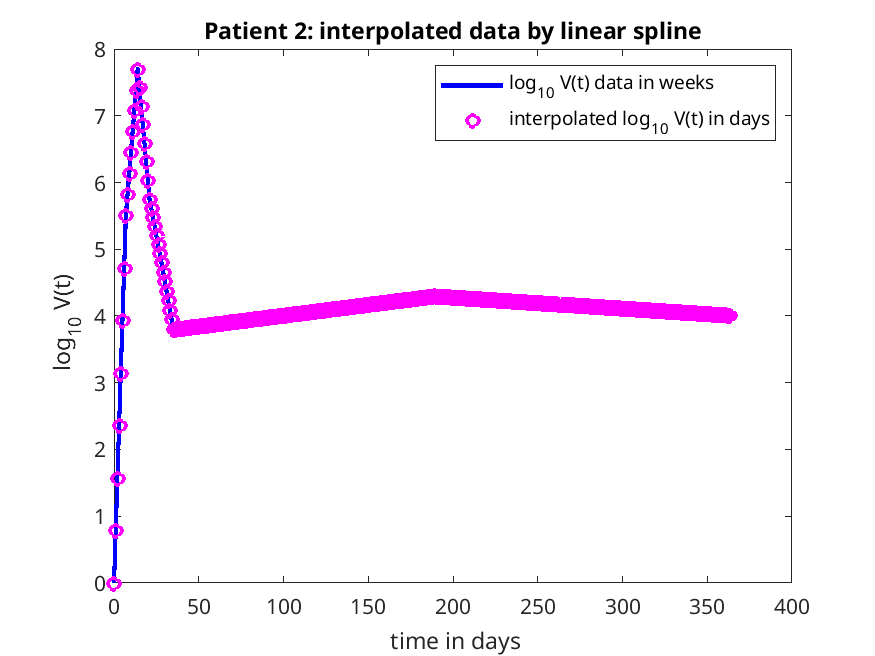}} &
  {\includegraphics[scale=0.5, clip = true, trim = 0.0cm 0.0cm 0.0cm 0.0cm]{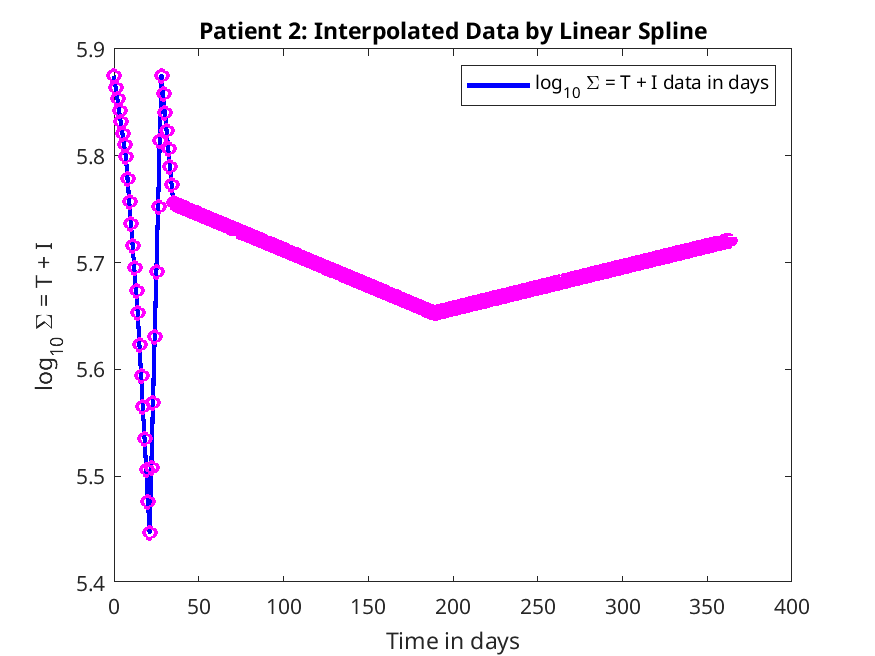}} \\
 $ V = g_2^0$ & $ \Sigma = g_1^0$  \\
\end{tabular}
\end{center}
\caption{Patient 2:  Interpolated clinical data by linear spline.}
 \label{fig:Pat2data}
 \end{figure}

For patient 2 we used interpolated clinical data  of the Table 2 presented in the Figure  \ref{fig:Pat2data}. 

For patient 2
 we define the parameters  in \eqref{CTLTest4}  as follows:
\begin{equation}\label{paramCTLtest2}
  \begin{split}
d_0 &= d = 0.26; \\
\beta_{CTL} &= 0.1;\\
\kappa  &=  1 + 10^5 \beta_{CTL}; \\
\delta T_1 &= 2.5; \\
\delta T_2 &= 5.0; \\
E_0^0 = 1,\\
t_1 &= 1,\\
t_2 &=50.
\end{split}
\end{equation}

\begin{figure}
\begin{center}
\begin{tabular}{ccc}
  {\includegraphics[scale=0.35, clip = true, trim = 0.0cm 0.0cm 0.0cm 0.0cm ]{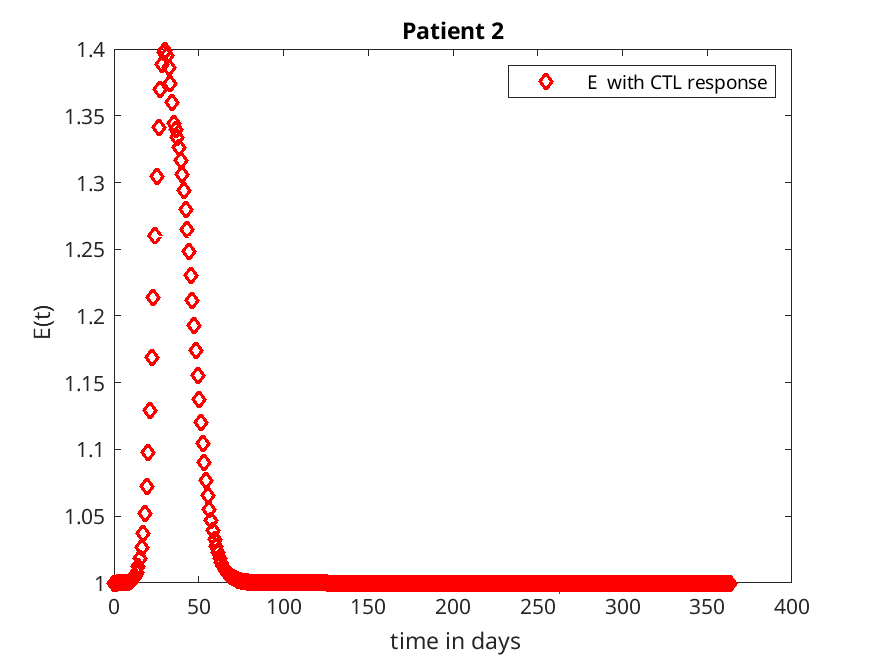}} &
  {\includegraphics[scale=0.35, clip = true, trim = 0.0cm 0.0cm 0.0cm 0.0cm]{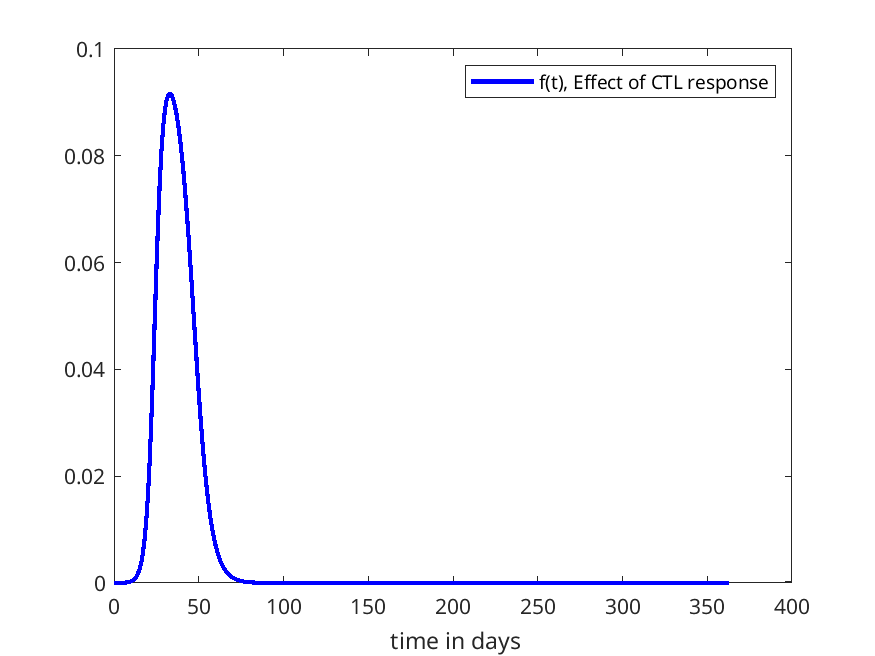}} &
    {\includegraphics[scale=0.35, clip = true, trim = 0.0cm 0.0cm 0.0cm 0.0cm ]{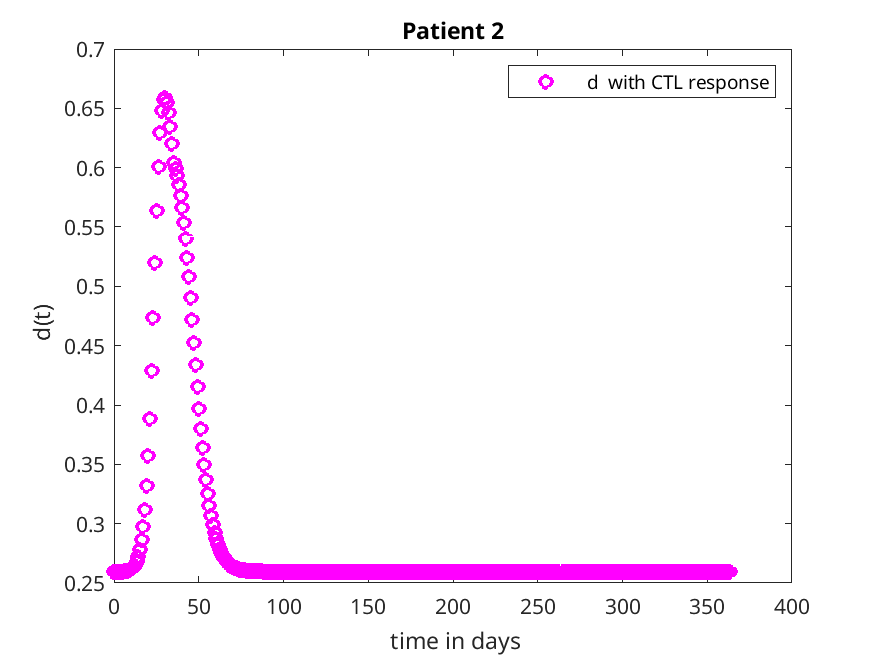}} \\
 $ E^0(t)$ & $f(t)$  & $d(t)$ \\
\end{tabular}
\end{center}
\caption{Patient 2: functions $E^0(t), f(t), d(t)$ with CTL response modelled as in \eqref{CTLTest4}  using parameters defined in \eqref{paramCTLPat1}.}
 \label{fig:Pat2CTLfunc}
 \end{figure}

\begin{figure}
\begin{center}
\begin{tabular}{cc}
  {\includegraphics[scale=0.5, clip = true, trim = 0.0cm 0.0cm 0.0cm 0.0cm ]{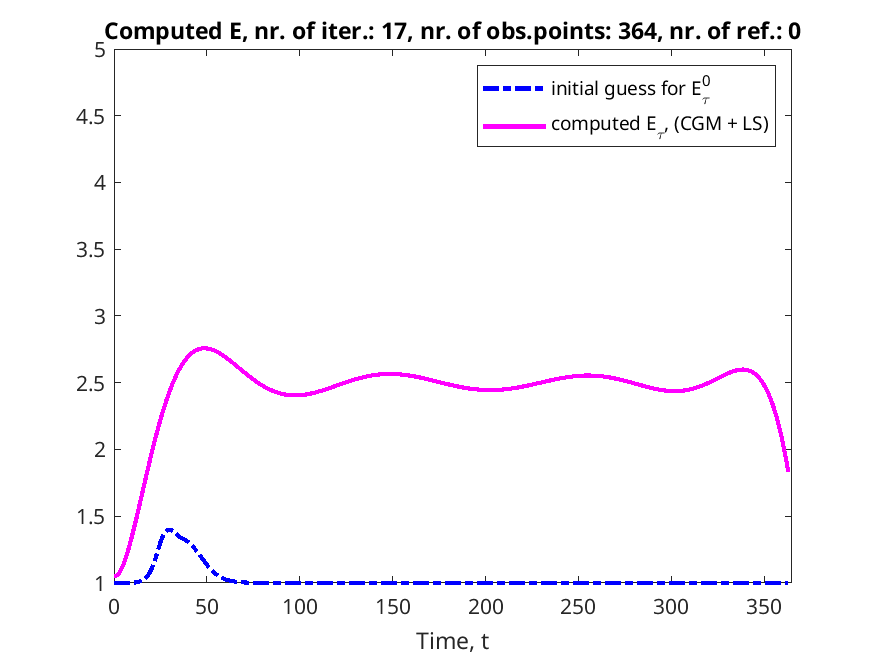}} &
  {\includegraphics[scale=0.5, clip = true, trim = 0.0cm 0.0cm 0.0cm 0.0cm]{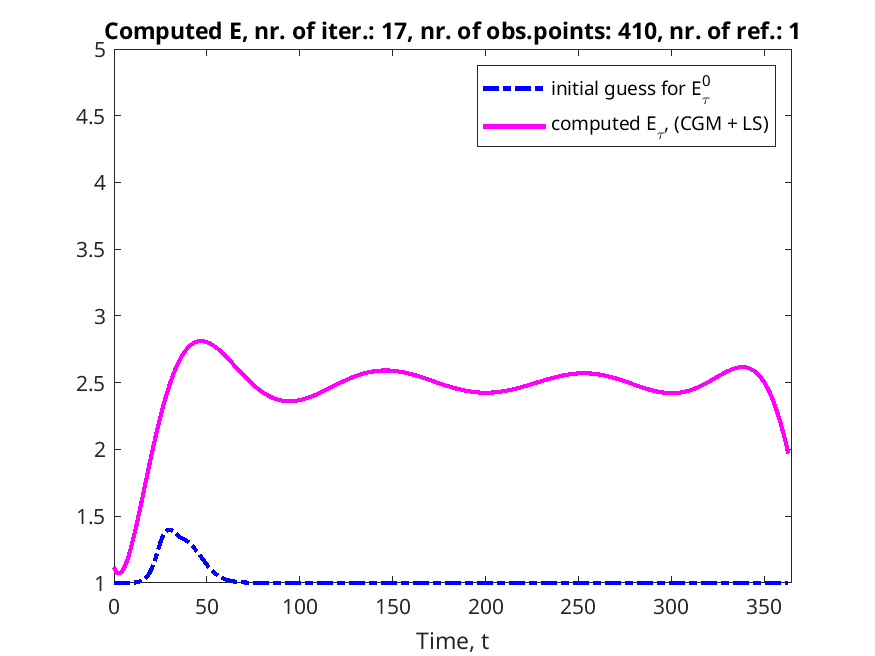}} \\
  $k=0$  & $k=1$  \\
\end{tabular}
\end{center}
\caption{Patient 2:  results of reconstruction of the function $E(t)$ on   $k$ times adaptivelly refined meshes $J_\tau^k$
 in ACGA algorithm. Computations are  performed for the mesh refinement parameter $\tilde{\beta_k}
  =0.875$  for all mesh refinements $k$.}
 \label{fig:Pat2RecE}
 \end{figure}

Figures \ref{fig:Pat2CTLfunc}  show  time-dependent  functions for initial guess $E^0(t)$ as well as $ f(t), d(t)$
 with CTL response  for parameters defined in \eqref{paramCTLPat1}.
These functions are used  now in optimization algorithms  for reconstruction of the immune response function $E(t)$.
 Figures \ref{fig:Pat2RecE} show
 computed functions $E_\tau^k$
 on   $k=0,1$ times adaptivelly refined meshes $J_\tau^k$.

\begin{figure}
\begin{center}
\begin{tabular}{ccc}
  {\includegraphics[scale=0.35, clip = true, trim = 0.0cm 0.0cm 0.0cm 0.0cm ]{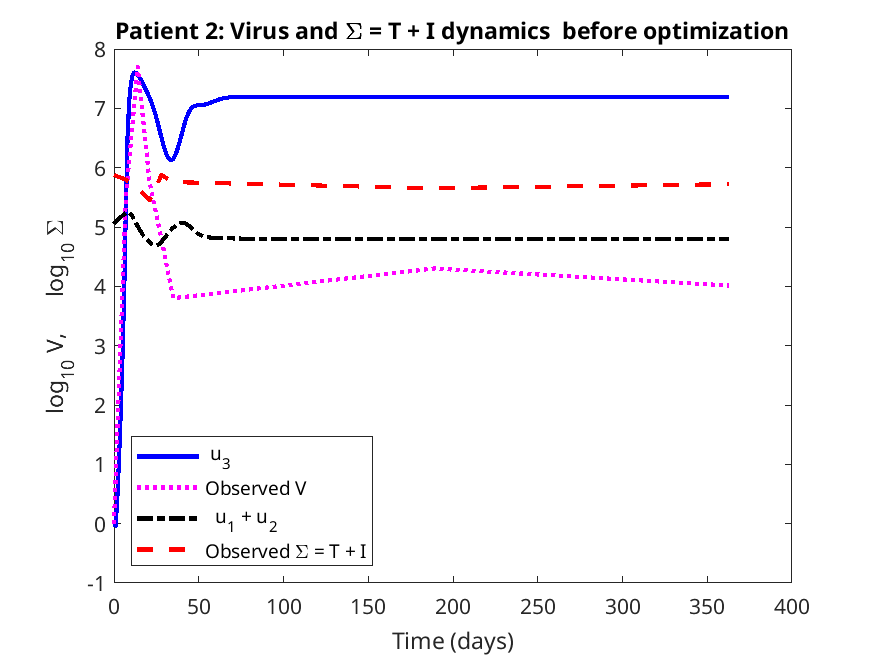}} &
  {\includegraphics[scale=0.35, clip = true, trim = 0.0cm 0.0cm 0.0cm 0.0cm]{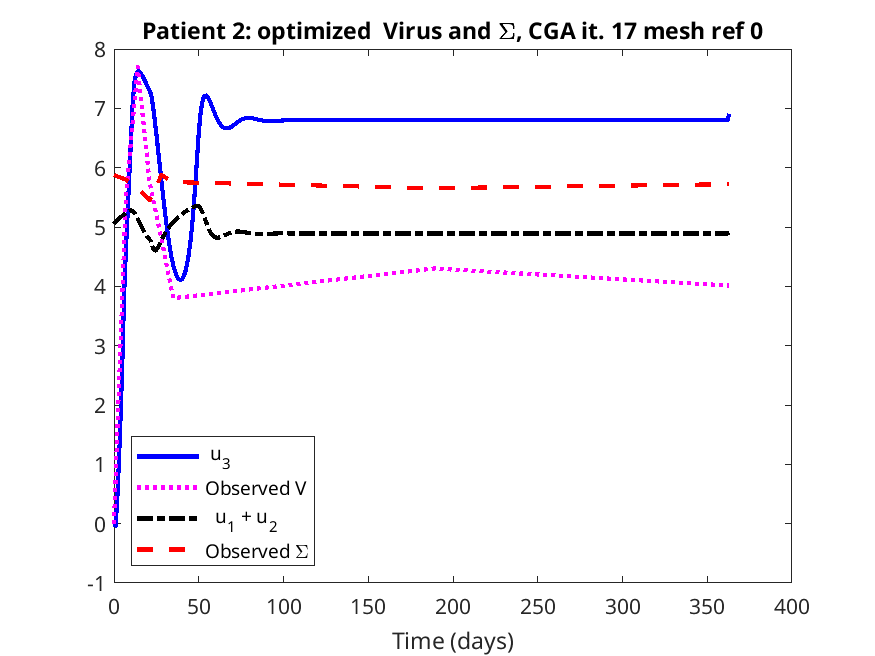}} &
   {\includegraphics[scale=0.35, clip = true, trim = 0.0cm 0.0cm 0.0cm 0.0cm]{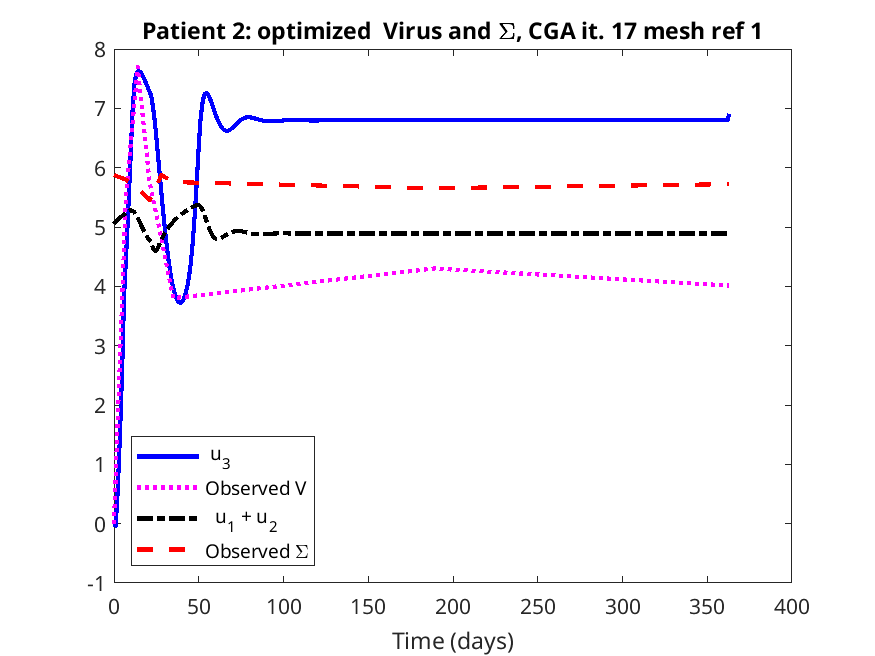}} \\
 a) $k=0$ & b) $k=0$  & c) $k=1$ \\
\end{tabular}
\end{center}
\caption{ Patient 2: Dynamics of the computed virus function $V_\tau^k$ before and after optimization corresponding to the
   computed $E_\tau^k$ on   $k, k=0,1$ times adaptivelly refined meshes versus interpolated clinical data ${g_1^0}_\tau,{g_2^0}_\tau$.}
 \label{fig:Pat2virus}
 \end{figure}

\begin{figure}
\begin{center}
\begin{tabular}{cc}
  {\includegraphics[scale=0.5, clip = true, trim = 0.0cm 0.0cm 0.0cm 0.0cm ]{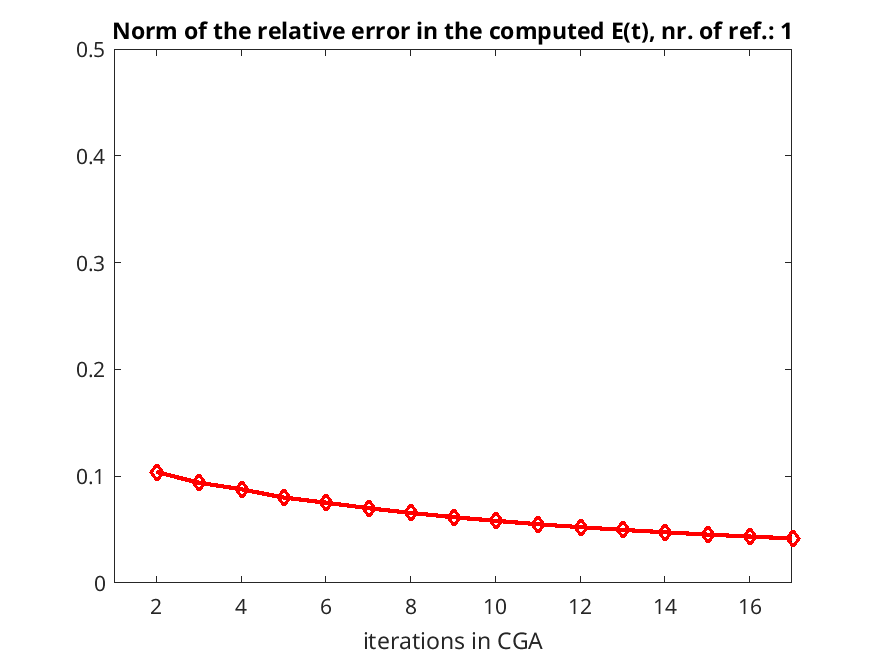}} &
  {\includegraphics[scale=0.5, clip = true, trim = 0.0cm 0.0cm 0.0cm 0.0cm]{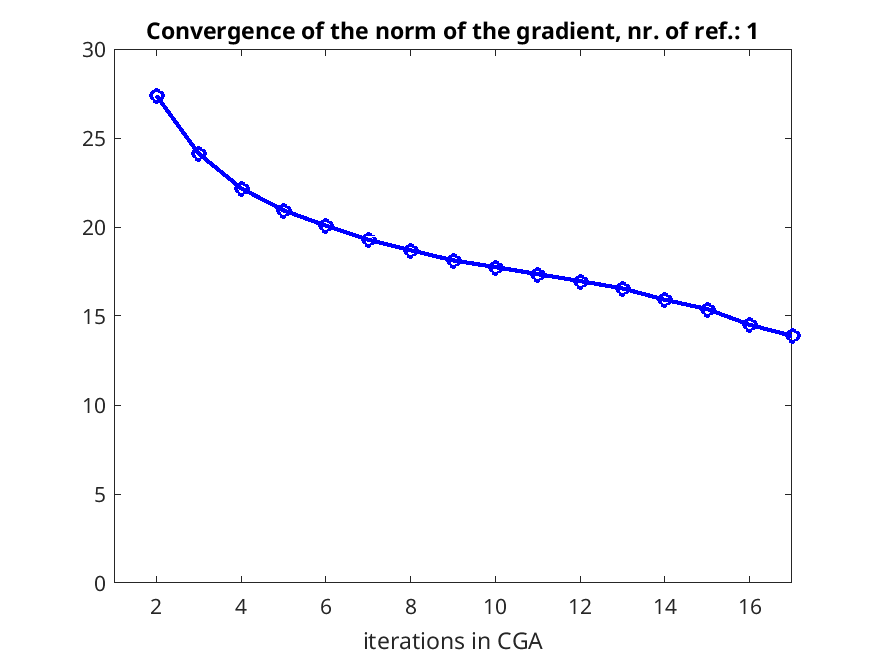}} \\
  $\frac{  \| E_\tau^{m} - E_\tau^{m-1} \|_{L_2(\Omega_t)}}{\| E_\tau^m\|_{L_2(\Omega_t)}}$   &  $\|G^m(t) \|_{L_2(\Omega_t)} $ 
\end{tabular}
\end{center}
\caption{Patient 2:   computed relative norms  $\frac{  \| E_\tau^{m} - E_\tau^{m-1} \|_{L_2(\Omega_t)}}{\| E_\tau^m\|_{L_2(\Omega_t)}}$
and $\|G^m(t) \|_{L_2(\Omega_t)}$  on the mesh $J_\tau^1$.}
 \label{fig:Pat2Residuals}
 \end{figure}

\begin{figure}
\begin{center}
\begin{tabular}{cc}
  {\includegraphics[scale=0.5, clip = true, trim = 0.0cm 0.0cm 0.0cm 0.0cm ]{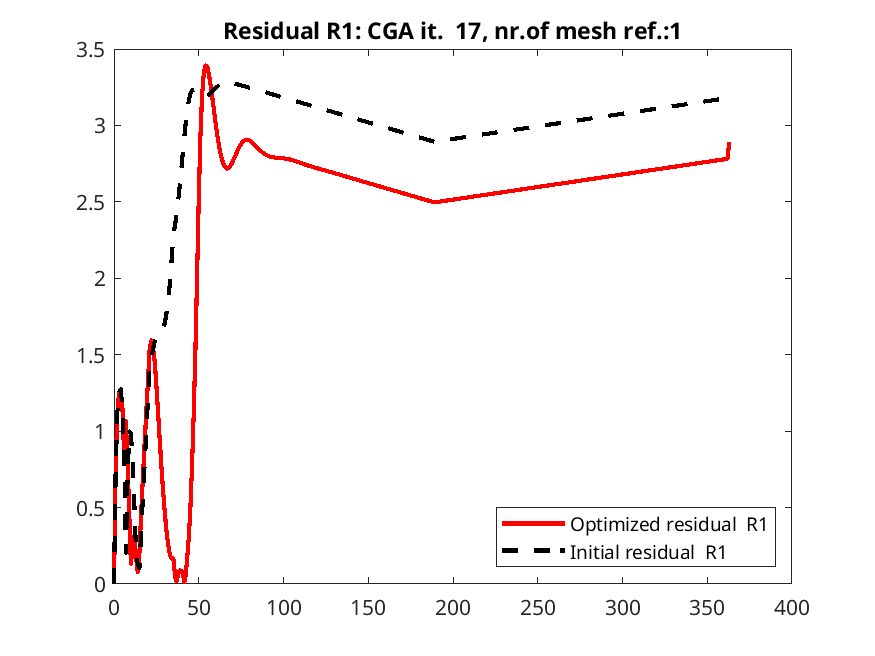}} &
  {\includegraphics[scale=0.5, clip = true, trim = 0.0cm 0.0cm 0.0cm 0.0cm]{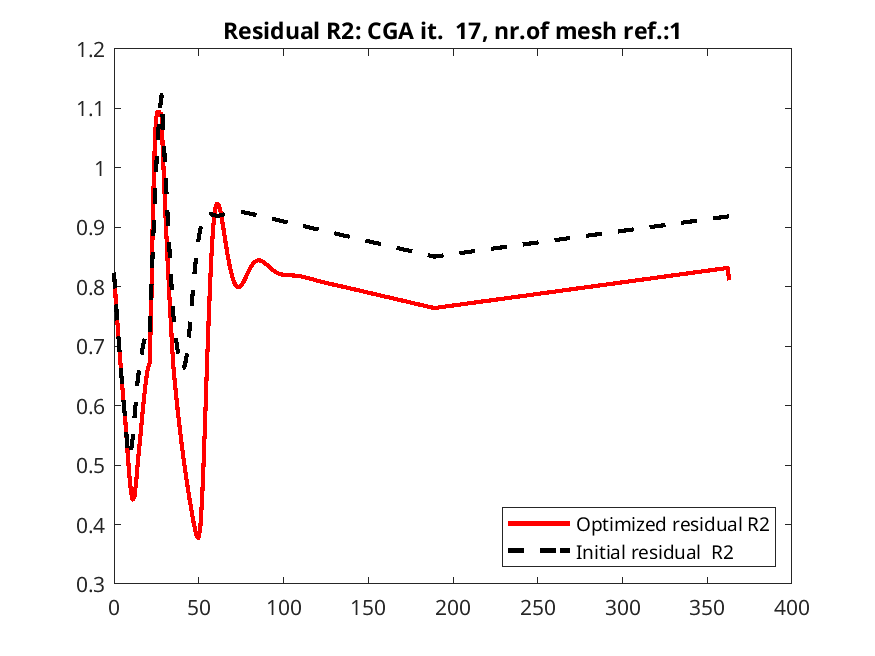}} \\
   a) $  R_1 $  & b) $ R_2  $  \\
 {\includegraphics[scale=0.5, clip = true, trim = 0.0cm 0.0cm 0.0cm 0.0cm]{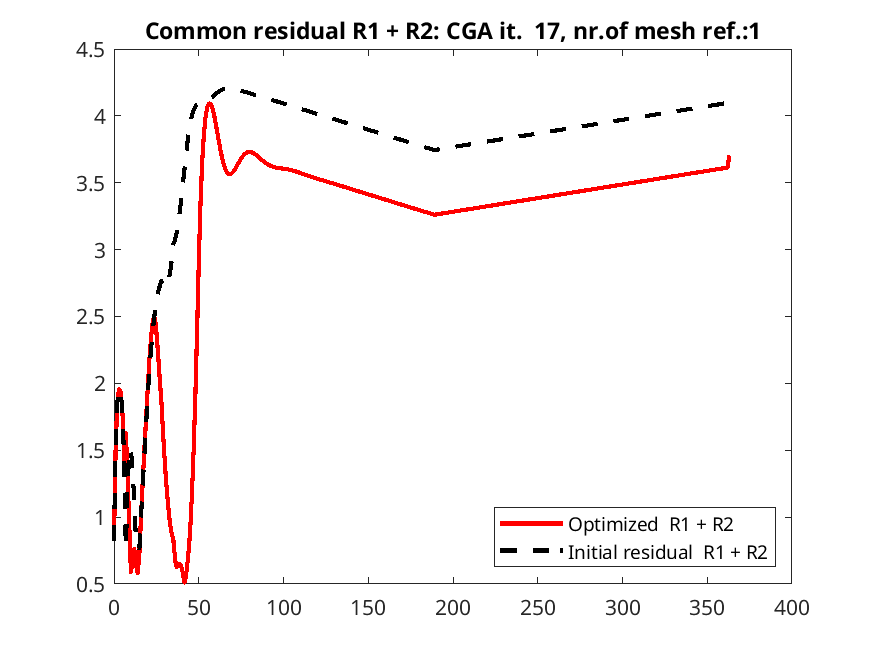}} &
  {\includegraphics[scale=0.5, clip = true, trim = 0.0cm 0.0cm 0.0cm 0.0cm]{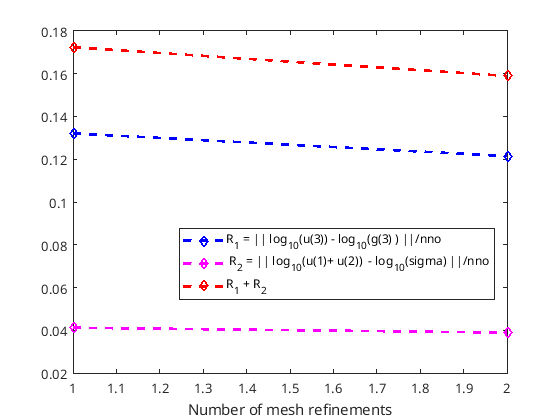}}  \\
 c)  $ R_1 + R_2  $ &  d) 
\end{tabular}
\end{center}
\caption{Patient 2: a), b), c)  computed residuals $R_1$ and $R_2$ on the mesh $J_\tau^1$; d) Comparison of related residuals
 $\|R_1\|, \|R_2\|$ on the meshes
$J_\tau^k, k = 0,1$. }
 \label{fig:Pat2R1R2}
 \end{figure}

Figures \ref{fig:Pat2virus}  show  time-dependent behaviour  of the computed virus function $V_\tau^k, k=0,1$ before and after applying ACGA.
Here,  the computed virus function $V_\tau^k$ corresponds to the
   computed $E_\tau^k$ presented in Figure \ref{fig:Pat2RecE}.
Results are compared with  clinical data  for the observed virus function  $V = g_2^0$  and observed total number of  the uninfected
and infected cells $\Sigma = g_1^0$.
 Computed optimized virus function on $J_\tau^0$ is presented in the
Figure \ref{fig:Pat2virus}-b).  Figure \ref{fig:Pat2virus}-c) shows
computed virus function on the one time locally refined time mesh
$J_\tau^1$.  Comparing figures  \ref{fig:Pat2virus}-b) and  \ref{fig:Pat2virus}-c)  we observe that local adaptive mesh refinement allows more exact fit data to achieve minimum of the virus load.
Figures   \ref{fig:Pat2Residuals} show
stabilization of the computed relative norms  
$\frac{  \| E_\tau^{m} - E_\tau^{m-1} \|_{L_2(\Omega_t)}}{\| E_\tau^m\|_{L_2(\Omega_t)}}$ 
  and behaviour of norms of the computed gradient $\|G^m(t) \|_{L_2(\Omega_t)} $  for patient 2.
 Relative residuals $R_1$ and $R_2$ and $\|R_1\|_{L_2(\Omega_t)} , \|R_2\|_{L_2(\Omega_t)} $ are presented in the Figures
 \ref{fig:Pat2R1R2} as well as in the Table 4.


\subsection{Results for Patient 3}

\begin{figure}
\begin{center}
\begin{tabular}{ccc}
  {\includegraphics[scale=0.5, clip = true, trim = 0.0cm 0.0cm 0.0cm 0.0cm ]{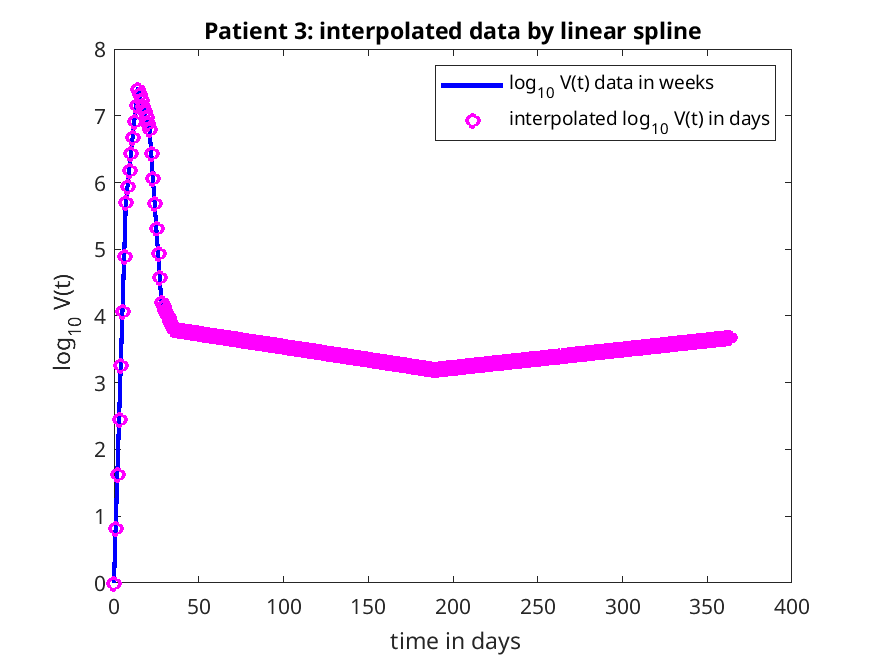}} &
  {\includegraphics[scale=0.5, clip = true, trim = 0.0cm 0.0cm 0.0cm 0.0cm]{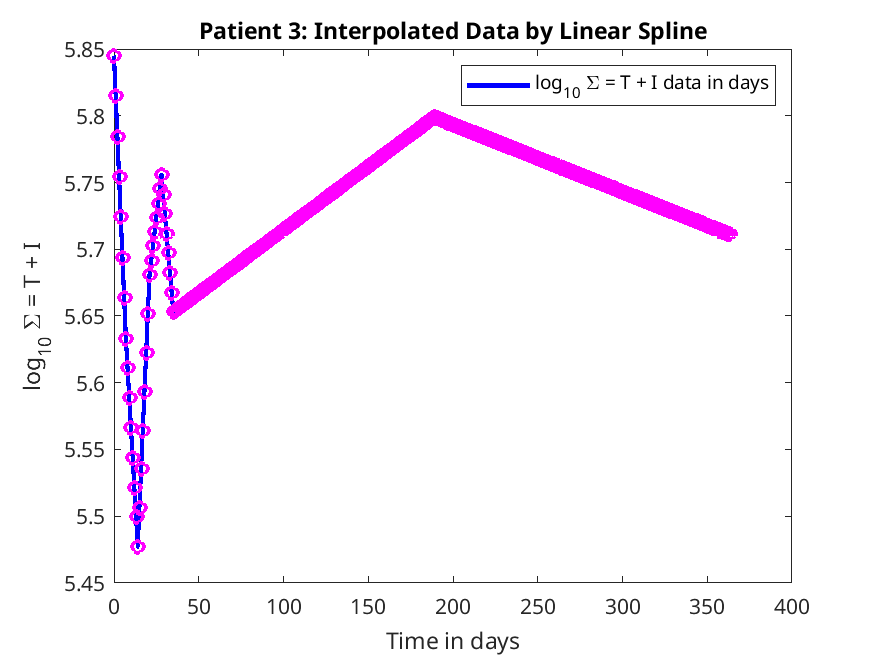}} \\
 $ V = g_2^0$ & $ \Sigma = g_1^0$  \\
\end{tabular}
\end{center}
\caption{Patient 3:  Interpolated clinical data by linear spline.}
 \label{fig:Pat3data}
 \end{figure}

For patient 3 we used interpolated clinical data  of the Table 2 presented in the Figure  \ref{fig:Pat3data}.
We set all parameters in \eqref{CTLTest4} the same as for patient 2.

\begin{figure}
\begin{center}
\begin{tabular}{ccc}
  {\includegraphics[scale=0.35, clip = true, trim = 0.0cm 0.0cm 0.0cm 0.0cm ]{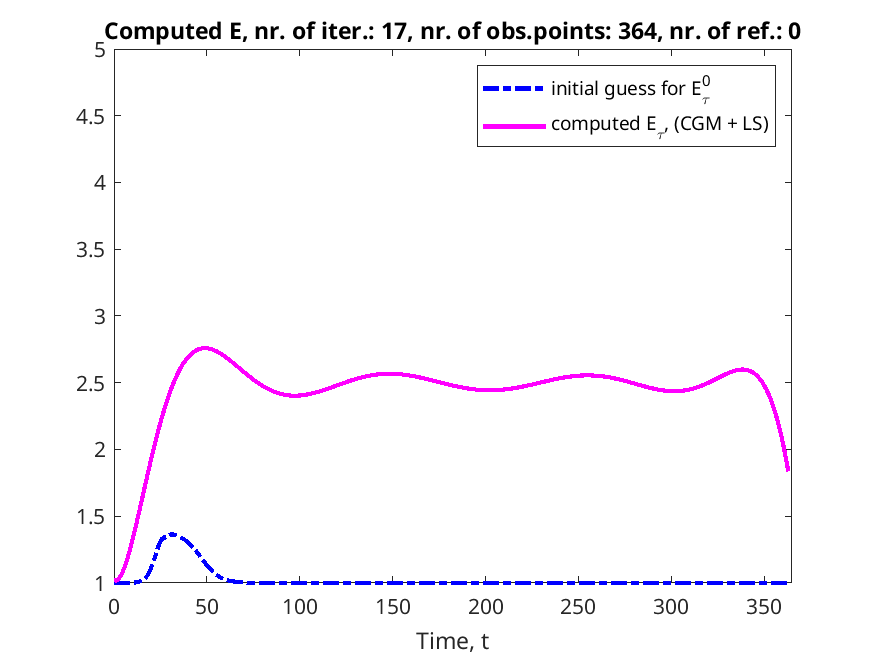}} &
  {\includegraphics[scale=0.35, clip = true, trim = 0.0cm 0.0cm 0.0cm 0.0cm]{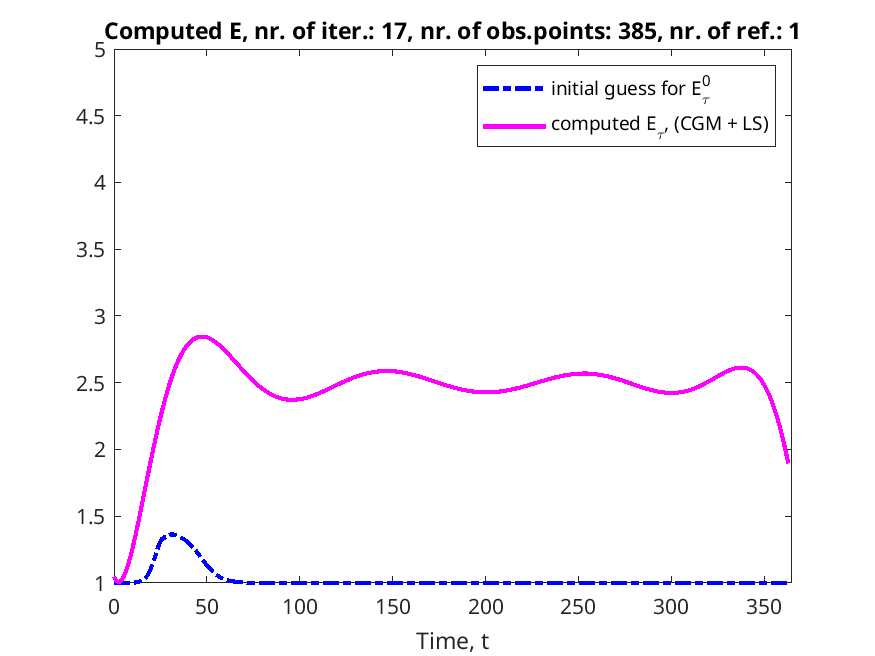}} &
    {\includegraphics[scale=0.35, clip = true, trim = 0.0cm 0.0cm 0.0cm 0.0cm]{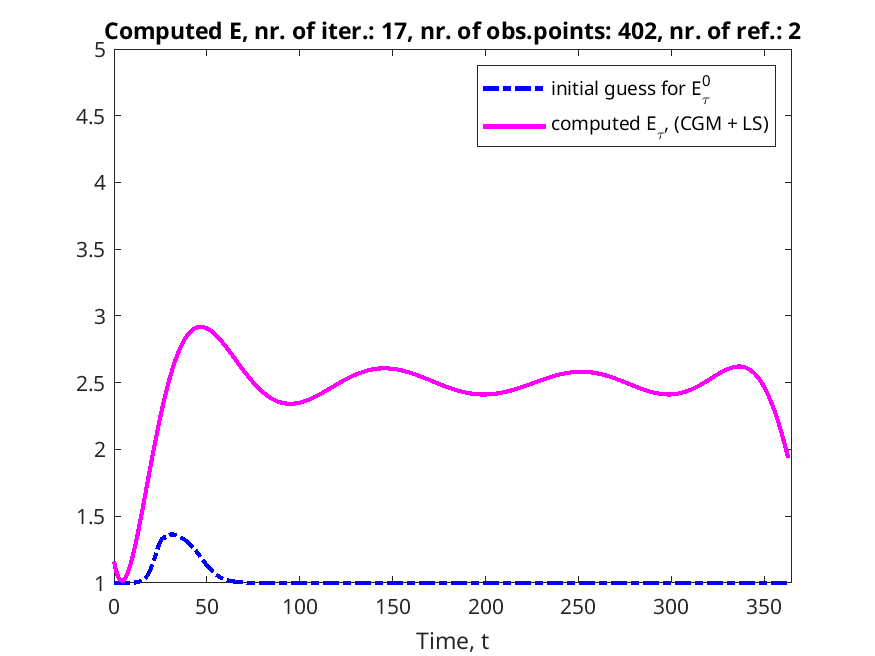}} \\
a)  $k=0$  & b) $k=1$ & c) $k=2$  \\
\end{tabular}
\end{center}
\caption{Patient 3:  results of reconstruction of the function $E(t)$ on   $k$ times adaptivelly refined meshes $J_\tau^k$
 in ACGA algorithm. Computations are  performed for the mesh refinement parameter $\tilde{\beta_k}
  =0.875$  for all mesh refinements $k$.}
 \label{fig:Pat3RecE}
 \end{figure}

Figures \ref{fig:Pat2CTLfunc}  show  time-dependent  functions for initial guess $E^0(t)$ as well as $ f(t), d(t)$
 with CTL response  for parameters defined in \eqref{paramCTLPat1} which are the same as for the patient 2.
These functions are used  now in optimization algorithms  for reconstruction of the immune response function $E(t)$.
 Figures \ref{fig:Pat3RecE} show
 computed functions $E_\tau^k$
 on   $k=0,1,2$ times adaptivelly refined meshes $J_\tau^k$.

\begin{figure}
\begin{center}
\begin{tabular}{ccc}
  {\includegraphics[scale=0.35, clip = true, trim = 0.0cm 0.0cm 0.0cm 0.0cm ]{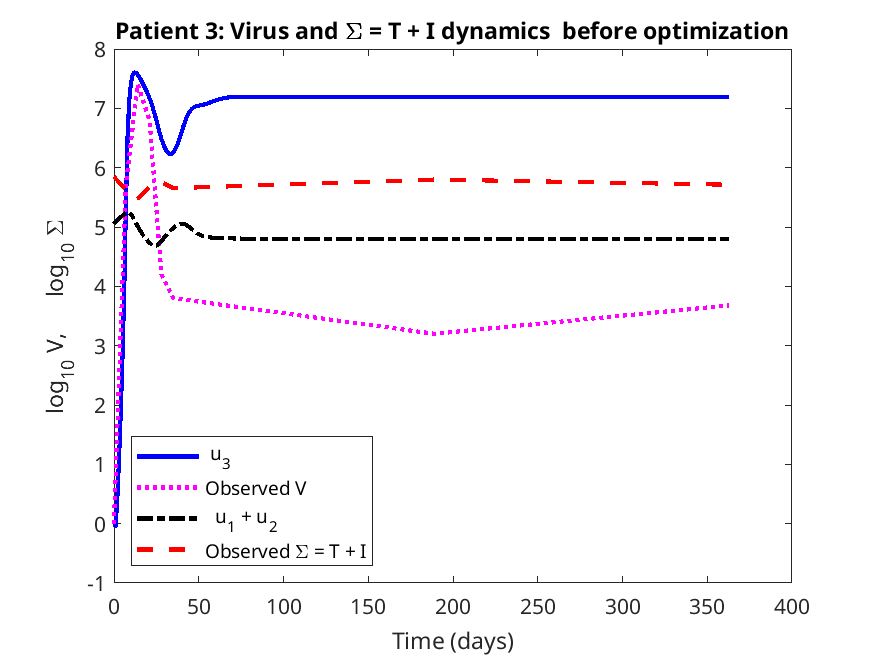}} &
  {\includegraphics[scale=0.35, clip = true, trim = 0.0cm 0.0cm 0.0cm 0.0cm]{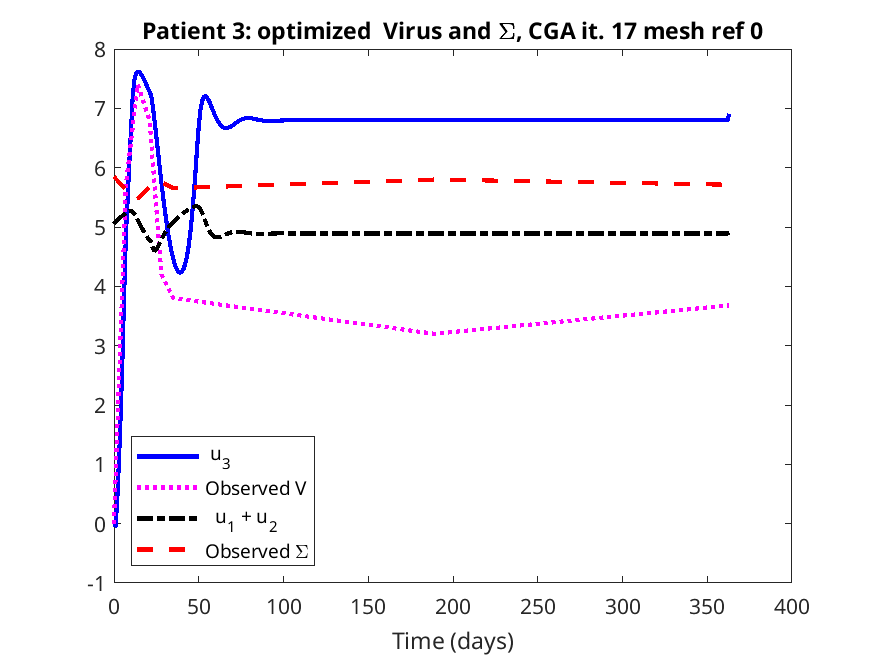}} &
   {\includegraphics[scale=0.35, clip = true, trim = 0.0cm 0.0cm 0.0cm 0.0cm]{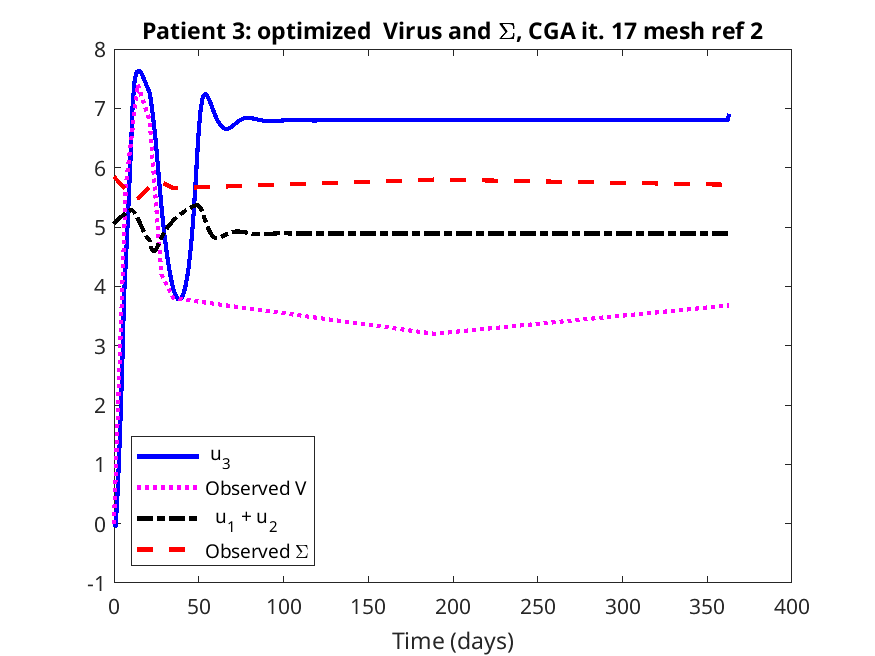}} \\
 a) $k=0$ & b) $k=0$  & c) $k=2$ \\
\end{tabular}
\end{center}
\caption{ Patient 3: Dynamics of the computed virus function $V_\tau^k$ before and after optimization corresponding to the
   computed $E_\tau^k$ on   $k, k=0,2$ times adaptivelly refined meshes versus interpolated clinical data ${g_1^0}_\tau,{g_2^0}_\tau$.}
 \label{fig:Pat3virus}
 \end{figure}

\begin{figure}
\begin{center}
\begin{tabular}{cc}
  {\includegraphics[scale=0.5, clip = true, trim = 0.0cm 0.0cm 0.0cm 0.0cm ]{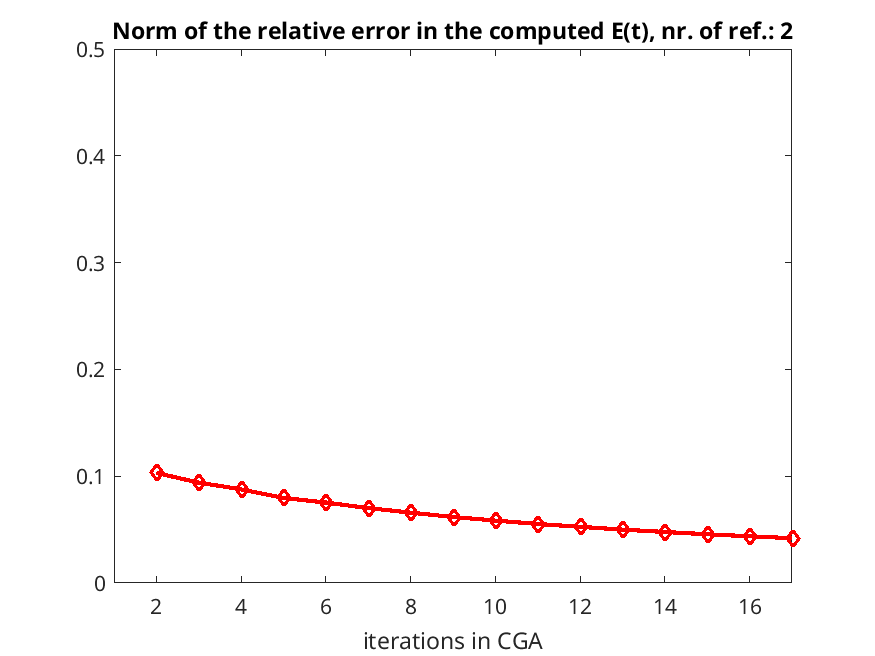}} &
  {\includegraphics[scale=0.5, clip = true, trim = 0.0cm 0.0cm 0.0cm 0.0cm]{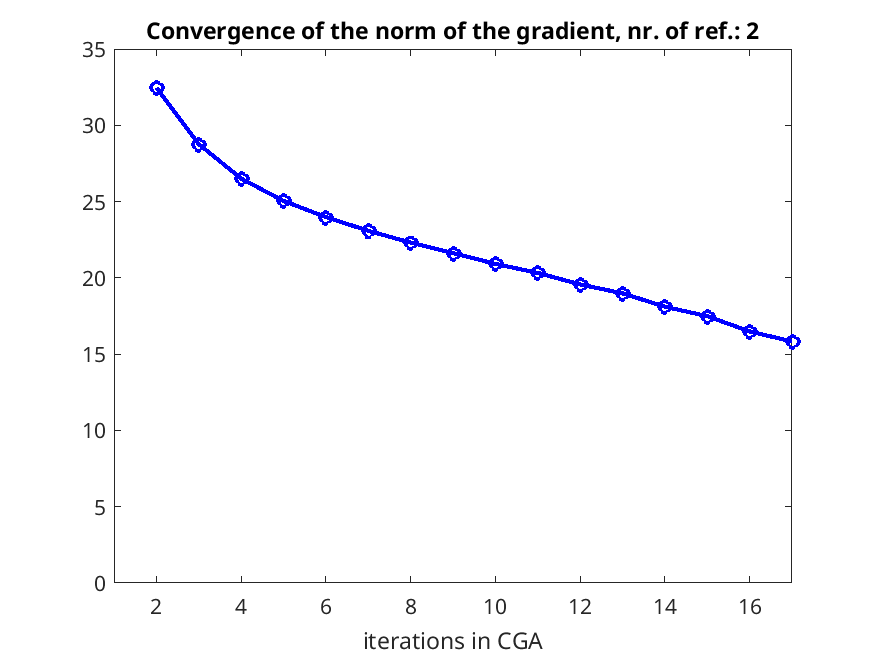}} \\
  $\frac{  \| E_\tau^{m} - E_\tau^{m-1} \|_{L_2(\Omega_t)}}{\| E_\tau^m\|_{L_2(\Omega_t)}}$   &  $\|G^m(t) \|_{L_2(\Omega_t)} $ 
\end{tabular}
\end{center}
\caption{Patient 3:   computed relative norms  $\frac{  \| E_\tau^{m} - E_\tau^{m-1} \|_{L_2(\Omega_t)}}{\| E_\tau^m\|_{L_2(\Omega_t)}}$
and $\|G^m(t) \|_{L_2(\Omega_t)}$  on the mesh $J_\tau^2$.}
 \label{fig:Pat3Residuals}
 \end{figure}

\begin{figure}
\begin{center}
\begin{tabular}{cc}
  {\includegraphics[scale=0.5, clip = true, trim = 0.0cm 0.0cm 0.0cm 0.0cm ]{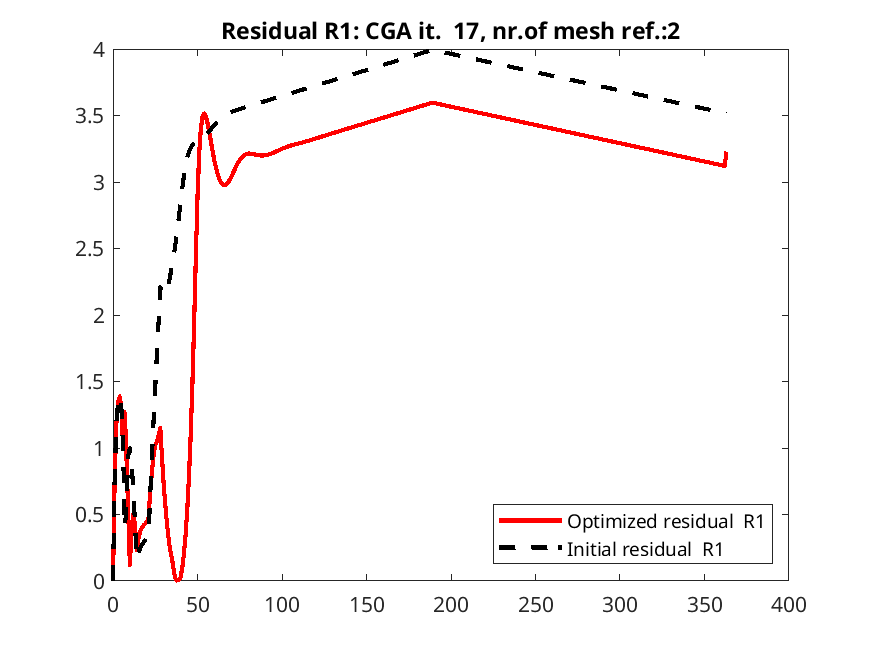}} &
  {\includegraphics[scale=0.5, clip = true, trim = 0.0cm 0.0cm 0.0cm 0.0cm]{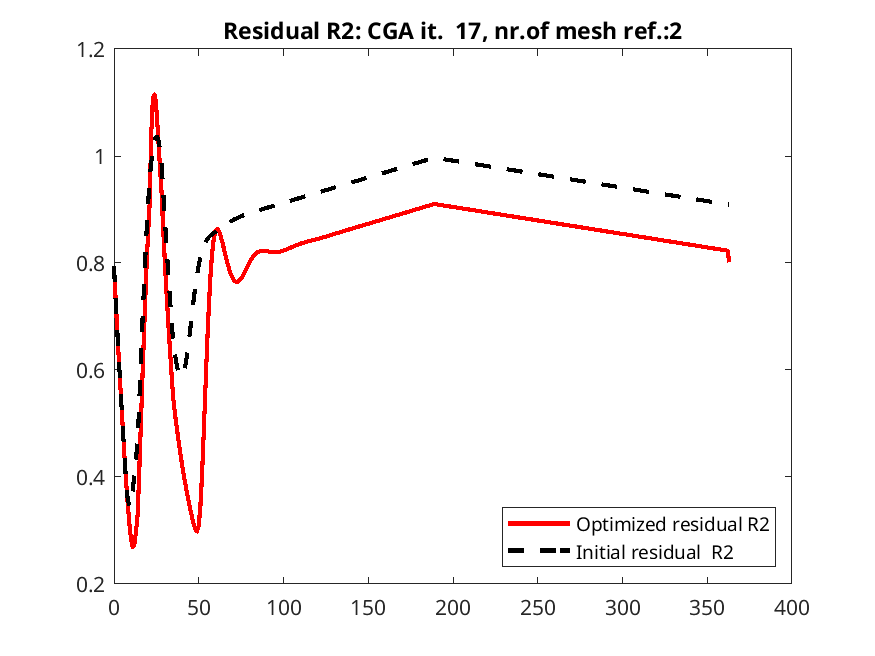}} \\
   a) $  R_1 $  & b) $ R_2  $  \\
 {\includegraphics[scale=0.5, clip = true, trim = 0.0cm 0.0cm 0.0cm 0.0cm]{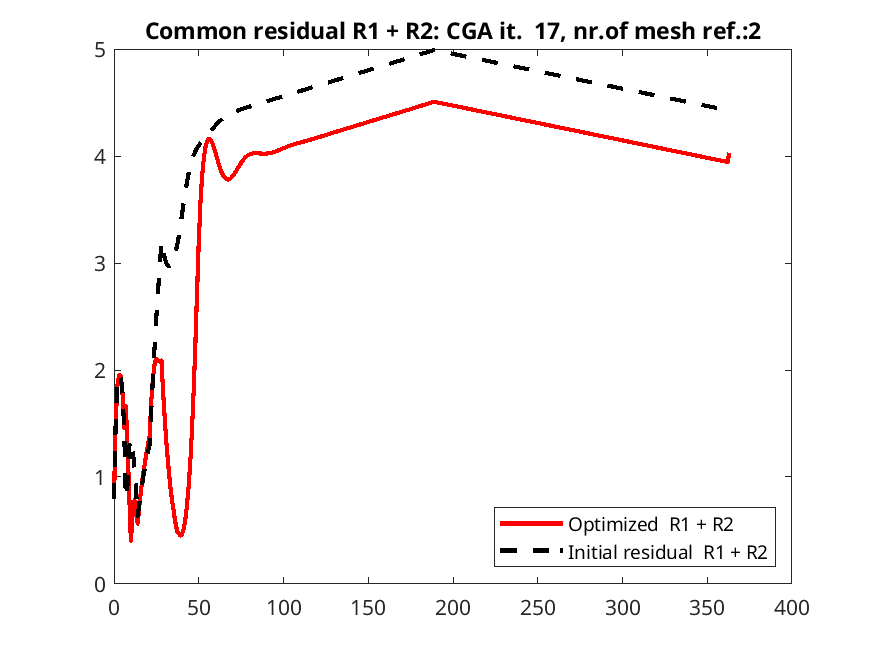}} &
  {\includegraphics[scale=0.5, clip = true, trim = 0.0cm 0.0cm 0.0cm 0.0cm]{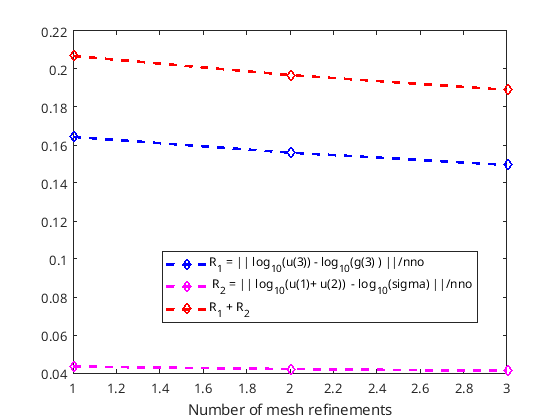}}  \\
 c)  $ R_1 + R_2  $ &  d) 
\end{tabular}
\end{center}
\caption{Patient 3: a), b), c)  computed residuals $R_1$ and $R_2$ on the mesh $J_\tau^2$; d) Comparison of related residuals
 $\|R_1\|, \|R_2\|$ on the meshes
$J_\tau^k, k = 0,1,2$. }
 \label{fig:Pat3R1R2}
 \end{figure}

Figures \ref{fig:Pat3virus}  show  time-dependent behaviour  of the computed virus function $V_\tau^k, k=0,2$ before and after applying ACGA.
Here,  the computed virus function $V_\tau^k$ corresponds to the
   computed $E_\tau^k$ presented in Figure \ref{fig:Pat3RecE}.
Results are compared with  clinical data  for the observed virus function  $V = g_2^0$  and observed total number of  the uninfected
and infected cells $\Sigma = g_1^0$.
 Computed optimized virus function on $J_\tau^0$ is presented in the
Figure \ref{fig:Pat3virus}-b).  Figure \ref{fig:Pat3virus}-c) shows
computed virus function on the twice  locally refined time mesh
$J_\tau^2$.  Comparing figures  \ref{fig:Pat3virus}-b) and  \ref{fig:Pat3virus}-c)
we  again observe that local adaptive mesh refinement allows more exact fit data to achieve minimum of the virus load.
Figures   \ref{fig:Pat3Residuals} show
stabilization of the computed relative norms  
$\frac{  \| E_\tau^{m} - E_\tau^{m-1} \|_{L_2(\Omega_t)}}{\| E_\tau^m\|_{L_2(\Omega_t)}}$ 
  and behaviour of norms of the computed gradient $\|G^m(t) \|_{L_2(\Omega_t)} $  for patient 3.
 Relative residuals $R_1$ and $R_2$ and $\|R_1\|_{L_2(\Omega_t)} , \|R_2\|_{L_2(\Omega_t)} $ are presented in the Figures
 \ref{fig:Pat3R1R2} as well as in the Table 4.


\subsection{Results for Patient 4}

\begin{figure}
\begin{center}
\begin{tabular}{ccc}
  {\includegraphics[scale=0.5, clip = true, trim = 0.0cm 0.0cm 0.0cm 0.0cm ]{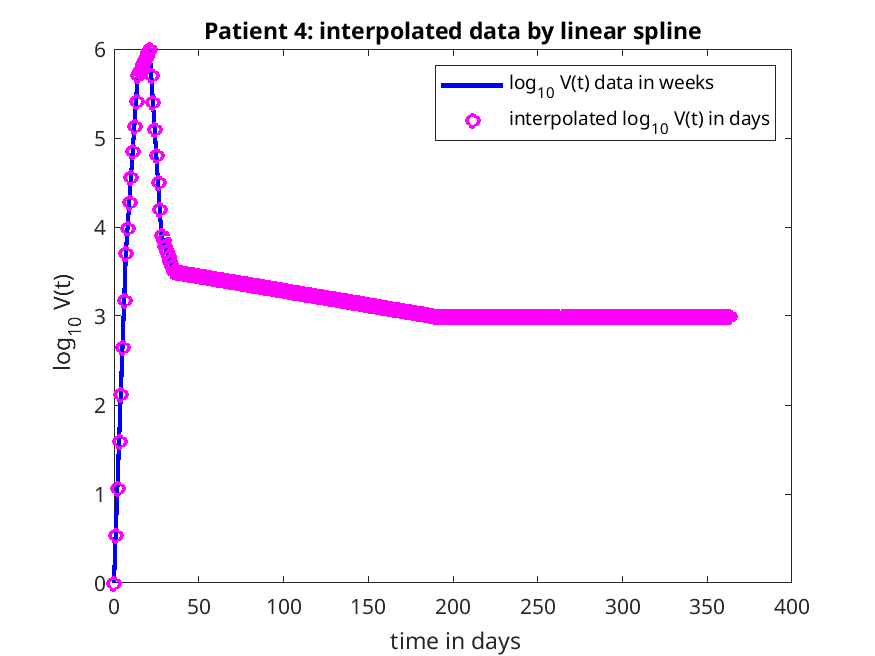}} &
  {\includegraphics[scale=0.5, clip = true, trim = 0.0cm 0.0cm 0.0cm 0.0cm]{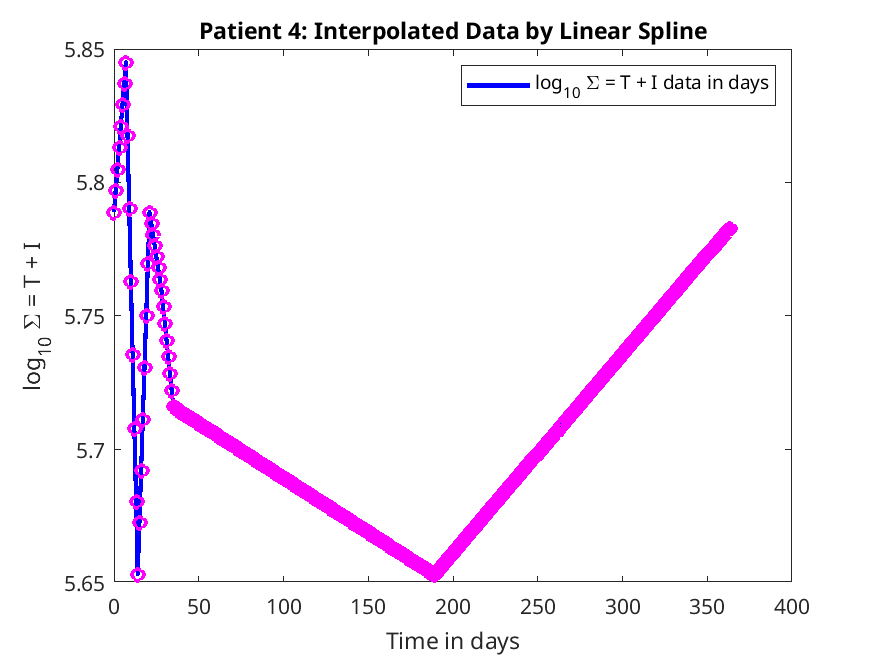}} \\
 $ V = g_2^0$ & $ \Sigma = g_1^0$  \\
\end{tabular}
\end{center}
\caption{Patient 4:  Interpolated clinical data by linear spline.}
 \label{fig:Pat4data}
 \end{figure}

For patient 4 we used interpolated clinical data of the Table 2
presented in the Figure \ref{fig:Pat4data}.  We observe that the data
is very discontinuos. However, we test the same ACGA algorithm for
this dataset to see performance of the code.

For patient 4
 we define the parameters  in \eqref{CTLTest4}  as follows:
\begin{equation}\label{paramCTLtest2}
  \begin{split}
d_0 &= d = 0.26; \\
\beta_{CTL} &= 0.1;\\
\kappa  &=  1 + 10^5 \beta_{CTL}; \\
\delta T_1 &= 2.5; \\
\delta T_2 &= 5.0; \\
E_0 = 1.4,\\
t_1 &= 1,\\
t_2 &= 1.
\end{split}
\end{equation}

\begin{figure}
\begin{center}
\begin{tabular}{ccc}
  {\includegraphics[scale=0.35, clip = true, trim = 0.0cm 0.0cm 0.0cm 0.0cm ]{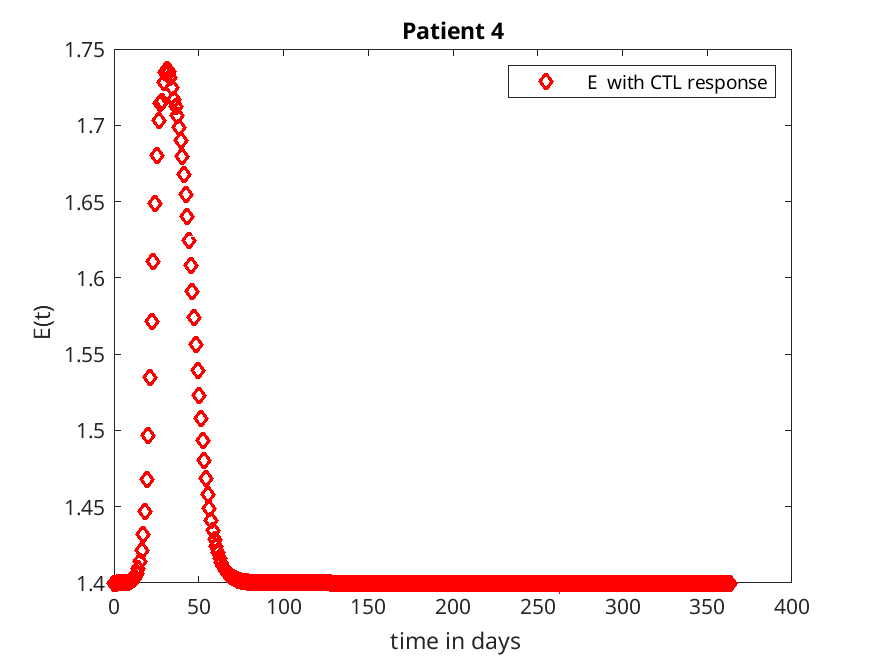}} &
  {\includegraphics[scale=0.35, clip = true, trim = 0.0cm 0.0cm 0.0cm 0.0cm]{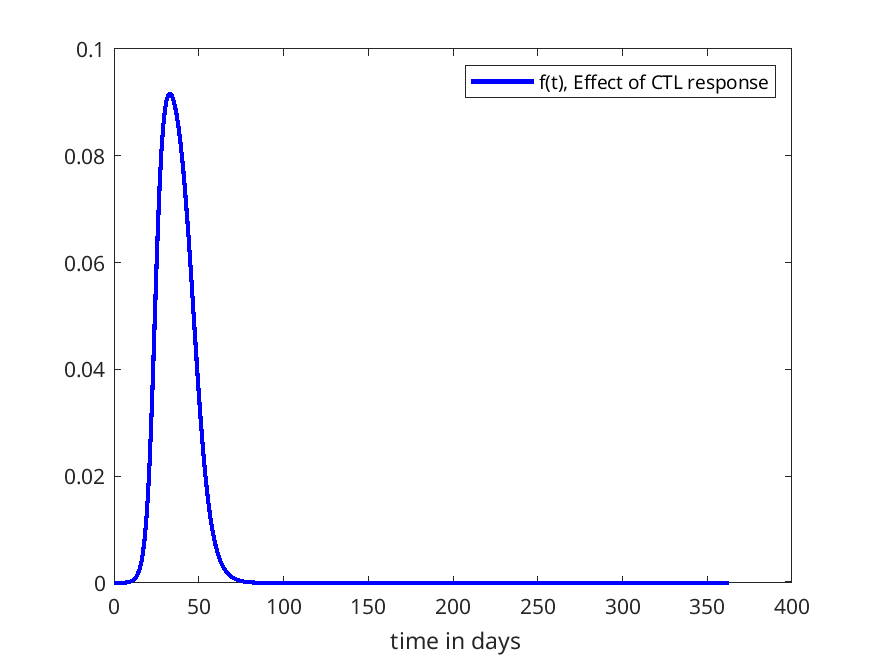}} &
    {\includegraphics[scale=0.35, clip = true, trim = 0.0cm 0.0cm 0.0cm 0.0cm ]{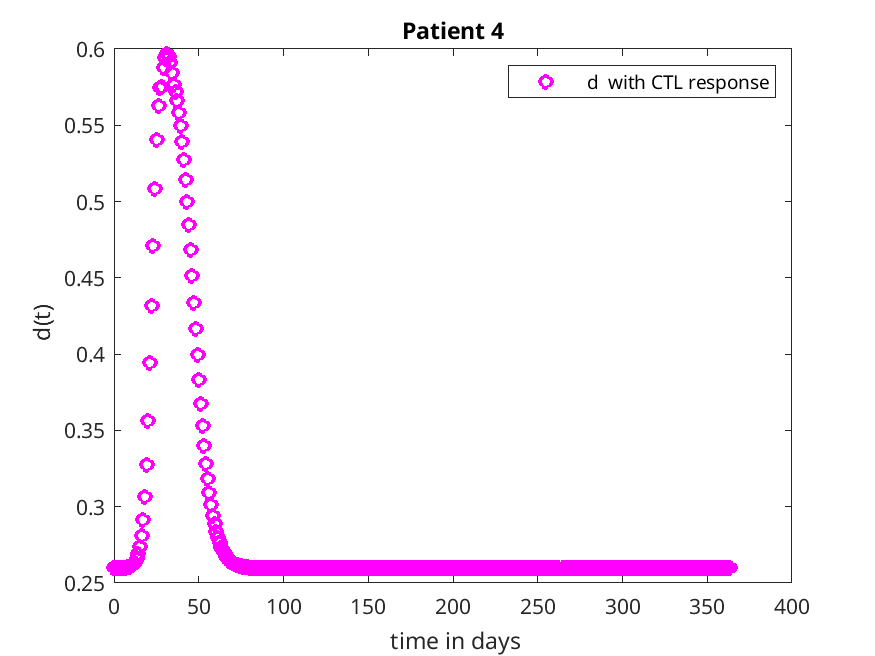}} \\
 $ E^0(t)$ & $f(t)$  & $d(t)$ \\
\end{tabular}
\end{center}
\caption{Patient 4: functions $E^0(t), f(t), d(t)$ with CTL response modelled as in \eqref{CTLTest4}  using parameters defined in \eqref{paramCTLPat1}.}
 \label{fig:Pat4CTLfunc}
 \end{figure}

\begin{figure}
\begin{center}
\begin{tabular}{ccc}
  {\includegraphics[scale=0.35, clip = true, trim = 0.0cm 0.0cm 0.0cm 0.0cm ]{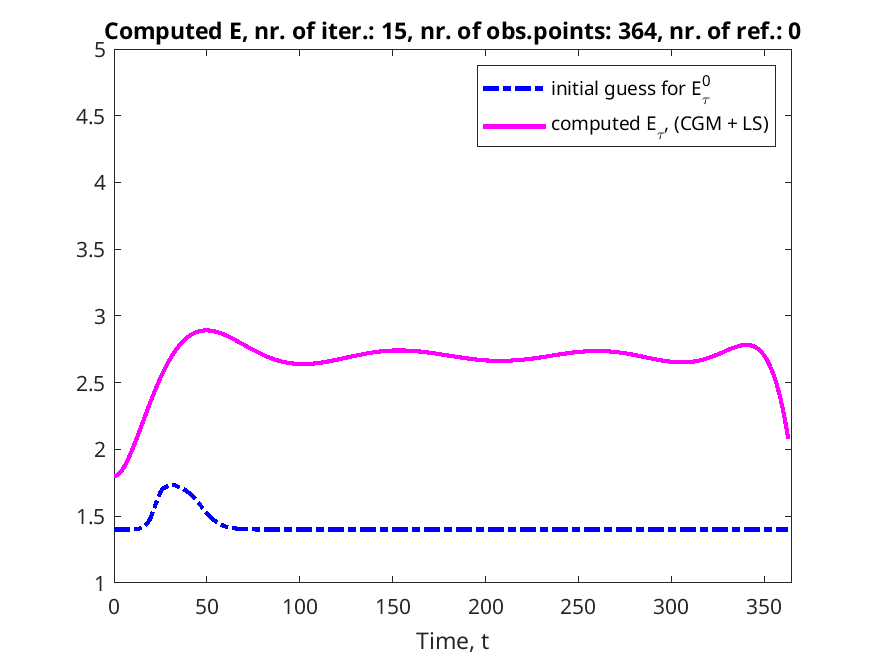}} &
  {\includegraphics[scale=0.35, clip = true, trim = 0.0cm 0.0cm 0.0cm 0.0cm]{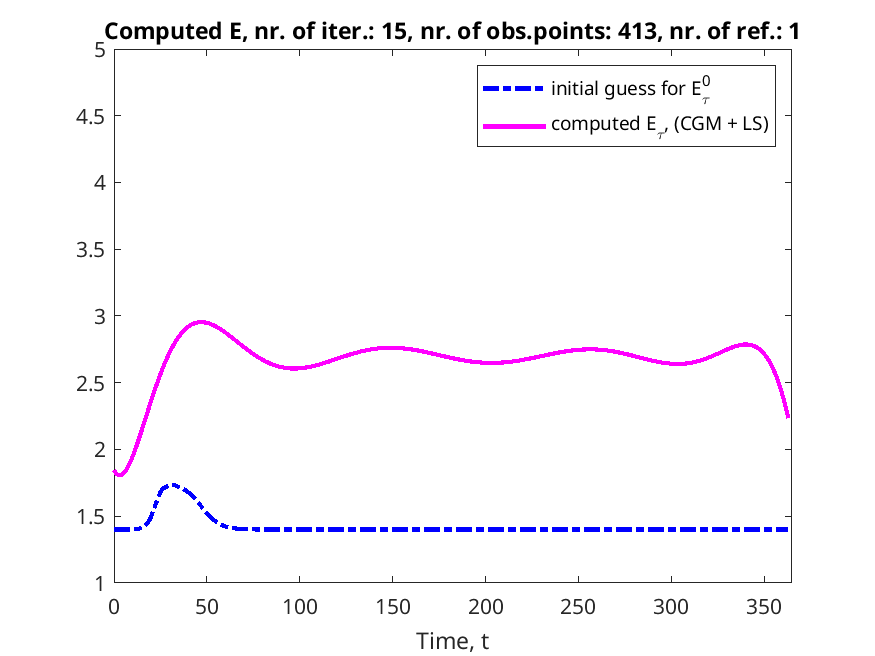}} &
    {\includegraphics[scale=0.35, clip = true, trim = 0.0cm 0.0cm 0.0cm 0.0cm]{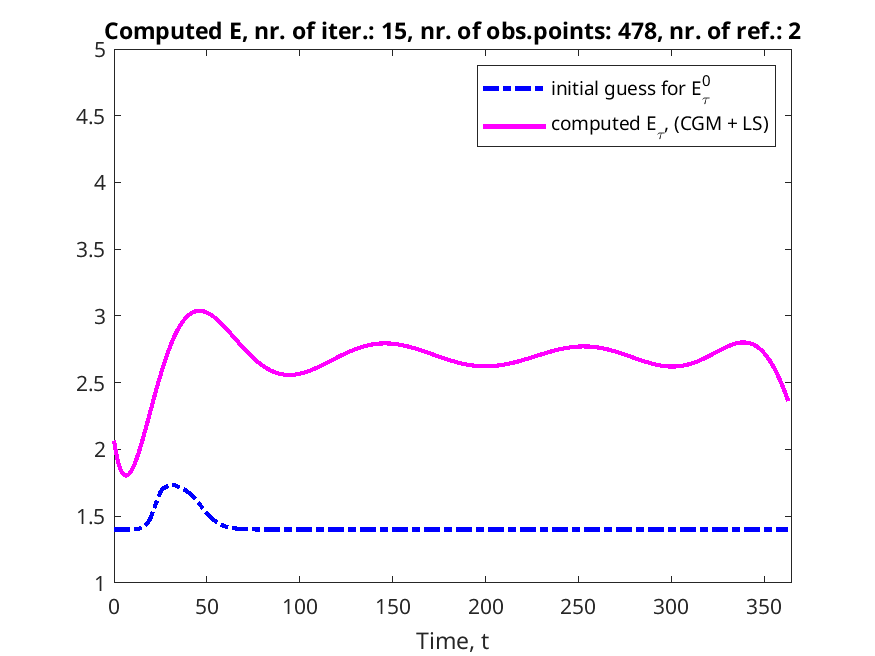}} \\
a)  $k=0$  & b) $k=1$ & c) $k=2$  \\
\end{tabular}
\end{center}
\caption{Patient 4:  results of reconstruction of the function $E(t)$ on   $k$ times adaptivelly refined meshes $J_\tau^k$
 in ACGA algorithm. Computations are  performed for the mesh refinement parameter $\tilde{\beta_k}
  =0.875$  for all mesh refinements $k$.}
 \label{fig:Pat4RecE}
 \end{figure}

Figures \ref{fig:Pat4CTLfunc}  show  time-dependent  functions for initial guess $E^0(t)$ as well as $ f(t), d(t)$
 with CTL response  for parameters defined in \eqref{paramCTLPat1} which are the same as for the patient 4.
These functions are used  now in optimization algorithms  for reconstruction of the immune response function $E(t)$.
 Figures \ref{fig:Pat4RecE} show
 computed functions $E_\tau^k$
 on   $k=0,1,2$ times adaptivelly refined meshes $J_\tau^k$.

\begin{figure}
\begin{center}
\begin{tabular}{ccc}
  {\includegraphics[scale=0.35, clip = true, trim = 0.0cm 0.0cm 0.0cm 0.0cm ]{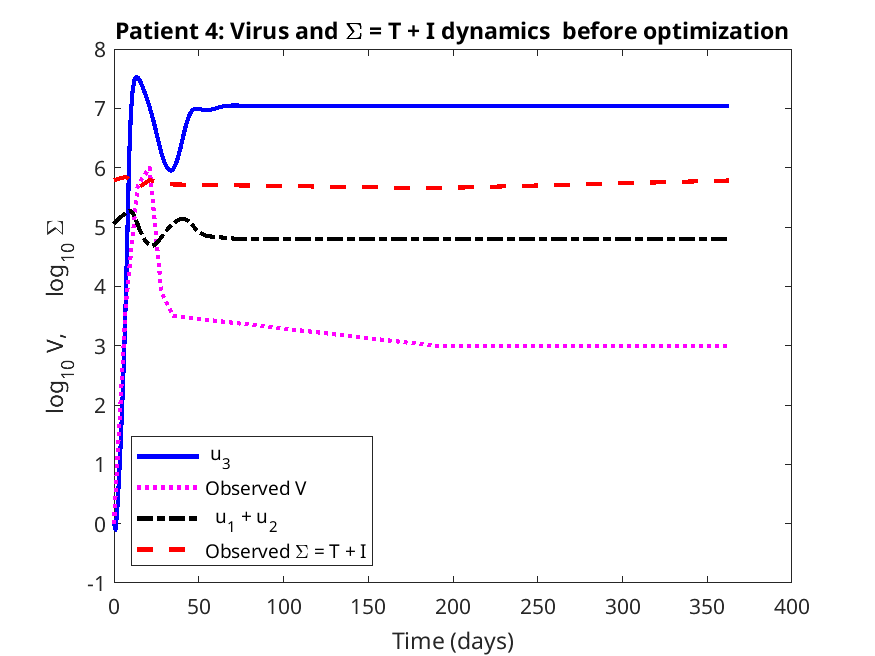}} &
  {\includegraphics[scale=0.35, clip = true, trim = 0.0cm 0.0cm 0.0cm 0.0cm]{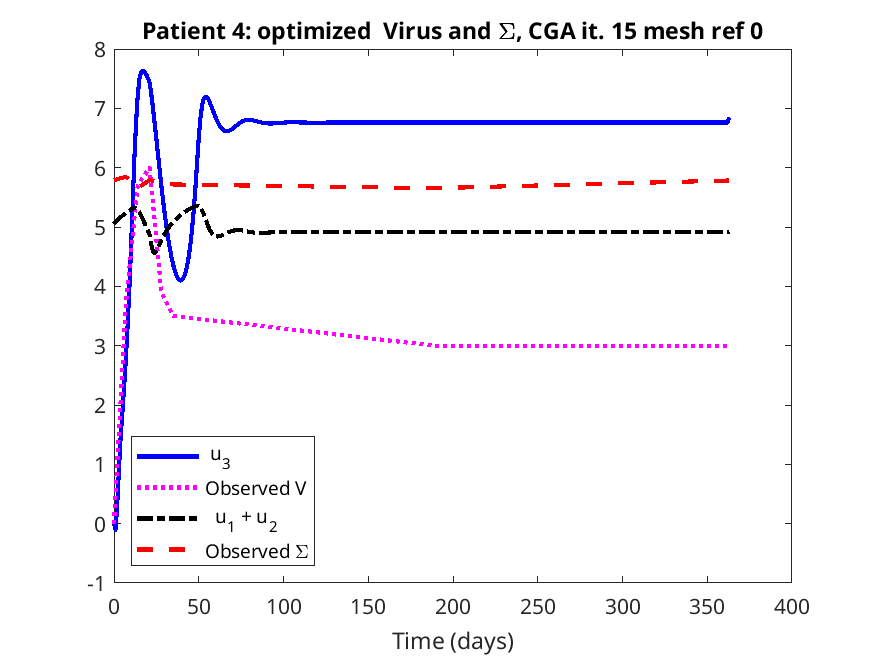}} &
   {\includegraphics[scale=0.35, clip = true, trim = 0.0cm 0.0cm 0.0cm 0.0cm]{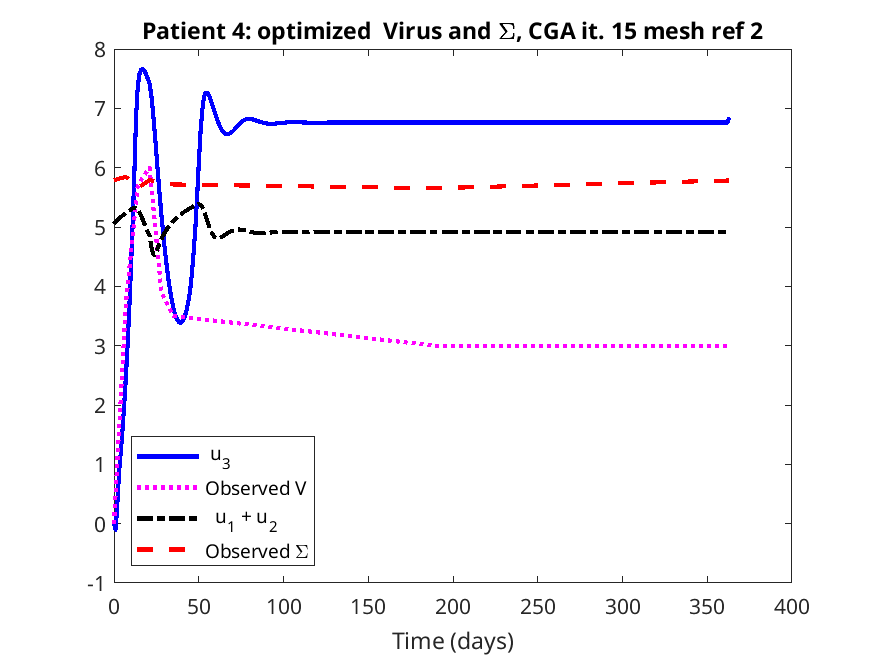}} \\
 a) $k=0$ & b) $k=0$  & c) $k=2$ \\
\end{tabular}
\end{center}
\caption{ Patient 4: Dynamics of the computed virus function $V_\tau^k$ before and after optimization corresponding to the
   computed $E_\tau^k$ on   $k, k=0,2$ times adaptivelly refined meshes versus interpolated clinical data ${g_1^0}_\tau,{g_2^0}_\tau$.}
 \label{fig:Pat4virus}
 \end{figure}

\begin{figure}
\begin{center}
\begin{tabular}{cc}
  {\includegraphics[scale=0.5, clip = true, trim = 0.0cm 0.0cm 0.0cm 0.0cm ]{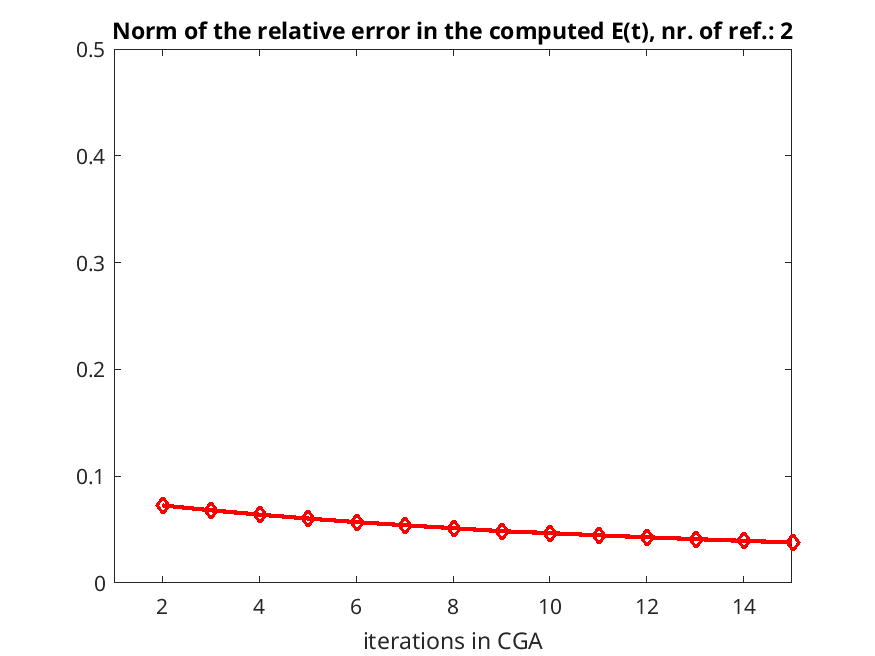}} &
  {\includegraphics[scale=0.5, clip = true, trim = 0.0cm 0.0cm 0.0cm 0.0cm]{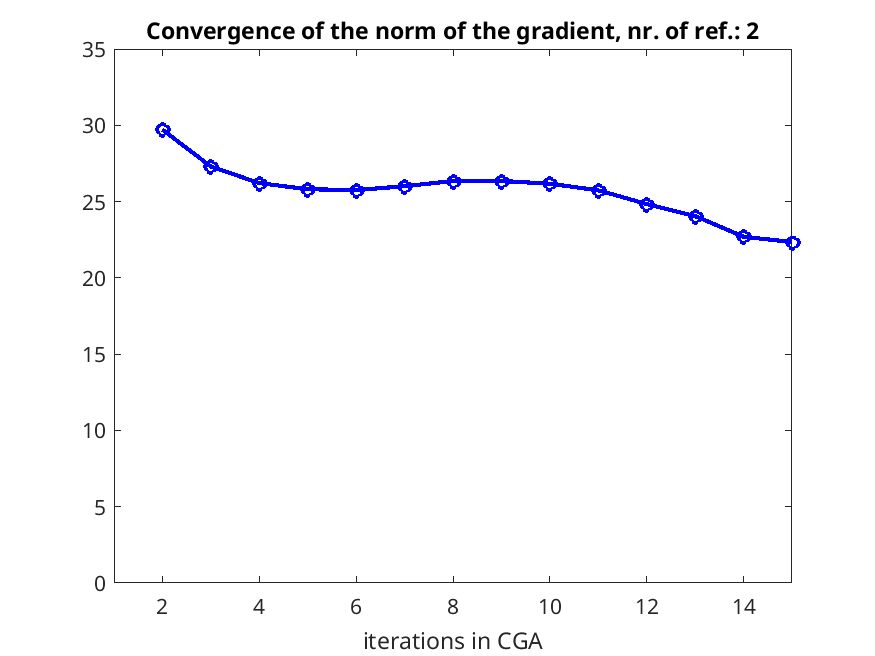}} \\
  $\frac{  \| E_\tau^{m} - E_\tau^{m-1} \|_{L_2(\Omega_t)}}{\| E_\tau^m\|_{L_2(\Omega_t)}}$   &  $\|G^m(t) \|_{L_2(\Omega_t)} $ 
\end{tabular}
\end{center}
\caption{Patient 4:   computed relative norms  $\frac{  \| E_\tau^{m} - E_\tau^{m-1} \|_{L_2(\Omega_t)}}{\| E_\tau^m\|_{L_2(\Omega_t)}}$
and $\|G^m(t) \|_{L_2(\Omega_t)}$  on the mesh $J_\tau^2$.}
 \label{fig:Pat4Residuals}
 \end{figure}

\begin{figure}
\begin{center}
\begin{tabular}{cc}
  {\includegraphics[scale=0.5, clip = true, trim = 0.0cm 0.0cm 0.0cm 0.0cm ]{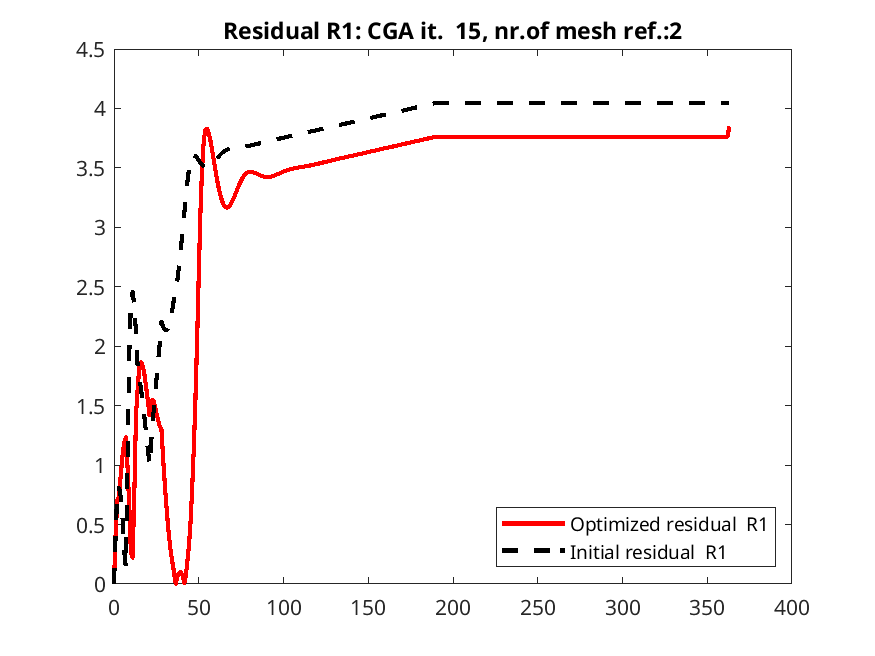}} &
  {\includegraphics[scale=0.5, clip = true, trim = 0.0cm 0.0cm 0.0cm 0.0cm]{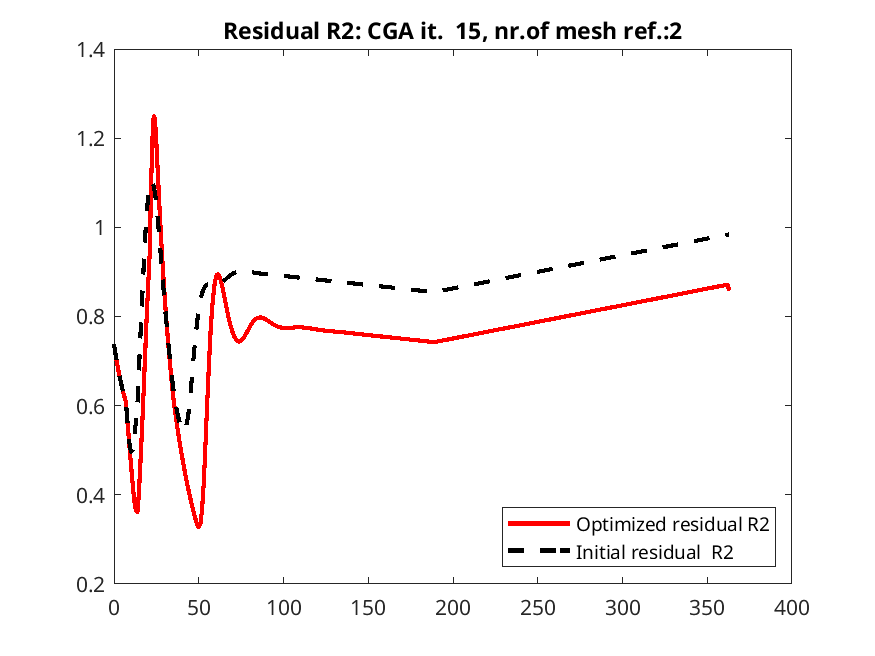}} \\
   a) $  R_1 $  & b) $ R_2  $  \\
 {\includegraphics[scale=0.5, clip = true, trim = 0.0cm 0.0cm 0.0cm 0.0cm]{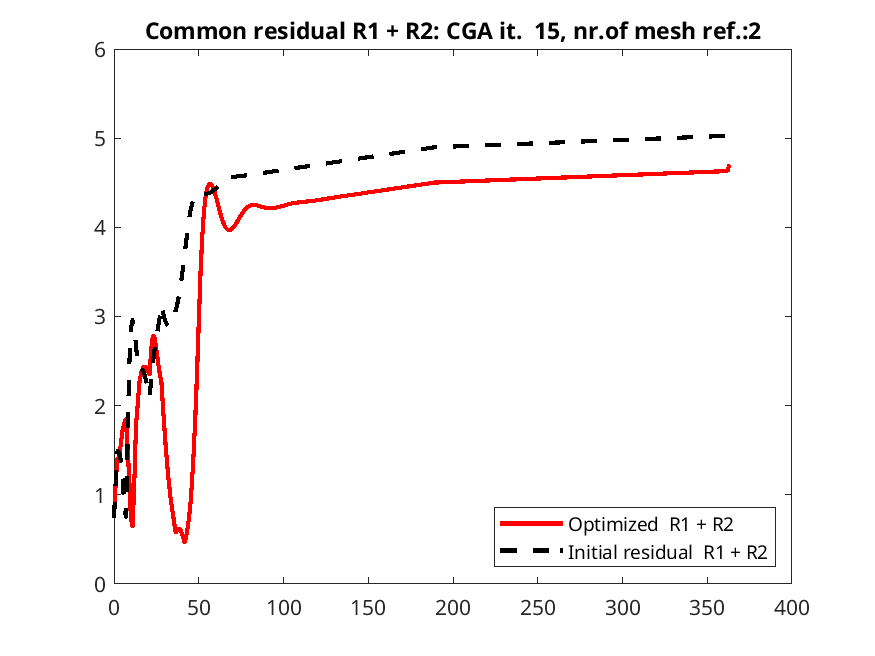}} &
  {\includegraphics[scale=0.5, clip = true, trim = 0.0cm 0.0cm 0.0cm 0.0cm]{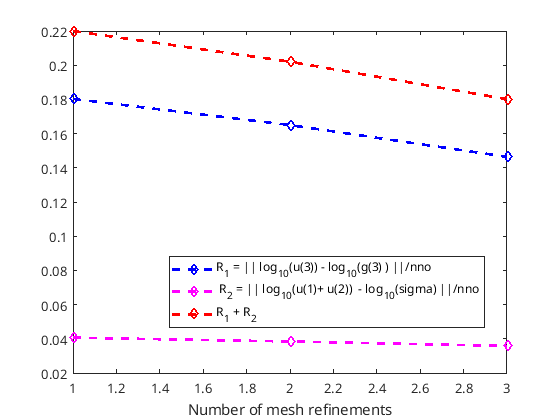}}  \\
 c)  $ R_1 + R_2  $ &  d) 
\end{tabular}
\end{center}
\caption{Patient 3: a), b), c)  computed residuals $R_1$ and $R_2$ on the mesh $J_\tau^2$; d) Comparison of related residuals
 $\|R_1\|, \|R_2\|$ on the meshes
$J_\tau^k, k = 0,1,2$. }
 \label{fig:Pat4R1R2}
 \end{figure}

Figures \ref{fig:Pat4virus}  show  time-dependent behaviour  of the computed virus function $V_\tau^k, k=0,2$ before and after applying ACGA.
Here,  the computed virus function $V_\tau^k$ corresponds to the
   computed $E_\tau^k$ presented in Figure \ref{fig:Pat4RecE}.
Results are compared with  clinical data  for the observed virus function  $V = g_2^0$  and observed total number of  the uninfected
and infected cells $\Sigma = g_1^0$.
 Computed optimized virus function on $J_\tau^0$ is presented in the
Figure \ref{fig:Pat4virus}-b).  Figure \ref{fig:Pat4virus}-c) shows
computed virus function on the twice  locally refined time mesh
$J_\tau^2$.  Comparing figures  \ref{fig:Pat4virus}-b) and  \ref{fig:Pat4virus}-c)
we  again observe that local adaptive mesh refinement allows more exact fit data to achieve minimum of the virus load.
Figures   \ref{fig:Pat4Residuals} show
stabilization of the computed relative norms  
$\frac{  \| E_\tau^{m} - E_\tau^{m-1} \|_{L_2(\Omega_t)}}{\| E_\tau^m\|_{L_2(\Omega_t)}}$ 
  and behaviour of norms of the computed gradient $\|G^m(t) \|_{L_2(\Omega_t)} $  for patient 4.
 Relative residuals $R_1$ and $R_2$ and $\|R_1\|_{L_2(\Omega_t)} , \|R_2\|_{L_2(\Omega_t)} $ are presented in the Figures
 \ref{fig:Pat4R1R2} as well as in the Table 4.


\section{Conclusion}

We propose a time-adaptive optimization framework for determining the
time-dependent immune response function within a mathematical model of
acute HIV infection, using clinical data from four untreated patients.

This approach begins with estimating the immune response on a coarse
initial time grid, based on several known values of the observable
functions. The time grid is then locally refined at points where the
residual $|R(E_\tau)(t)|$ reaches its maximal values. The immune
response function is subsequently recomputed on the refined time mesh,
enhancing the resolution where it is most needed.

Our methodology uses Lagrangian approach, from which we derive the
optimality conditions and a numerical scheme to solve the forward
problem, adjoint problem, and parameter identification
problem. Furthermore, we establish three distinct a posteriori error
estimates and formulate an adaptive optimization algorithm tailored to
this dynamic framework.

Numerical experiments demonstrate the effectiveness of the proposed
adaptive method in reconstructing the immune response during the acute
phase of HIV infection, using patient-specific clinical data. These
results highlight that local adaptive mesh refinement yields a more
precise data fit to achieve
 minimum of the virus function at the acute phase of HIV infection.

However, further computational investigations are required to
accurately fit clinical data over the clinical latency stage of HIV,
using the same dataset.

The proposed
 time-adaptive optimization strategy
 has the potential to aid clinicians by enabling individualized parameter identification in specific PIP.


\end{document}